\documentclass[11pt]{article}

\usepackage{latexsym}
\usepackage{amsfonts}
\usepackage{amssymb}
\usepackage{amsmath}
\usepackage{mathrsfs}
\input xypic
\usepackage{hyperref}

\newcommand{\bldr}{\mathbf{r}}
\newcommand{\blds}{\mathbf{s}}

\newcommand{\bldm}{\mathbf{m}}
\newcommand{\bldn}{\mathbf{n}}
\newcommand{\bldz}{\mathbf{z}}

\newcommand{\sgn}{\mbox{sign\,}}

\newcommand{\be}{\begin{equation}}
\newcommand{\ee}{\end{equation}}
\newcommand{\bea}{\begin{eqnarray}}
\newcommand{\eea}{\end{eqnarray}}
\newcommand{\bean}{\begin{eqnarray*}}
\newcommand{\eean}{\end{eqnarray*}}
\newcommand{\brray}{\begin{array}}
\newcommand{\erray}{\end{array}}

\newtheorem{dfn}{Definition}[section]
\newtheorem{thm}[dfn]{Theorem}
\newtheorem{lmma}[dfn]{Lemma}
\newtheorem{ppsn}[dfn]{Proposition}
\newtheorem{crlre}[dfn]{Corollary}
\newtheorem{xmpl}[dfn]{Example}
\newtheorem{rmrk}[dfn]{Remark}

\newcommand{\bdfn}{\begin{dfn}\rm}
\newcommand{\bthm}{\begin{thm}}
\newcommand{\blmma}{\begin{lmma}}
\newcommand{\bppsn}{\begin{ppsn}}
\newcommand{\bcrlre}{\begin{crlre}}
\newcommand{\bxmpl}{\begin{xmpl}}
\newcommand{\brmrk}{\begin{rmrk}\rm}

\newcommand{\edfn}{\end{dfn}}
\newcommand{\ethm}{\end{thm}}
\newcommand{\elmma}{\end{lmma}}
\newcommand{\eppsn}{\end{ppsn}}
\newcommand{\ecrlre}{\end{crlre}}
\newcommand{\exmpl}{\end{xmpl}}
\newcommand{\ermrk}{\end{rmrk}}

\newcommand{\bbc}{\mathbb{C}}
\newcommand{\bbz}{\mathbb{Z}}
\newcommand{\bbm}{\mathbb{M}}
\newcommand{\bbn}{\mathbb{N}}
\newcommand{\bbr}{\mathbb{R}}

\newcommand{\bbt}{\mathbb{T}}

\newcommand{\scrb}{\mathscr{B}}

\newcommand{\scrg}{\mathscr{G}}

\newcommand{\scrt}{\mathscr{T}}

\newcommand{\cla}{\mathcal{A}}
\newcommand{\clb}{\mathcal{B}}
\newcommand{\cle}{\mathcal{E}}
\newcommand{\clf}{\mathcal{F}}
\newcommand{\clh}{\mathcal{H}}

\newcommand{\clk}{\mathcal{K}}
\newcommand{\cll}{\mathcal{L}}
\newcommand{\clp}{\mathcal{P}}

\newcommand{\cls}{\mathcal{S}}

\newcommand{\prf}{\noindent{\it Proof\/}: }
\newcommand{\seq}{\subseteq}
\newcommand{\one}{{1\!\!1}}

\newcommand{\id}{\mbox{id}}

\def \qed { \mbox{}\hfill
$\Box$\vspace{1ex}}

\newcommand{\half}{\frac{1}{2}}

\newcommand{\trace}{\mbox{\textit{Trace}\,}}
\begin{document}

%%%%%%%%%%%%%%%%%%%%%%%%%%%%%%%%%

%%%%%%%%%%%%%%%%%%%%%%%%%%%%%%%%%
\author{{\sc Arupkumar Pal} and
{\sc S. Sundar}}
\title{Regularity and dimension spectrum
of the equivariant spectral triple
for the odd dimensional quantum spheres}
\maketitle
%%%%%%%%%%%%%%%%%%%%%%%%%%%%%%%%%%
%%%%%  ABSTRACT
%%%%%%%%%%%%%%%%%%%%%%%%%%%%%%%%%%
 \begin{abstract}
The odd dimensional quantum sphere $S_q^{2\ell+1}$ is a homogeneous
space for the quantum group $SU_q(\ell+1)$. A generic equivariant spectral
triple for $S_q^{2\ell+1}$ on its $L_2$ space was constructed by
Chakraborty \& Pal in \cite{cha-pal-2008a}. We prove
regularity for that spectral triple here. We also compute its
dimension spectrum and show that it is simple.
We give detailed construction of its smooth function
algebra and some related algebras that help proving
regularity and in the computation of the dimension spectrum.
Following the idea of Connes for $SU_q(2)$, 
we first study another spectral triple for $S_q^{2\ell+1}$
equivariant under torus group action
constructed by Chakraborty \& Pal in \cite{cha-pal-2007a}.
We then derive the results for the $SU_q(\ell+1)$-equivariant
triple in the $q=0$ case from those for the torus equivariant
triple. For the $q\neq 0$ case, we deduce regularity and dimension
spectrum from the $q=0$ case.
 \end{abstract}
%%%%%%%%%%%%%%%%%%%%%%%%%%%%%%%%%%
{\bf AMS Subject Classification No.:} {\large 58}B{\large 34}, {\large
46}L{\large 87}, {\large
  19}K{\large 33}

%%%%%%%%%%%%%%%%%%%%%%%%%%%%%%%%%%%%%%%%%%%%%%%%%%%

\tableofcontents

%%%%%%%%%%%%%%%%%%%%%%%%%%%%%%%%%%%%%%%%%%%%%%%%%%%%%%%%%%%%%%%%%%%
\section{Introduction}
%%%%%%%%%%%%%%%%%%%%%%%%%%%%%%%%%%%%%%%%%%%%%%%%%%%%%%%%%%%%%%%%%%%
%%%%%%%%%%%%%%%%%%%%%%%%%%%%%%%%%%%%%%%%%%%%%%%%%%%%
In noncommutative geometry, the starting point is usually
a separable unital $C^*$-algebra $A$
which is the noncommutative version of a compact Hausdorff
space. Associated to this, one  has certain invariants like 
the $K$-groups and the $K$-homology groups.
In geometry, what one does next is to equip the topological space with a
smooth structure so that in particular one can then talk about its de Rham
cohomology. In the noncommutative situation, the parallel
is to look for an appropriate dense subalgebra of $A$ that will
play the role of smooth functions on the space. Given this
dense subalgebra, one can  compute various cohomology groups
associated with it, namely the Hochschild cohomology, cyclic cohomology
and the periodic cyclic cohomology, which are noncommutative and
far-reaching generalizations of ordinary de Rham homology and cohomology.
Question is: how to get hold of this dense subalgebra? 
One answer to this lies in the notion of a spectral triple,
which plays a central role in Connes' formulation of noncommutative geometry.

In ordinary differential geometry, with just a smooth structure
on a manifold, one can hardly go very far.
In order that one can talk about shapes and sizes of spaces,
one needs to bring in extra structure. One example is the Riemannian
structure, which gives rise to a Riemannian connection,
which in turn enables one to talk about curvature and so on.
Other examples of such extra structures are $Spin$ and $Spin^c$ structures.
In the presence of these extra structures, one has an operator 
theoretic data that completely encodes the geometry.
In noncommutative geometry, one takes this operator theoretic data
as the initial data and this is what goes by the name \textbf{spectral triple}.

\bdfn
Let $\cla$ be an associative unital  *-algebra.
An \textbf{even spectral triple} for $\cla$ is a triple $(\clh,\pi,D)$  together
with a $\bbz_2$-grading $\gamma$ on $\clh$ such that
\begin{enumerate}
 \item
a (complex separable) Hilbert space $\clh$,
\item
a *-representation $\pi:\cla\longrightarrow \cll(\clh)$ (usually assumed faithful),
\item
a self-adjoint operator $D$ with compact resolvent such that
$[D,\pi(a)]\in\cll(\clh)$ for all $a\in \cla$,
\item
$\pi(a)\gamma=\gamma\pi(a)$ for all $a\in \cla$ and
$D\gamma=-\gamma D$.
\end{enumerate}
If no grading is present, one calls it an \textbf{odd spectral triple}.
\edfn
The algebra $\cla$ appearing in 
this definition is in general  different from the $C^*$-algebra $A$
one starts with. Typically it is a dense subalgebra
in $A$, big enough so that the $K$-groups of $A$ and $\cla$
coincide.

Since $D$ has compact resolvent, it has finite dimensional kernel. 
Hence by making a finite rank perturbation, one can
make $D$ invertible. Now if one replaces $D$ with $\sgn D$,
then one gets the notion of a Fredholm module over the algebra $\cla$.
If $\cla$ is a dense *-subalgebra of a $C^*$-algebra $A$,
then this Fredholm module extends uniquely and gives a Fredholm
module over $A$. In other words, one gets an element in the 
$K$-homology group of $A$.
This gives a map from the $K$-theory of $A$ to the set of integers
via the $K$-theory--$K$-homology pairing.
Starting from this Fredholm module, one can construct its Chern character,
which gives an element in the periodic cyclic cohomology
of $\cla$, which, in turn, gives a map from the $K$-theory of $A$ to $\bbc$
via the periodic cyclic cohomology--$K$-theory pairing.
The two maps thus obtained are the same. This is the content
of the index theorem.

The Chern character is often difficult to compute. And that is where
the spectral triple comes into the picture. Under certain hypothesis
on the spectral triple, one can construct a cyclic cocycle,
i.e.\ an element in the periodic cyclic cohomology 
that differs from the Chern character by a coboundary, so that
it gives rise to the same map from $K$-theory to $\bbz$.
Under mild hypothesis on an invariant known as the dimension spectrum,
these cocycles are given in terms of certain residue functionals,
that can be relatively easier to compute.
This is Connes--Moscovici local index theorem (\cite{con-mos-1995a})
which is one of the major results in noncommutative geometry.
This was first proved in the context of transverse geometry
of foliations, but is much more general in nature and has
wider applicability. This was illustrated by Connes
in \cite{con-2004a}, where he made a detailed analysis
of the equivariant spectral triple for the quantum
$SU(2)$ group constructed in \cite{cha-pal-2003a}.
A similar analysis was later done by Dabrowski et al in
\cite{v-d-l-s-v-2005a} for the spectral triple
constructed in \cite{d-l-s-v-v-2005a}.
Typically, the $C^*$-algebras associated to quantum groups
or their homogeneous spaces are given by a set of generators
and relations. While constructing spectral triples, one does
it for the associated coordinate function algebra, i.e.\ the
*-subalgebra generated by these generators. This algebra is
 not closed under the holomorphic function calculus
of the $C^*$-algebra. Therefore one needs to construct the smooth
function algebra, prove regularity and compute the dimension spectrum
in order to be able to apply the Connes--Moscovici theorem.
This was done in \cite{con-2004a}, where he also defined
a symbol map and formulae for computing the residue functionals
in terms of the symbol maps were given.

Odd dimensional quantum spheres are higher dimensional
analogues of the quantum $SU(2)$. The $(2\ell+1)$ dimensional sphere
$S_q^{2\ell+1}$ is a homogeneous space of the quantum group
$SU_q(\ell+1)$. In~\cite{cha-pal-2008a}, Chakraborty and Pal
constructed a generic spectral triple on the $L_2$ space of the sphere
with nontrivial K-homology class and equivariant under the action
of $SU_q(\ell+1)$.
The main aim of the present article is to prove that this spectral triple is
regular. We also introduce the smooth function algebra and 
compute the dimension spectrum. The dimension spectrum is shown
to be simple so that the
Connes--Moscovici local index theorem is applicable to
this triple. The local index computation will be taken up
in a separate article.

Here is a brief outline of the contents of this paper.
In the next section, we recall a few basic notions
from \cite{con-mos-1995a}. We then collect together
a few observations and remarks on tensor products of Fr\'echet algebras
and fix some of the notations.
In section~3, we first look at the torus equivariant spectral triple
for the  spheres and introduce a smooth function algebra,
prove regularity and compute the dimension spectrum.
The key here is the short exact sequence
\[
  0\longrightarrow \clk\otimes C(\bbt)\longrightarrow C(S_q^{2\ell+1})
 \longrightarrow C(S_q^{2\ell-1}) \longrightarrow 0,
\]
and  results by Schweitzer (\cite{sch-1992a}, \cite{sch-1993a}) on spectral invariance.
Using these and using the idea employed by Connes in \cite{con-2004a},
we build the smooth function algebra over $S_q^{2\ell+1}$
recursively starting from $C^\infty(\bbt)$.

In section~4, we deal with the $SU_q(\ell+1)$-equivariant spectral triple.
We first treat the $q=0$ case. Using a decomposition of the $L_2$-space, 
we relate it to the torus equivariant
triple with a certain multiplicity. Regularity, smooth function
algebra and properties of the dimension spectrum
all then follow from the results in section~3.
Here again the idea is exactly as in \cite{con-2004a} for $SU_q(2)$. 
In subsection~4.4 we take up the $q\neq 0$ case.
We take a close look at the representation
of the algebra, in particular the images of the generating elements
and after a careful analysis, we prove that modulo
operators that one can neglect for the purpose of
computing the dimension spectrum, things can be deduced
from the torus equivariant case again.

%%%%%%%%%%%%%%%%%%%%%%%%%%%%%%%%%%%%%%%%%%%%%%%%%%%%%%%%%%%%%%%%%%%
\section{Preliminaries}
%%%%%%%%%%%%%%%%%%%%%%%%%%%%%%%%%%%%%%%%%%%%%%%%%%%%%%%%%%%%%%%%%%%
%%%%%%%%%%%%%%%%%%%%%%%%%%%%%%%%%%%%%%%%%%%%%%%%%%%
\subsection{Regular spectral triples}
%%%%%%%%%%%%%%%%%%%%%%%%%%%%%%%%%%%%%%%%%%%%%%%%%%%
In this subsection we recall some definitions and notions from
\cite{con-mos-1995a}.
Let $D$ be a selfadjoint
operator on a Hilbert space $\clh$ which is invertible. Define $\clh_{s}=Dom(| D
|^{s})$ for $s \geq 0$. Then $\clh_{s}$ is a decreasing family of vector
subspaces of $\clh$. Let $\clh_{\infty}:= \cap_{s\geq 0}\clh_{s}$. The subspace
$\clh_{\infty}$ is a dense subspace of $\clh$.
\bdfn
 An operator $T:\clh_{\infty} \to \clh_{\infty}$ is said to be smoothing if for every
$m,n \geq 0$ the operator $| D |^{m}T | D |^{n}$ is bounded. The
vector space of smoothing operators is denoted by $OP^{-\infty}$.
\edfn
For $T \in OP^{-\infty}$, define $\| T \|_{m,n} = \| | D
|^{m}T| D |^{n} \|$ for $m,n \geq 0$.
%%%%%%%%%%%%%%%%%%%%%%%%%%%%%%%%%%%%%%%%%%
\blmma
 The vector space $OP^{-\infty}$ is an involutive subalgebra of $\cll(\clh)$ and equipped
with the family of seminorms $\|\cdot\|_{m,n}$, it is a Fr\'echet algebra.
\elmma
%%%%%%%%%%%%%%%%%%%%%%%%%%%%%%%%%%%%%%%%%%
Let $\delta$ be the unbounded derivation $[| D |,\cdot]$. More precisely,
$Dom(\delta)$ consists of all bounded operators $T$ which leaves $Dom(| D
|)$ invariant and for which $\delta(T):=[| D |,T]$ extends to a bounded
operator.
%%%%%%%%%%%%%%%%%%%%%%%%%%%%%%%%%%%%%%%%%%
\blmma[\cite{con-2008a}]
The unbounded derivation $\delta$ is a closed derivation i.e. if $T_{n}$ is a
sequence in $ Dom(\delta)$ such that $T_{n} \to T$ and $\delta(T_{n}) \to S$
then $T \in Dom(\delta)$ and $\delta(T)=S$.
\elmma
%%%%%%%%%%%%%%%%%%%%%%%%%%%%%%%%%%%%%%%%%%

Define $OP^{0}:=\{T \in \cll(\clh):T \in \cap_{n} Dom(\delta^{n})\}$. The
following lemma says that elements of $OP^{0}$ are operators on $\clh_{\infty}$.
%%%%%%%%%%%%%%%%%%%%%%%%%%%%%%%%%%%%%%%%%%
\blmma
\label{OP}
Let $T$ be a bounded operator on $\clh$. Then the following are equivalent.
\begin{enumerate}
\item The operator $T \in OP^{0}$.
\item The operator $T$ leaves $\clh_{\infty}$ invariant and
$\delta^{n}(T):\clh_{\infty} \to \clh_{\infty}$ is bounded for every $n \in
\mathbb{N}$
\end{enumerate}
\elmma
%%%%%%%%%%%%%%%%%%%%%%%%%%%%%%%%%%%%%%%%%%

It is easy to see from lemma \ref{OP} that $OP^{0}$ is an algebra and that for
every $m \in \mathbb{Z}$ that $| D |^{-m}T | D |^{m}$ is bounded if
$T \in OP^{0}$. As a consequence it follows that $OP^{-\infty}$ is an ideal in
$OP^{0}$. Now we recall the notions of regularity and dimension spectrum for a
spectral triple.

%%%%%%%%%%%%%%%%%%%%%%%%%%%%%%%%%%%%%%%%%%
\bdfn
Let $(\mathcal{A},\clh,D)$ be a spectral triple. We say that $(\mathcal{A},\clh,D)$ is
\textbf{regular} if $\mathcal{A}+[D,\mathcal{A}] \subset OP^{0}$.
\edfn
%%%%%%%%%%%%%%%%%%%%%%%%%%%%%%%%%%%%%%%%%%
A spectral triple $(\mathcal{A},\clh,D)$ is $p+$~summable if $| D |^{-p}$ is in 
the ideal of Dixmier traceable operators
$\cll^{(1,\infty)}$.
In particular, if $(\mathcal{A},\clh,D)$ is $p+$~summable, then
$| D |^{-s}$ is trace class for $s>p$. 
Let $(\mathcal{A},\clh,D)$ be a regular spectral triple which is $p+$~summable for
some $p$. Let $\mathcal{B}$ be the algebra generated by
$\delta^{n}(\mathcal{A})$ and $\delta^{n}([D,\mathcal{A}])$. We say that the
spectral triple $(\mathcal{A},\clh,D)$ has \textbf{discrete dimension spectrum}
$\Sigma \subset \mathbb{C}$ if $\Sigma$ is discrete and for every $b \in
\mathcal{B}$, the function $Trace(b| D |^{-z})$ initially defined for
$Re(z)>p$ extends to a meromorphic function with poles only in $\Sigma$. We say
the dimension spectrum is simple if all the poles are simple.

%%%%%%%%%%%%%%%%%%%%%%%%%%%%%%%%%%%%%%%
\subsection{Topological tensor products}
%%%%%%%%%%%%%%%%%%%%%%%%%%%%%%%%%%%%%%%
The $C^*$-algebras involved in tensor products
that we deal with in this paper are all nuclear.
Therefore no ambiguities arise due to nonuniqueness
of tensor products.
Apart from $C^*$-algebras and their tensor products,
we will also deal with Fr\'echet algebras and their
tensor products. Suppose $A_1$ and $A_2$ are two
Fr\'echet algebras with topologies coming from the families of seminorms
$(\|\cdot\|_{\lambda})_{\lambda\in\Lambda}$
and $(\|\cdot\|_{\lambda'})_{\lambda'\in\Lambda'}$.
For each pair $(\lambda,\lambda')\in\Lambda\times\Lambda'$,
one forms the projective cross norm
$\|\cdot\|_{\lambda,\lambda'}$ which is a seminorm
on the algebraic tensor product $A_1\otimes_{alg}A_2$.
The family 
$(\|\cdot\|_{\lambda,\lambda'})_{(\lambda,\lambda')\in\Lambda\times\Lambda'}$
then gives rise to a  topology on
$A_1\otimes_{alg}A_2$. Completion with respect to this
is a Fr\'echet algebra and
is called the projective tensor product of $A_1$ and $A_2$.
While talking about tensor product of two Fr\'echet algebras,
we will always mean their projective tensor product
and will denote it by $A_1\otimes A_2$.

We will mainly be concerned with Fr\'echet algebras
sitting inside some $\cll(\clh)$ with Fr\'echet topology
finer than the norm topology. In other words, we will be
dealing with Fr\'echet algebras with faithful continuous
representations on Hilbert spaces.
Let $A_{1},A_{2}$ be Fr\'echet algebras. If $\rho_{i}:A_{i} \to \cll(\clh_{i})$ are
continuous representations for $i=0,1$ where $\clh_{i}$'s are Hilbert
spaces, then by the universality of the projective tensor product it follows
that there exists a unique continuous representation $\rho_{1} \otimes
\rho_{2}:A_{1} {\otimes} A_{2} \to \cll(\clh_{1}\otimes \clh_{2})$ such that
$(\rho_{1} \otimes \rho_{2})(a_{1}\otimes a_{2})=\rho_{1}(a_{1})\otimes
\rho_{2}(a_{2})$. If $A_{i}$'s are subalgebras of $\cll(\clh_{i})$ then we will call
the tensor product representation of $A_{1} \otimes A_{2}$ on
$\clh_{1}\otimes \clh_{2}$ as the natural representation.

%%%%%%%%%%%%%%%%%%%%%%%%%%%%%%%%%%%%%%%%%%%%%%%%%%%%%%%
\blmma
\label{ProductofSpectralTriples}
Let $(A_{1},\clh_{1},D_{1})$ and $(A_{2},\clh_{2},D_{2})$ be regular spectral triples.
Assume that the following conditions hold
\begin{enumerate}
\item The algebras $A_{1}$ and $A_{2}$ are Fr\'echet algebras represented
faithfully on $\clh_{1}$ and $\clh_{2}$ respectively.
\item The selfadjoint operators $D_{1}$ and $D_{2}$ are positive with compact
resolvent.
\item For $i=0,1$, the unbounded derivations $\delta_{i}=[D_{i},.]$  leave
$A_{i}$ invariant and $\delta_{i}:A_{i}\to A_{i}$ is continuous.
\end{enumerate} 
Let $D:=D_{1}\otimes 1 + 1 \otimes D_{2}$. Suppose that the natural
representation of $A_{1}{\otimes}A_{2}$ on $\clh_{1}\otimes \clh_{2}$ is
faithful. Then the triple $(A_{1}{\otimes} A_{2},\clh_{1}\otimes \clh_{2},D)$ is a regular
spectral triple. More precisely the unbounded derivation $\delta:=[D,.]$ leaves
the algebra $A_{1}{\otimes} A_{2}$ invariant and the map
$\delta:A_{1}{\otimes}A_{2}\to A_{1}{\otimes}A_{2}$ is
continuous.
\elmma
%%%%%%%%%%%%%%%%%%%%%%%%%%%%%%%%%%%%%%%%%%%%%%%%%%%%%%%
\prf
Let $\delta'=\delta_{1}\otimes 1 + 1\otimes \delta_{2}$. Then
$\delta'$ is a continuous linear operator on $A_{1}{\otimes}A_{2}$.
Clearly $A_{1}\otimes_{alg}A_{2} \subset Dom(\delta)$ and $\delta=\delta'$ on
$A_{1}\otimes_{alg}A_{2}$. Now let $a \in A_{1}{\otimes} A_{2}$ be
given. Choose a sequence $(a_{n}) \in A_{1}\otimes_{alg}A_{2}$ such that
$a_{n}\to a$ in $A_{1}{\otimes}A_{2}$. Then $a_{n} \to a$ in
$\cll(\clh_{1}\otimes \clh_{2})$. Since $\delta'$ is continuous and because the
inclusion $A_{1}{\otimes}A_{2} \subset \cll(\clh_{1}\otimes \clh_{2})$ is continuous,
it follows that the sequence $(\delta'(a_{n}))=(\delta(a_{n}))$ is Cauchy in
$\cll(\clh_{1}\otimes \clh_{2})$. As $\delta$ is closed, it follows that $a \in
Dom(\delta)$ and $\delta(a)=\delta'(a)$. Now the proposition follows.
\qed

The above lemma can  be extended to tensor product of finite number of
spectral triples with the appropriate assumptions.

%%%%%%%%%%%%%%%%%%%%%%%%%%%%%%%%%%%%%%%%%%%%%%%%%%%%%%
\paragraph{Notations.}
Let us now collect together some of the notations that will
be used throughout the paper.
The symbol $\clh$, with or without subscripts, will
denote a Hilbert space. The space of bounded
linear operators on $\clh$ will be denoted by $\cll(\clh)$
and the space of compact linear operators on $\clh$ will be
denoted by $\clk(\clh)$.
We will denote by
$\Sigma$ the set $\{1,2,\ldots,2\ell+1\}$ and
by $\Sigma_{\ell}$ and $\Sigma_{j,\ell}$ the subsets $\{1,2,\ldots,\ell+1\}$
and $\{\ell-j+1,\ell-j+2,\ldots,\ell+1\}$ respectively, where  $0\leq j\leq\ell$.

Let $\Gamma\equiv\Gamma_\Sigma$ denote the set of maps $\gamma$ from $\Sigma$
to $\bbz$ such that $\gamma_{i}\in\bbn$ for all $i\in\Sigma\backslash\{\ell+1\}$,
i.e.\ $\Gamma_\Sigma=\bbn^\ell\times \bbz\times \bbn^\ell$.
For a subset $A$ of $\Sigma$, we will denote by 
$\gamma_A$ the restriction $\gamma|_A$ of $\gamma$ to $A$.
Let $\Gamma_A$ denote the
set $\{\gamma_A: \gamma\in\Gamma\}$ and $\clh_A$ be
the Hilbert space $\ell_2(\Gamma_A)$.
We will denote $\clh_\Sigma$ by just $\clh$,
and $\clh_{\Sigma_{j,\ell}}$ by $\clh_j$.
Thus 
\[
\clh_\Sigma= \underbrace{\ell_2(\bbn)\otimes\cdots\otimes
    \ell_2(\bbn)}_{\ell \mbox{ copies}}\otimes \ell_2(\bbz)\otimes 
   \underbrace{\ell_2(\bbn)\otimes\cdots\otimes
         \ell_2(\bbn)}_{\ell \mbox{ copies}},\qquad
\clh_j=\underbrace{\ell_2(\bbn)\otimes\cdots\otimes
    \ell_2(\bbn)}_{j \mbox{ copies}}\otimes \ell_2(\bbz)
\]
Note that $\clh_j$ and $\clh_{\{j\}}$ are different.

Let $A\seq\Sigma$. We will denote by $\{e_\gamma\}_\gamma$
the natural orthonormal basis for $\clh_A=\ell_2(\Gamma_A)$
and by $p_\gamma$ the rank one projection 
$|e_\gamma\rangle\langle e_\gamma|$.
For $i\in A$, we will denote by $N_i$ the number operator
on the $i$\raisebox{.4ex}{th} coordinate on $\clh_A$, i.e.\ 
\[
N_i\equiv \sum_\gamma \gamma_i p_\gamma: e_\gamma\mapsto \gamma_i e_\gamma 
  (\text{defined on $\clh_A$ with $i\in A$}).
\]
We will denote by $|D_A|$ the operator $\sum_{i\in A}|N_i|$ on $\clh_A$.

Let $F_0$ be the following operator on $\ell_2(\bbz)$:
\[
 F_0 e_k=\begin{cases}
          e_k & \text{ if }k\geq 0,\cr
          -e_k & \text{ if }k<0.
         \end{cases}
\]
For $1\leq j\leq 2\ell+1$, let $V_j$ be the operator on $\clh_{\{j\}}$
defined by
\[
 V_j:=\begin{cases}
       F_0 & \text{ if } j=\ell+1,\cr
       I & \text{ otherwise.}
      \end{cases}
\]
Let $F_A$ denote the operator  $\otimes_{j\in A}V_j$ on $\clh_A$
and let $D_A=F_A|D_A|$.
Thus
\[
D_A e_\gamma= \begin{cases}
               -\left(\sum_{i\in A}|\gamma_i|\right) e_\gamma & 
                 \text{ if } \ell+1\in A \text{ and }\gamma_{\ell+1}<0,\cr
            \left(\sum_{i\in A}|\gamma_i|\right) e_\gamma & \text{ otherwise.}
              \end{cases}
\]
We will denote $F_{\Sigma_{j,\ell}}$ by $F_j$ and  $D_{\Sigma_{j,\ell}}$
by $D_j$.

Recall that $\clh_{\{j\}}$ is $\ell_2(\bbn)$ if $j\neq \ell+1$ and is
$\ell_2(\bbz)$ if $j=\ell+1$. Suppose for each $j\in\Sigma$,
$\clf_j$ is a subspace of $\cll(\clh_{\{j\}})$.
For $A\seq\Sigma$, define
\[
 \clf_{j,A}=\begin{cases}
             \clf_j & \text{ if $j\in A$,}\cr
             \bbc.I  & \text{ if $j\not\in A$,}
            \end{cases}
\]
and $\clf_A$ to be the tensor product $\otimes_{j\in\Sigma}\clf_{j,A}$
in $\cll(\clh_\Sigma)$
(the type of the tensor product will depend on the specific
$\clf_j$'s we look at). This tensor product will often be
identified with $\otimes_{j\in A}\clf_j\seq \cll(\clh_A)$.

On both $\ell_{2}(\mathbb{N})$ and
$\ell_{2}(\mathbb{Z})$, we will denote by $N$  the number operator defined by
$Ne_{n}=ne_{n}$ and by $S$  the left shift  defined by
$Se_{n}=e_{n-1}$. For $k\in \bbz$ (for $k\in\bbn$ in case of 
$\ell_2(\bbn)$), let $p_k$ denote the projection
$|e_k\rangle \langle e_k|$. We will freely identify
$\ell_2(\bbz)$ with $L_2(\bbt)$. Thus the right shift on $\ell_{2}(\mathbb{Z})$ 
will be multiplication by the function $t\mapsto t$ and
will be denoted by $\bldz$. Let $\scrt$ be the Toeplitz algebra, i.e.\ 
the $C^{*}$-subalgebra of $\cll(\ell_{2}(\mathbb{N}))$
generated by $S$. 
For a positive integer $k$, we will denote by $\scrt_k$ the $k$-fold tensor
product of $\scrt$, embedded in $\cll(\ell_2(\bbn^k))$.
Denote by $\sigma$ the symbol map from $\scrt$ to $C(\bbt)$ that
sends $S^*$ to $\bldz$ and all compact operators to 0.

%%%%%%%%%%%%%%%%%%%%%%%%%%%%%%%%%%%%%%%%%%%%%%%%%%%%%%%%%%%%%%%%%%%
\section{Torus equivariant spectral triple}
%%%%%%%%%%%%%%%%%%%%%%%%%%%%%%%%%%%%%%%%%%%%%%%%%%%%%%%%%%%%%%%%%%%

%%%%%%%%%%%%%%%%%%%%%%%%%%%%%%%%%%%%%%%%%%%%%%%%%%%%%%%%%%%%%%%%%%
\subsection{The spectral triple}
%%%%%%%%%%%%%%%%%%%%%%%%%%%%%%%%%%%%%%%%%%%%%%%%%%%%%%%%%%%%%%%%%%
 In this section we recall the spectral triple for the odd dimensional quantum
spheres given in \cite{cha-pal-2007a}.
We begin with some known facts about odd dimensional quantum spheres. 
Let $q\in[0,1]$.
The $C^*$-algebra $C(S_q^{2\ell+1})$ of the quantum
sphere $S_q^{2\ell+1}$
is the universal $C^*$-algebra generated by
elements
$z_1, z_2,\ldots, z_{\ell+1}$
satisfying the following relations (see~\cite{hon-szy-2002a}):
\bean
z_i z_j & =& qz_j z_i,\qquad 1\leq j<i\leq \ell+1,\\
z_i^* z_j & =& q z_j z_i^* ,\qquad 1\leq i\neq j\leq \ell+1,\\
z_i z_i^* - z_i^* z_i +
(1-q^{2})\sum_{k>i} z_k z_k^* &=& 0,\qquad \hspace{2em}1\leq i\leq \ell+1,\\
\sum_{i=1}^{\ell+1} z_i z_i^* &=& 1.
\eean
We will denote by $\cla(S_q^{2\ell+1})$ the *-subalgebra of $A_\ell$
generated by the $z_j$'s. Note that for $\ell=0$, the $C^*$-algebra
$C(S_q^{2\ell+1})$ is the algebra of continuous functions
$C(\bbt)$ on the torus and for $\ell=1$, it is $C(SU_q(2))$.

There is a natural torus group $\mathbb{T}^{\ell +1}$ action $\tau$ on
$C(S_{q}^{2\ell+1})$ as follows. For $w=(w_{1},\ldots,w_{\ell +1})$, define
an automorphism $\tau_{w}$ by $\tau_{w}(z_{i})=w_{i}z_{i}$.  
Let $Y_{k,q}$ be the following operators on $\clh_\ell$:
\be\label{eq:ykq}
 Y_{k,q}=\begin{cases}
 \underbrace{q^N\otimes\ldots\otimes q^N}_{k-1 \mbox{
copies}}\otimes
      \sqrt{1-q^{2N}}S^*\otimes 
   \underbrace{I \otimes\cdots\otimes I}_{\ell+1-k \mbox{ copies}}, & \mbox{ if } 1\leq k\leq \ell,\cr
   &\cr
    \underbrace{q^N\otimes\cdots\otimes q^N}_{\ell \mbox{ copies}}
       \otimes S^*, &  \mbox{ if } k=\ell+1.
         \end{cases}
\ee
Here for $q=0$, $q^N$ stands for the rank one projection $p_0=|e_0\rangle\langle e_0|$.
Then $\pi_\ell:z_k\mapsto Y_{k,q}$ gives a faithful representation
of $C(S_q^{2\ell+1})$ on $\clh_\ell$ for $q\in[0,1)$ 
(see lemma~4.1 and remark~4.5, \cite{hon-szy-2002a}).
We will denote the image $\pi_\ell(C(S_q^{2\ell+1}))$ by $A_\ell(q)$ or by just $A_\ell$.

 Let  $\{e_{\gamma}: \gamma \in \Gamma_{\Sigma_\ell} \}$ 
be the standard orthonormal basis for $\clh_{\ell }$. For
$w=(w_{1},w_{2},\cdots,w_{\ell +1})$ we define the unitary $U_{w}$ on $\clh_{\ell }$ by
$U_{w}(e_{\gamma})=w_{1}^{\gamma_{1}}w_{2}^{\gamma_{2}}\ldots
w_{\ell +1}^{\gamma_{\ell +1}}e_{\gamma}$ where
$\gamma=(\gamma_{1},\gamma_{2},\cdots,\gamma_{\ell +1})\in \Gamma_{\Sigma_\ell}$. 
Then $(\pi_{\ell },U)$ is a
covariant representation of $(C(S_{q}^{2\ell+1}),\mathbb{T}^{\ell +1},\tau)$.  Note that 
$A_{\ell } \subset \scrt_\ell\otimes C(\mathbb{T})$.

In \cite{cha-pal-2007a} all spectral triples equivariant with respect to
this covariant representation were characterised and an optimal one was
constructed. We recall the following theorem from \cite{cha-pal-2007a}.

%%%%%%%%%%%%%%%%%%%%%%%%%%%%%%%%%%%%%%%%%%%%%
\bthm[\cite{cha-pal-2007a}]  
Let $D_{\ell}$ be the operator $e_{\gamma} \to d(\gamma)e_{\gamma}$ on $\clh_{\ell }$
where the $d_{\gamma}$'s are given by
\begin{displaymath} 
\begin{array}{lll}
d(\gamma)&=&\left\{\begin{array}{ll}
                            \gamma_{1}+\gamma_{2}+\cdots
\gamma_{\ell }+|\gamma_{\ell +1}| & \text{ if  }\gamma_{\ell +1} \geq 0 \\
                            -(\gamma_{1}+\gamma_{2}+\cdots
\gamma_{\ell }+|\gamma_{\ell +1}|) & \text{ if }   \gamma_{\ell +1} < 0
                       \end{array} \right.
\end{array}
\end{displaymath}
Then $(\cla(S_{q}^{2\ell+1}),\clh_{\ell },D_{\ell})$ is a non-trivial $(\ell+1)$ summable
spectral triple. Also $D_{\ell}$ commutes with $U_{w}$ for every $w \in
\mathbb{T}^{\ell +1}$.

The operator $D_{\ell}$ is optimal i.e.\ if $(\cla(S_{q}^{2\ell+1}),\clh_{\ell },D)$ is a
spectral triple such that  $D$  commutes with
$U_{w}$ for every $w$, then there exist positive reals $a$ and $b$ such that
$| D | \leq a + b | D_{\ell} |$.
\ethm
%%%%%%%%%%%%%%%%%%%%%%%%%%%%%%%%%%%%%%%%%%%%%

In the next few subsections, 
we will introduce a dense subalgebra $\cla_\ell^\infty$
of $A_\ell(q)$ closed under its holomorphic function calculus
and establish regularity of the spectral triple $(\cla_\ell^\infty,\clh_\ell,D_\ell)$.
We will also compute its  dimension spectrum.

%%%%%%%%%%%%%%%%%%%%%%%%%%%%%%%%%%%%%%%%%%%%%%%%%%%%%%%%%%%%%%%%%%
\subsection{The smooth function algebra $\cla_\ell^\infty$}
%%%%%%%%%%%%%%%%%%%%%%%%%%%%%%%%%%%%%%%%%%%%%%%%%%%%%%%%%%%%%%%%%%
 In this section we associate a dense Fr\'echet $*$-subalgebra 
of $A_\ell(q)=\pi_\ell(C(S_{q}^{2\ell+1}))$
which is closed under holomorphic functional calculus. 
We will first show that the $C^*$-algebra $A_\ell(q)$ is independent of $q$.
% Observe that $Y_{i,q}$ is norm continuous in the variable  $q$.
%%%%%%%%%%%%%%%%%%%%%%%%%%%%%%%%%%%%%%%%%%%%%
\blmma\label{lem:Aell0}
For any $q\in(0,1)$, one has $A_\ell(0)=A_\ell(q)$.
\elmma
%%%%%%%%%%%%%%%%%%%%%%%%%%%%%%%%%%%%%%%%%%%%%
\prf
Let us first show that $A_\ell(q)\seq A_\ell(0)$.
We will prove this by induction on $\ell$.
Let us denote the generators $Y_{j,q}$ of $A_\ell(q)$ by
$Y_{j,q}^{(\ell+1)}$. Note that for $\ell=0$, one has
$Y_{1,q}^{(1)}=Y_{1,q}^{(1)}$ and $A_0(q)=A_0(0)=C(\bbt)$ so that
the inclusion is trivial.
Next, assume the inclusion for $\ell-1$.
Observe that for $1\leq j\leq \ell$, 
we have $Y_{j+1,q}^{(\ell+1)}=q^N\otimes Y_{j,q}^{(\ell)}$
and
$Y_{j+1,0}^{(\ell+1)}=p_0\otimes Y_{j,0}^{(\ell)}$.
From this last equality and from the induction hypothesis,
it follows that
$p_0\otimes Y_{j,q}^{(\ell)}\in A_\ell(0)$ for $1\leq j\leq \ell$.
Since for $1\leq j\leq \ell$,
\[
 Y_{j+1,q}^{(\ell+1)}=q^N\otimes Y_{j,q}^{(\ell)}
    =\sum_{n\in\bbn}q^n (Y_{1,0}^{(\ell+1)})^n (p_0\otimes Y_{j,q}^{(\ell)})
           (Y_{1,0}^{(\ell+1)})^{*n},
\]
it follows that $Y_{j,q}^{(\ell+1)}\in A_\ell(0)$ for $2\leq j\leq\ell+1$.
So it remains to show that $Y_{1,q}^{(\ell+1)}\in A_\ell(0)$.
Note that $Y_{1,q}^{(\ell+1)}=(\sqrt{I-q^{2N}}\otimes I)Y_{1,0}^{(\ell+1)}$,
and
\[
 q^N\otimes I=\sum_{n\in\bbn}q^n (Y_{1,0}^{(\ell+1)})^n (p_0\otimes I)
           (Y_{1,0}^{(\ell+1)})^{*n}
   =\sum_{n\in\bbn}q^n (Y_{1,0}^{(\ell+1)})^n 
              (Y_{2,0}^{(\ell+1)})^*(Y_{2,0}^{(\ell+1)})(p_0\otimes I)
           (Y_{1,0}^{(\ell+1)})^{*n}.
\]
Therefore we have the required inclusion.

% % % ALTERNATIVE PROOF:
% % % 
% % For $j,k\in\bbn$ with $j+k\leq\ell$, let
% % \[
% %  X_{j,k}=\underbrace{p_0\otimes\cdots\otimes p_0}_{j}\otimes
% %        \underbrace{q^N\otimes\cdots\otimes q^N}_{k}\otimes S^*\times
% %        \underbrace{I\otimes\cdots\otimes I}_{\ell-j-k}.
% % \]
% % Then
% % $X_{j,0}=Y_{j+1,0}$ for $0\leq j\leq \ell$,
% % and
% % $X_{j,k}=\sum_{n\in\bbn} q^n Y_{j,0}^n X_{j+1,k-1} Y_{j,0}^{*n}$.
% % Therefore it follows that $X_{j,k}\in A_\ell(0)$ for 
% % $0\leq j,k\leq \ell$ with $j+k\leq \ell$.
% % Now observe that
% % $|Y_{k,q}|=(|X_{0,k-1}|^2-|X_{0,k}|^2)^\half$
% % for $1\leq k\leq \ell$ and 
% % $Y_{\ell+1,q}=X_{0,\ell}$.
% % Thus $|Y_{k,q}|\in A_\ell(0)$ for all $1\leq k\leq \ell+1$.
% % Let $A_n=X_{0,k-1}|X_{0,k-1}|^{-1+\frac{1}{n}}$.
% % By lemma~4.4 in \cite{lan-1995a}, one has
% % $A_n\in A_\ell(0)$ for all $n>1$.
% % Observe that $q^{\frac{1}{n}N}$ converges in strong operator topology
% % to $I$ as $n$ goes to infinity. Since $q^N$ is compact, this implies
% % $q^{\frac{1}{n}N}q^N$ converges to $q^N$ in norm.
% % It follows from this that $A_n|Y_{k,q}|$ converges to $Y_{k,q}$ in norm.
% % Since $A_n|Y_{k,q}|\in A_\ell(0)$, one has $Y_{k,q}\seq A_\ell(0)$.
% % Thus $A_\ell(q)\seq A_\ell(0)$.

For the other inclusion, we will use the following
fact: if $B$ denotes the $C^*$-subalgebra of $\cll(\ell_2(\bbn))$
generated by the operator $X=(1-q^{2N})^\half S^*$,
then $B$ contains the shift operator $S$.
This is because the operator $|X|$ is invertible and
$S^*=X|X|^{-1}$.
Using this fact for the first copy of $\ell_2(\bbn)$,
since $Y_{1,q}\in A_\ell(q)$, one gets $Y_{1,0}\in A_\ell(q)$.
Next assume  that $Y_{i,0}\in A_\ell(q)$ for
$1\leq i\leq j-1$, where $2\leq j\leq\ell$.
Then  $P_{j-1}:=I-\sum_{k=1}^{j-1}Y_{k,0}Y_{k,0}^*\in A_\ell(q)$.
Observe that
\[ 
P_{j-1}Y_{j,q}=\underbrace{p_0\otimes\cdots\otimes p_0}_{j-1}
    \otimes X\otimes\underbrace{I\otimes\cdots\otimes I}_{\ell+1-j},\qquad
Y_{j,0}=\underbrace{p_0\otimes\cdots\otimes p_0}_{j-1}
    \otimes S^*\otimes\underbrace{I\otimes\cdots\otimes I}_{\ell+1-j}.
\]
Therefore using the above fact for the j\raisebox{.4ex}{th} copy of $\ell_2(\bbn)$,
we get $Y_{j,0}\in A_\ell(q)$.
Finally, since $Y_{\ell+1,0}=Y_{\ell+1,q}(I-\sum_{k=1}^{\ell}Y_{k,0}Y_{k,0}^*)$,
one has $Y_{\ell+1,0}\in A_\ell(q)$.
\qed

Let us write $\alpha_{i}$ for $Y_{i,0}^{*}$. 
Note that the $C^{*}$-subalgebra of $A_\ell$ generated by
$\alpha_{2},\cdots,\alpha_{\ell +1}$ is isomorphic to $A_{\ell-1}$ where the
map $a \mapsto p_0 \otimes a$ gives the isomorphism.  
We define the Fr\'echet subalgebras $\cla_\ell^\infty$
inductively as follows.

The algebra 
\[
 \cla_0^{\infty}:=\left\{~\sum_{n \in \mathbb{Z}}a_{n}\bldz^{n}:
\text{$(a_{n})$ is rapidly decreasing} \right\}
\]
 is the algebra of smooth functions
on $\mathbb{T}$ together with the increasing family of seminorms 
$\| \cdot\|_{p}$ given by 
$\|(a_{n})\|_{p}=\sum~ (1+| n |)^{p} | a_{n} |$. 
Then $\cla_0^\infty$ is a dense $*$ Fr\'echet
subalgebra of $A_0=C(\bbt)$. Note that $\| a \| \leq \| a
\|_{0}$ for $a \in \cla_0^\infty$. Now assume that
$(\cla_{\ell-1}^{\infty}, \|\cdot\|_{m})$ be defined such that
\begin{enumerate}
\item the seminorms $\| \cdot \|_{m}$ are increasing and
$(\cla_{\ell-1}^{\infty}, \|\cdot\|_{m})$ is a Fr\'echet algebra,
\item the subalgebra $\cla_{\ell-1}^{\infty}$ is $*$ closed and dense in
$A_{\ell-1}$. For every $a \in \cla_{\ell-1}^{\infty}$, one has
$\| a^{*} \|_{m}=\| a \|_{m}$,
\item for every $a \in \cla_{\ell-1}^{\infty}$, one has
$\| a \| \leq \| a \|_{0}$ where $\|\cdot\|$ denotes the $C^{*}$ norm
of $A_{\ell-1}$.
\end{enumerate}
Now define
\bea
\cla_\ell^\infty &:= & \left\{ \sum_{j,k \in \mathbb{N}}\alpha_{1}^{*j}(p_0
\otimes a_{jk})\alpha_{1}^{k} + \sum_{k\geq 0}\lambda_{k}\alpha_{1}^{k}+
\sum_{k >0}\lambda_{-k}\alpha_{1}^{*k}~:~a_{jk} \in \cla_{\ell-1}^\infty,\right.\nonumber\\
&& \hspace{2em}\left.\sum_{j,k}(1+j+k)^{n} \| a_{jk} \|_{m}<\infty,~(\lambda_{k})
\text{ is rapidly decreasing} \right\}.\label{eq:Aellinfty}
\eea
Let $a:= \sum_{j,k}\alpha_{1}^{*j}(p_0 \otimes
a_{jk})\alpha_{1}^{k}+\sum_{k\geq
0}\lambda_{k}\alpha_{1}^{k}+\sum_{k >0}\lambda_{-k} \alpha_{1}^{*k}$ ~be an
element of $\cla_\ell^\infty$. Define for $m \in \mathbb{N}$, the seminorms
$\|a \|_{m}$ as follows:
\[
\| a \|_{m} = \max_{r,s \leq m}( \sum_{j,k}(1+j+k)^{r} \|
a_{jk} \|_{s}) + \sum_{k \in \mathbb{Z}}(1+| k |)^{m} |
\lambda_{k} |.
\]
%%%%%%%%%%%%%%%%%%%%%%%%%%%%%%%%%%%%%%%%%%%%%%%%
\bppsn\label{pro:propAellinfty}
The pair $(\cla_\ell^\infty,\|\cdot \|_{m})$ has the following
properties:
\begin{enumerate}
\item the seminorms $\| \cdot \|_{m}$ are increasing and
$(\cla_\ell^\infty,\|\cdot \|_{m})$ is a Fr\'echet algebra,
\item the subalgebra $\cla_\ell^\infty$ is $*$ closed and dense in
$A_\ell$. For every $a \in \cla_\ell^\infty$, one has
$\| a^{*} \|_{m}=\| a \|_{m}$,
\item for every $a \in \cla_\ell^\infty$, one has
$\| a \| \leq \| a \|_{0}$ where $\|\cdot\|$ denotes the $C^{*}$ norm
of $A_\ell$.
\end{enumerate}
\eppsn
%%%%%%%%%%%%%%%%%%%%%%%%%%%%%%%%%%%%%%%%%%%%%%%%
\prf
The proof is by induction on $\ell$. Parts $(2)$ and $(3)$ 
and the fact that the seminorms $\| \cdot \|_{m}$ are 
increasing follow from the definition and  the induction 
hypothesis. One verifies directly that 
$(\mathcal{A}_{\ell}^{\infty},\| \cdot \|_{m})$ is a Fr\'echet algebra 
using induction and the following relations.
\begin{align*}
 \alpha_{1}\alpha_{1}^{*}&=1\\
 \alpha_{1}^{*j}(p_{0}\otimes a_{jk})\alpha_{1}^{k} 
     \alpha_{1}^{*r}(p_{0}\otimes a_{rs})\alpha_{1}^{s}
  &= \delta_{kr}\alpha_{1}^{*j}(p_{0}\otimes a_{jk}a_{rs})\alpha_{1}^{s}
\end{align*}
\[
\alpha_{1}^{*m}\alpha_{1}^{n}=\begin{cases}
    (\alpha_{1}^*)^{m-n} - \sum_{k=0}^{n-1} (\alpha_{1}^*)^{m-n+k}(p_{0}\otimes 1)\alpha_{1}^{k}
                        & \text{ if $m\geq n$}\\
  \alpha_{1}^{n-m} - \sum_{k=0}^{m-1} \alpha_{1}^{*k}(p_{0}\otimes 1)\alpha_{1}^{n-m+k}
             & \text{ if $m< n$}.
                         \end{cases}
\]
\qed
%%%%%%%%%%%%%%%%%%%%%%%%%%%%%%%%%%%%%%%%%%%%%%%%

Denote the generators $z_{1},z_{2},\cdots
z_{\ell +1}$ of $C(S_{q}^{2\ell+1})$ by 
$z_{1}^{(\ell +1)},z_{2}^{(\ell +1)},\cdots,z_{\ell +1}^{(\ell +1)}$.
Let $\sigma_\ell:C(S_{q}^{2\ell+1}) \to C(S_{q}^{2\ell-1})$ be the homomorphism
given by $\sigma_\ell(z_{\ell +1}^{(\ell +1)})=0$ and 
$\sigma_\ell(z_{i}^{(\ell +1)})=z_{i}^{(\ell) }$ for $1 \leq i \leq \ell$.
Let us denote by the same symbol $\sigma_\ell$ the induced homomorphism
from $A_\ell$ to $A_{\ell-1}$.
Observe that if one applies the map $\sigma$ on the
$\ell$\raisebox{.4ex}{th} copy of $\scrt$ in $\scrt_\ell\otimes C(\bbt)$
followed by evaluation at 1 in the $(\ell+1)$\raisebox{.4ex}{th} copy,
then the restriction of the resulting map to $A_\ell$ is precisely $\sigma_\ell$.

%%%%%%%%%%%%%%%%%%%%%%%%%%%%%%%%%%%%%%%%%%%%%%%%
\bppsn
 The dense Fr\'echet $*$-subalgebra $\cla_\ell^{\infty}$ of $A_\ell$
is closed under holomorphic functional calculus in $A_\ell$. Moreover, 
the algebra
$\cla_\ell^{\infty}$ contains the generators
$Y_{1,q}^{(\ell +1)},\cdots, Y_{\ell +1,q}^{(\ell +1)}$.
\eppsn
%%%%%%%%%%%%%%%%%%%%%%%%%%%%%%%%%%%%%%%%%%%%%%%%
\prf
We prove this proposition by induction on $\ell$. For $\ell =0$, by
definition $\cla_0^\infty=C^{\infty}(\mathbb{T})$. Hence the proposition
is clear in this case. Now assume that the algebra $\cla_{\ell-1}^\infty$ is
closed under holomorphic functional calculus in $A_{\ell-1}$ and contains
$Y_{1,q}^{(\ell)},\cdots, Y_{\ell ,q}^{(\ell )}$. 
The homomorphism $\sigma_\ell:A_\ell\to A_{\ell-1}$ gives the following exact
sequence 
\begin{displaymath}
0 \longrightarrow \clk(\ell_{2}(\mathbb{N}^{\ell })) \otimes C(\mathbb{T}) \longrightarrow A_\ell \longrightarrow
A_{\ell-1} \longrightarrow 0
\end{displaymath}
One also has at the smooth algebra level the ``sub''
extension  
\begin{displaymath}
0 \longrightarrow \cls(\ell_2(\mathbb{N}^{\ell })) \otimes C^{\infty}(\mathbb{T}) \longrightarrow
\cla_\ell^\infty \longrightarrow \cla_{\ell-1}^\infty \longrightarrow 0.
\end{displaymath}
Since $\cls(\ell_2(\mathbb{N}^{\ell }))\otimes
C^{\infty}(\mathbb{T}) \subset \clk(\ell_{2}(\mathbb{N}^{\ell }))\otimes C(\mathbb{T})$
and $\cla_{\ell-1}^\infty \subset A_{\ell-1}$ are closed under
the respective holomorphic functional calculus, it follows 
from theorem~3.2, part~2, \cite{sch-1993a} that 
$\cla_\ell^\infty$ is spectrally invariant in $A_\ell$.
Since $\|a\|\leq \|a\|_0$ for all $a\in\cla_\ell^\infty$,
it follows that the Fr\'echet topology of $\cla_\ell^\infty$
is finer than the norm topology.
Therefore $\cla_\ell^\infty$
 is closed under holomorphic functional calculus in $A_\ell$.
Observe that for $i \geq 2$, we have
$Y_{i,q}^{(\ell +1)}= \sum_{n \geq 0} q^{n}
\alpha_{1}^{*n}(p_0 \otimes Y_{i-1,q}^{(\ell) })\alpha_{1}^{n}$. 
Hence $Y_{i,q}^{(\ell +1)} \in \cla_\ell^\infty$ for $i=2,\cdots, \ell +1$. 
Also note that 
$q^{N}\otimes I= \sum_{n \geq 0}q^{n}\alpha_{1}^{*n}(p_0 \otimes 1)\alpha_{1}$. Since
$\cla_\ell^\infty$ is closed under holomorphic functional calculus, it
follows that $\sqrt{1-q^{2N+2}}\otimes I \in \cla_\ell^\infty$. As
$Y_{1,q}^{(\ell +1)}=\alpha_{1}^{*}(\sqrt{1-q^{2N+2}} \otimes I)$ it follows that
$Y_{1,q}^{(\ell +1)} \in \cla_\ell^\infty$. This completes the proof. 
\qed

Next we proceed to prove that  the spectral triple
$(\cla_\ell^\infty,\clh_{\ell },D_{\ell})$ is regular and compute its
dimension spectrum. The proof is by induction. We start with the case $\ell =0$
to start the induction.

%%%%%%%%%%%%%%%%%%%%%%%%%%%%%%%%%%%%%%%%%%%%%%%%%%%%%%%%%%%%%%%%%%
%%%%%%%%%%%%%%%%%%%%%%%%%%%%%%%%%%%%%%%%%%%%%%%%%%%%%%%%%%%%%%%%%%

%%%%%%%%%%%%%%%%%%%%%%%%%%%%%%%%%%%%%%%%%%%%%%%%%%%%%%%%%%%%%%%%%%
\subsection{The case $\ell =0$}
%%%%%%%%%%%%%%%%%%%%%%%%%%%%%%%%%%%%%%%%%%%%%%%%%%%%%%%%%%%%%%%%%%
%%%%%%%%%%%%%%%%%%%%%%%%%%%%%%%%%%%%%%%%%%%%%%%%%%%%%%%%%%%%%%%%%%
 For $\ell =0$, the spectral triple $(\cla_0^\infty,\clh_{0},D_{0})$ is
unitarily equivalent to the spectral triple
$(C^{\infty}(\mathbb{T}),L_2(\mathbb{T}),\frac{1}{i} \frac{d}{d\theta})$.
For $f \in C^{\infty}(\mathbb{T})$ one has $[D_0,f]=\frac{1}{i}f'$. Let
$(e_{k})$ be the standard orthonormal basis and let $p_{k}$ be the projection
onto $e_{k}$. Let $F_0:=sign(D_{0})$. Note that $[F_0,\bldz]=2p_{0}\bldz$ and hence by
induction $[F_0,\bldz^{n}]= 2 \sum_{k=0}^{n-1}p_{k}\bldz^{n}p_{k-n}$ for $n \geq 0$.
Thus $[F_0,\bldz^{n}]$ is smoothing for $n \geq 0$. Also 
$\| | D_0 |^{r}[F_0,\bldz^{n}]| D_0 |^{s} \| \leq 2 (1+n)^{r+s+1}$. Since
$[F_0,\bldz^{-| n |}]^{*}=-[F_0,\bldz^{| n |}]$, it follows that $[F_0,\bldz^{n}] \in
OP^{-\infty}$ for every $n$. Also $\| [F_0,\bldz^{n}] \|_{r,s} \leq
2(1+| n |)^{r+s+1}$. Hence we observe that $[F_0,f] \in OP^{-\infty}$ and
$\| [F_0,f] \|_{r,s} \leq 2 \| f \|_{r+s+1}$. Let
$\delta$ be the unbounded derivation $[| D_0 |,\cdot]$.

%%%%%%%%%%%%%%%%%%%%%%%%%%%%%%%%%%%%%%%%%%%%%%%%%
\blmma
\label{l=0 case}
 Let $\mathcal{B}:=\{ f_{0}+f_{1}F_0 +R:  f_{0},f_{1} \in C^{\infty}(\mathbb{T}),
R \in OP_{D_0}^{-\infty}\}$. Then
\begin{enumerate}
\item If $f_{0}+f_{1}F_0$ is smoothing then $f_{0}=f_{1}=0$. Hence $\mathcal{B}$
is isomorphic to the direct sum $C^{\infty}(\mathbb{T})\oplus
C^{\infty}(\mathbb{T})\oplus OP_{D_0}^{-\infty}$. We give $\mathcal{B}$ the Fr\'echet
space structure coming from this decomposition. This topology on
$\mathcal{B}$ is generated by the seminorms $(\|\cdot\|_{m})_{m \in \mathbb{N}}$
which are defined by $\| f_{0}+f_{1}F_0+R \|_{m}:=
\| f_{0} \|_{m} + \| f_{1} \|_{m}+ \sum_{r+s \leq
m}\| R \|_{r,s}$.
\item The vector space $\mathcal{B}$ is closed under $\delta$ and the derivation
$[D_0,\cdot]$.
\item For every $b \in \mathcal{B}$, $[F_0,b] \in OP^{-\infty}$. Also the map $b
\to [F_0,b] \in OP^{-\infty}$ is continuous. The derivations $\delta$ and $[D_0,\cdot]$
are continuous.
\item The vector space $\mathcal{B}$ is an algebra and contains
$C^{\infty}(\mathbb{T})$.
\end{enumerate} 
\elmma
%%%%%%%%%%%%%%%%%%%%%%%%%%%%%%%%%%%%%%%%%%%%%%%%%
\prf
First observe that a bounded operator $T$ on
$\ell_{2}(\mathbb{Z})$ is smoothing if and only if $(\langle Te_{m},e_{n}\rangle)_{m,n}$ 
is rapidly decreasing. Now suppose that $R:=f_{0}+f_{1}F_0$ be
smoothing. Fix an integer $r$. Observe that $\langle R(e_{n}), e_{r+n}\rangle$ 
converges to $\hat{f_{0}}(r)+\hat{f_{1}}(r)$ as $n \to +\infty$ and
converges to $\hat{f_{0}}(r)-\hat{f_{1}}(r)$ as $n \to -\infty$. But since $R$
is smoothing it follows that
$\hat{f_{0}}(r)+\hat{f_{1}}(r)=0=\hat{f_{0}}(r)-\hat{f_{1}}(r)$. Hence
$\hat{f_{0}}(r)=\hat{f_{1}}(r)=0$ for every integer $r$. Thus $f_{0}=f_{1}=0$.
This proves part~(1).

Parts~(2), (3) and~(4) follow from the observations that
$[D_0,f]=\frac{1}{i}f'$, $[F_0,f] \in OP^{-\infty}$, $\| [F_0,f]
\|_{r,s} \leq 2 \| f \|_{r+s+1}$ and
$\delta(b)=[D_0,b]F_0+D_0[F_0,b]$. This completes the proof.
\qed\\
In particular, it follows from parts~(2) and (4) of the above lemma
that the  spectral triple $(\cla_0^\infty,\clh_0,D_0)$ is regular.

Let $\cle$ be the $C^{*}$-subalgebra of $\cll(\ell_{2}(\mathbb{Z}))$ generated by
$C(\mathbb{T})$ and $F_0$.  
Note that the algebra $\clb$ plays the role of smooth function subalgebra
for the $C^*$-algebra $\cle$. Therefore  $\cle^\infty$ will stand for
the algebra $\clb$.
%%%%%%%%%%%%%%%%%%%%%%%%%%%%%%%%%%%%%%%%%%%%%%%%%%%%%%%%%%%%%%%%%%
%%%%%%%%%%%%%%%%%%%%%%%%%%%%%%%%%%%%%%%%%%%%%%%%%%%%%%%%%%%%%%%%%%

%%%%%%%%%%%%%%%%%%%%%%%%%%%%%%%%%%%%%%%%%%%%%%%%%%%%%%%%%%%%%%%%%%
\subsection{Regularity and the dimension spectrum}
%%%%%%%%%%%%%%%%%%%%%%%%%%%%%%%%%%%%%%%%%%%%%%%%%%%%%%%%%%%%%%%%%%
%%%%%%%%%%%%%%%%%%%%%%%%%%%%%%%%%%%%%%%%%%%%%%%%%%%%%%%%%%%%%%%%%%
In this subsection we prove regularity and calculate the dimension spectrum for
the  spectral triple 
$(\cla_\ell^{\infty},\clh_{\ell},D_{\ell})$. 
The proof is by induction on $\ell$.
Let us denote the derivation 
$[| D_{\ell } |,\cdot]$ by $\delta_{\ell }$ and let $F_{\ell }$ stand 
for the sign of the operator $D_{\ell }$. Observe that
 $F_{\ell }=1^{\otimes \ell } \otimes F_{0}=1\otimes F_{\ell-1}$.

%%%%%%%%%%%%%%%%%%%%%%%%%%%%%%%%%%%%%%%%%%%%%%%%%%%%%%%%%%%%%%%%%%
\bppsn
\label{regularity}
 Let $\mathcal{B}_{\ell }:=\{A_{0}+A_{1}F_{\ell}+R:~A_{0},A_{1}\in
\cla_\ell^{\infty},R\in OP^{-\infty}\}$. Then
\begin{enumerate}
\item if $A_{0}+A_{1}F_{\ell}$ is smoothing then $A_{0}=A_{1}=0$. Hence
$\mathcal{B}_{\ell }$ is isomorphic to the direct sum
$\cla_\ell^\infty\oplus \cla_\ell^\infty\oplus OP^{-\infty}$. Equip
$\mathcal{B}_{\ell }$ with the Fr\'echet space structure coming from this
decomposition. This topology on $\mathcal{B}_{\ell }$ is induced by the seminorms
$(\|\cdot\|_{m})_{m \in \mathbb{N}}$ which are defined by 
$\| A_{0}+A_{1}F_{\ell}+R \|_{m}:= \| A_{0} \|_{m} +
\| A_{1} \|_{m}+ \sum_{r+s \leq m}\| R \|_{r,s}$.
\item For every $b \in \mathcal{B}_{\ell }$, $[F_{\ell},b] \in OP^{-\infty}$. Also the map
$b \to [F_{\ell},b] \in OP^{-\infty}$ is continuous. 
\item The vector space $\mathcal{B}_{\ell }$ is closed under the derivations
$\delta_{\ell }$ and $[D_{\ell },\cdot]$. Moreover the derivations $\delta_{\ell }$ and
$[D_{\ell },\cdot]$ are continuous.
\item The vector space $\mathcal{B}_{\ell }$ is an algebra and contains
$\cla_\ell^\infty$.
\end{enumerate}
\eppsn
%%%%%%%%%%%%%%%%%%%%%%%%%%%%%%%%%%%%%%%%%%%%%%%%%%%%%%%%%%%%%%%%%%
\prf
The proof is by induction on $\ell$. For $\ell =0$, the proposition is
just lemma \ref{l=0 case}. Now assume that the proposition is true for $\ell-1$.
Suppose that $A_{0}+A_{1}F_{\ell}$ is
smoothing for some $A_{0},A_{1}\in \cla_\ell^\infty$. 
Then $A_{0}+A_{1}F_{\ell}\in \scrt_\ell \otimes \cle$ 
and $A_{0}+A_{1}F_{\ell}$ is compact. Therefore 
$(\sigma\otimes id)(A_{0}+A_{1}F_{\ell})=0$. Now let 
\[
A_{i}= \sum_{j,k \geq 0}
\alpha_{1}^{*j}(p_0\otimes a_{jk}^{(i)})\alpha_{1}^{k} + \sum_{k \geq 0}
\lambda_{k}^{(i)}\alpha_{1}^{k} + \sum_{k >0}\lambda_{-k}^{(i)}\alpha_{1}^{*k}
\]
for $i=0,1$. Let $f_{i}(z)=\sum_{k \in \mathbb{Z}}\lambda_{k}^{(i)}z^{k}$ for
$i=0,1$. Now 
$(\sigma\otimes id)(A_{0}+A_{1}F_{\ell})=f_{0} \otimes I + f_{1}\otimes F_{\ell-1} $.
So  we have $f_{0} \otimes I + f_{1}\otimes F_{\ell-1}=0$. 
Writing $F_{\ell}=2P_{\ell}-I$, it follows that 
$(f_{0}+f_{1})\otimes P_{\ell-1} + (f_{0}-f_{1})\otimes (1-P_{\ell-1})=0$. 
Hence $f_{0}=f_{1}=0$. This shows that
$\lambda_{k}^{(i)}=0$ for $i=0,1$. Let $b_{jk}=a_{jk}^{0}+a_{jk}^{1}F_{\ell-1}$. Since
$R:=A_{0}+A_{1}F_{\ell}$ is smoothing, it follows that for every $j,k$, the matrix
entries $\langle e_{(j,\gamma)}, R(e_{(k,\gamma')})\rangle$ are rapidly decreasing in
$(\gamma,\gamma')$. Hence $b_{jk}$ is smoothing for every $j,k$. By
induction hypothesis $a_{jk}^{(i)}=0$ for every $j,k \geq 0$ and for $i=0,1$.
Thus $A_{0}=A_{1}=0$. This proves part~(1).

Observe that 
\begin{displaymath}
 \delta_{\ell }(\alpha_{1})= -\alpha_{1},\qquad
 | D_{\ell } |^{r} \alpha_{1}^{*k}=\alpha_{1}^{*k}(| D_{\ell } | +k)^{r}, \qquad
 \alpha_{1}^{k}| D_{\ell }|^{s}=(| D_{\ell } | + k)^{s}\alpha_{1}^{k}. 
 \end{displaymath}
Also $F_\ell$ commutes with $\alpha_{1}$. To prove $(2)$, it is enough to show that
$[F_\ell,a]$ is smoothing for every $a \in \cla_\ell^\infty$ and the map 
$a\mapsto [F_\ell,a]$ is continuous. Let 
\[
a=\sum_{m,n \geq 0} \alpha_{1}^{*m}(p_0 \otimes
a_{mn})\alpha_{1}^{n} + \sum_{m \geq 0} \lambda_{m} \alpha_{1}^{m} +
\sum_{m>0} \lambda_{-m} \alpha_{1}^{*m}
\]
be an element in $\cla_\ell^\infty$. Then 
$[F_\ell,a] = \sum_{m,n \geq 0} \alpha_{1}^{*m}(p_0\otimes [F_{\ell-1},a_{mn}])\alpha_{1}^{n}$. 
By induction hypothesis, it follows that 
$p_0\otimes [F_{\ell-1},a_{mn}]$ is smoothing for every $m,n \geq 0$. Since
$(OP_{D_\ell}^{-\infty},\|\cdot\|_{r,s})$ is a Fr\'echet space, 
to show that $[F_\ell,a]$ is smoothing it is enough to show that
the infinite sum $\sum_{m,n \geq 0} \alpha_{1}^{*m}(p_0 \otimes
[F_{\ell-1},a_{mn}])\alpha_{1}^{n}$ converges absolutely in every seminorm
$\|\cdot\|_{r,s}$. Now observe that 
\begin{equation}\label{smoothing}
 | D_{\ell } |^{r}\alpha_{1}^{*m}(p_0 \otimes
[F_{\ell-1},a_{mn}])\alpha_{1}^{n} | D_{\ell } |^{s} = 
\alpha_{1}^{*m}(|
D_{\ell } | + m)^{r}(p_0 \otimes [F_{\ell-1},a_{mn}])(| D_{\ell } | +n)^{s} \alpha_{1}^{n}.
\end{equation}
Since the map $ a'\in \cla_{\ell-1}^\infty \mapsto [F_{\ell-1},a']\in OP^{-\infty}$
is continuous, there exist $p\in\bbn$ and $C_{p}>0$ such that 
$\| [F_{\ell-1},a']\|_{i,j} \leq C_{p} \| a' \|_{p}$ 
for every $a' \in \cla_{\ell-1}^{\infty}$ 
and for $i,j \leq max\{r,s\}$. Hence by equation~(\ref{smoothing}), it follows that 
\begin{align*}
\sum_{m,n} \| \alpha_{1}^{*m}(p_0 \otimes [F_{\ell-1},a_{mn}])\alpha_{1}^{n}
\|_{r,s} & \leq \sum_{m,n}
\sum_{i=0}^{r}\sum_{j=0}^{s}\binom{r}{i}\binom{s}{j}m^{r-i}n^{s-j} \|
[F_{\ell-1},a_{mn}] \|_{i,j}\\
        & \leq \sum_{i=0}^{r}\sum_{j=0}^{s}
\binom{r}{i}\binom{s}{j}C_{p}
\left(\sum_{m,n} m^{r}n^{s} \| a_{mn}\|_{p}\right).  
\end{align*}
This shows that $[F_\ell,a]$ is smoothing and the above inequality also shows that
for every $r,s \geq 0$, there exists $t \geq 0$ and a $C_{t}>0$ such that
$\| [F_\ell,a] \|_{r,s} \leq C_{t} \| a \|_{t}$. Hence
the map $a \mapsto [F_\ell,a]$ is continuous. This proves $(2)$.

To show (3),  it is enough to show that the map 
$a  \mapsto \delta_{\ell }(a)$ from $\cla_\ell^\infty$ to $\mathcal{B}_{\ell }$ 
makes sense and is continuous. We will use the fact that the unbounded derivation $\delta_{\ell }$
is a closed derivation. Let $a=\sum_{m,n \geq 0} \alpha_{1}^{*m}(p_0 \otimes
a_{mn})\alpha_{1}^{n} + \sum_{m \geq 0} \lambda_{m} \alpha_{1}^{m} +
\sum_{m>0} \lambda_{-m} \alpha_{1}^{*m}$ be an element in
$\cla_\ell^\infty$. Since $\alpha_{1}$ and $p_0 \otimes a_{mn} \in
Dom(\delta_{\ell })$ it follows that each of the terms in the infinite sum is an
element in $Dom(\delta_{\ell })$. Hence in order to show $a \in Dom(\delta_{\ell })$, it is
enough to show that the sum
\[
\sum_{m,n} \delta_{\ell }(\alpha_{1}^{*m}(p_0 \otimes
a_{mn})\alpha_{1}^{n}) + \sum_{m \geq 0}\lambda_{m} \delta_{\ell }(\alpha_{1}^{m})
+ \sum_{n>0}\lambda_{-n} \delta_{\ell }(\alpha_{1}^{*n})
\]
converges. Observe that 
$\delta_{\ell }(\alpha_{1}^{*m})=m\alpha_{1}^{*m}$,
$\delta_{\ell }(\alpha_{1}^{n})=-n\alpha_{1}^{n}$  and
\[
 \delta_{\ell }\bigl(\alpha_{1}^{*m}(p_0 \otimes a_{mn})\alpha_{1}^{n}\bigr)=
(m-n)\alpha_{1}^{*m}(p_0 \otimes a_{mn})\alpha_{1}^{n} + \alpha_{1}^{*m}(p_0
\otimes \delta_{\ell -1}(a_{mn})) \alpha_{1}^{n}
\]
Since $\delta_{\ell -1}$ is continuous, it follows that 
$\| \delta_{\ell -1}(a_{mn}) \|$  is rapidly decreasing where
$\|\cdot\|$ is the operator norm. (Note that for $b \in
\clb_{\ell }$, one has $\| b \| \leq \| b \|_{0}$.) Hence
the infinite sum  
\[
\sum_{m,n} \delta_{\ell }\bigl(\alpha_{1}^{*m}(p_0 \otimes
a_{mn})\alpha_{1}^{n}\bigr) 
+ \sum_{m \geq 0}\lambda_{m} \delta_{\ell }(\alpha_{1}^{m})
+ \sum_{n>0}\lambda_{-n} \delta_{\ell }(\alpha_{1}^{*n})
\]
converges absolutely in the operator norm.
Therefore $a \in Dom(\delta_{\ell })$. Since $\delta_{\ell -1}$ is continuous for every $r$
there exists $p$ and $C_{p}$ such that $\| \delta_{\ell -1}(a')
\|_{r} \leq C_{p} \| a' \|_{p}$.  Write
$\delta_{\ell -1}(a_{mn})$ as $\delta_{\ell -1}(a_{mn})=a'_{mn}+a''_{mn}F_{\ell}+R_{mn}$.
 Let 
%%%%%%%%%%%%%%%%%%%%%%%%%%%%%%%%%%%%%%%%%%%%%%%%%%%%
\begin{align*}
A_{0}& = \sum_{m,n}\alpha_{1}^{*m}(p_0
\otimes((m-n)a_{mn}+a'_{mn}))\alpha_{1}^{n}+ \sum_{m \geq 0}
m\lambda_{m}\alpha_{1}^{m} + \sum_{n>0} (-n)\lambda_{-n}\alpha_{1}^{*n},\\
A_{1}&=\sum_{m,n} \alpha_{1}^{*m}(p_0 \otimes a''_{mn})\alpha_{1}^{n}, \\
R &= \sum_{m,n}\alpha_{1}^{*m}(p_0 \otimes R_{mn})\alpha_{1}^{n}.
\end{align*}
%%%%%%%%%%%%%%%%%%%%%%%%%%%%%%%%%%%%%%%%%%%%%%%%%%%%
Then $\delta_{\ell }(a)=A_{0}+A_{1}F_{\ell}+R$.
In every seminorm of $\cla_{\ell-1}^\infty$ the double sequence
$(a'_{mn})$ and $(a''_{mn})$ are rapidly decreasing. Also $R_{mn}$ is
rapidly decreasing in every seminorm of $OP_{D_\ell}^{-\infty}$. 
Hence $A_{0},A_{1} \in \cla_{\ell}^\infty$ and 
as in the proof of $(2)$, it follows that $R$ is
smoothing and given $r,s$ there exists $t$ and $C_{t}$ such that $\| R
\|_{r,s} \leq C_{t} \| a \|_{t}$. Fix an $r \geq 0$ and
choose $t>1+r$ and $C_{t}>1$ such that $\| \delta_{\ell -1}(a')
\|_{r} \leq C_{t} \| a'\|_{t}$ for every 
$a' \in \cla_{\ell-1}^\infty$. Now $\| A_{0}\|_{r} \leq C_{t}
\| a \|_{t}$ and $\| A_{1} \|_{r} \leq C_{t}
\| a \|_{t}$. This shows that the map $a \to \delta_{\ell }(a) \in
\mathcal{B}_{\ell }$ is continuous. Since $[D_{\ell },b]=\delta_{\ell }(b)F_\ell+| D_{\ell } |
[F_\ell,b]$, the second part of (3) follows as $[F_\ell,b]$ is smoothing by $(2)$. This
proves $(3)$.

Part~(4) follows from $(2)$ and $(3)$.
\qed 

We next prove a lemma that will be crucial in the computation of the dimension spectrum. 
For an $r$ tuple $n=(n_{1},n_{2},\cdots,n_{r}) \in \mathbb{N}^{r}$, 
we will write $| n |$ for $\sum_{i=1}^{r}n_{i}$. 
For $r=0$, we let $\mathbb{N}^{0}=\{0\}$. 

%%%%%%%%%%%%%%%%%%%%%%%%%%%%%%%%%%%%%%%%%%%%%%%%%%%%
\blmma
\label{dimension}
Let $r \geq 0$ and $s \geq 1$ be integers. Let $(a(n))_{n \in \mathbb{N}^{r}}$ be
rapidly decreasing. Then the function 
\[
\xi(z):=\sum_{\substack{n\in\bbn^r,m\in\bbn^s\\| n | + | m |\geq 1}}
               \frac{a(n)}{(| n | + | m |)^{z}}
\]
is meromorphic with simple
poles at $\{1,2,\cdots,s\}$ and $Res_{z=s} \xi(z)= \frac{1}{(s-1)!}\sum_{n}a(n)$.
\elmma
%%%%%%%%%%%%%%%%%%%%%%%%%%%%%%%%%%%%%%%%%%%%%%%%%%%%
\prf
First observe that for $Re\, z>r+s$, 
%%%%%%%%%%%%%%%%%%%%%%%%%%%%%%%%%%%%%%%%%%%%%%%%%%%%
\begin{align*}
\xi(z)&= \sum_{N \geq 1}\frac{1}{N^{z}}\left(\sum_{| n | + | m | =
N}a(n)\right)\\
& = \sum_{N \geq 1} \frac{1}{N^{z}}\left(\sum_{| n | \leq N}
a(n)\sum_{m:| m | = N-| n | } 1\right)\\
&= \sum_{N \geq 1} \frac{1}{N^{z}}\left( \sum_{| n | \leq N} a(n)
\binom{N-| n | + s-1}{s-1}\right).
\end{align*}
%%%%%%%%%%%%%%%%%%%%%%%%%%%%%%%%%%%%%%%%%%%%%%%%%%%%
Note that for a function $(b(n))_{n \in \mathbb{N}^{r}}$ of rapid decay,
the sequence $\left(\sum_{| n | \geq N}b(n)\right)_{N \in \mathbb{N}}$
is of rapid decay. Now $\binom{N-| n |+s-1}{s-1}= \sum_{k=0}^{s-1}
g_{k}(n) N^{k}$ where $g_{k}(n)$ is a polynomial in $(n_{1},n_{2},\ldots,n_{r})$
and $g_{s-1}(n)=\frac{1}{(s-1)!}$. Hence modulo a holomorphic function 
$\xi(z) = \sum_{k=0}^{s-1}(\sum_{n}g_{k}(n)a(n))\zeta(z-k)$. 
Now the result follows from the fact that
$\zeta(z)$ is meromorphic with a simple pole at $z=1$ with residue $1$.
\qed

We will next prove that the spectral triple
$(\cla_\ell^\infty,\clh_\ell,D_{\ell})$ is regular and has discrete dimension spectrum
with simple poles at $\{1,2,\cdots, \ell +1\}$. 
%%%%%%%%%%%%%%%%%%%%%%%%%%%%%%%%%%%%%%%%%%%%%%%%%%
\brmrk
\label{remark}
Recall that the unitaries $U_{w}$ for
$w=(w_{1},w_{2},\cdots,w_{\ell +1}) \in \mathbb{T}^{\ell +1}$ are given by
$U_{w}e_{\gamma}=w_{1}^{\gamma_{1}}w_{2}^{\gamma_{2}}\cdots
w_{\ell +1}^{\gamma_{\ell +1}}e_{\gamma}$. 
A bounded operator $T$ on $\clh_\ell$ is said to be homogeneous of degree
$(m_{1},m_{2},\cdots,m_{\ell +1})$ if
$U_{w}TU_{w}^{*}=w_{1}^{m_{1}}w_{2}^{m_{2}}\cdots w_{\ell +1}^{m_{\ell +1}} T$. If $T$
is homogeneous of degree $(m_{1},m_{2},\cdots,m_{\ell +1}) \neq (0,\ldots,0)$ then 
$Trace(T |D_\ell|^{-z})=0$ if $Re(z)>\ell+1$ since $U_{w}$'s 
commute with the operator $| D_\ell|$.
\ermrk
%%%%%%%%%%%%%%%%%%%%%%%%%%%%%%%%%%%%

%%%%%%%%%%%%%%%%%%%%%%%%%%%%%%%%%%%%
\bppsn
The spectral triple $(\cla_\ell^\infty,\clh_\ell,D_\ell)$ is regular 
and has $\{1,2,\cdots, \ell +1 \}$ as the
dimension spectrum with only simple poles.
\eppsn
%%%%%%%%%%%%%%%%%%%%%%%%%%%%%%%%%%%%
\prf
Regularity of the spectral triple follows from proposition~\ref{regularity}. 
We now prove that for $b\in\clb_\ell$, the function $Trace(b| D_\ell |^{-z})$ is meromorphic with
simple poles at $\{1,2,\cdots, \ell +1\}$. 
Since $\trace(b|D_\ell|^{-z})$ is holomorphic for $b\in OP^{-\infty}$,
we need only to show that for $a\in \cla_\ell^\infty$, the
functions $Trace(a| D_\ell|^{-z})$ and $Trace(aF_\ell| D_\ell |^{-z})$ extend to meromorphic functions with
simple poles at $\{1,2,\cdots, \ell +1\}$. Now any element 
$a \in \cla_\ell^\infty$ can be written as $a=a^{0}+a^{1}$ where $a^{0}$ is
homogeneous of degree 0 and $a^{1}$ is an infinite sum of homogeneous elements
of non zero degrees. Hence by remark~\ref{remark}, 
$Trace(a| D_\ell|^{-z})=Trace(a^{0}| D_\ell |^{-z})$ and 
$Trace(aF_\ell| D_\ell |^{-z})=Trace(a^{0}F_\ell|
D_\ell |^{-z})$. Thus it is enough to consider the functions $Trace(a | D_\ell
|^{-z})$ and $Trace(aF_\ell| D_\ell |^{-z})$ where $a$ is homogeneous of degree $0$.

It is easy to see that the set of homogeneous elements of degree $0$ in
$\cla_\ell^\infty$ is 
\[
\left\{ \sum_{i=0}^{\ell }\left(\sum_{n \in
\mathbb{N}^{i}}\lambda_{n}^{i}(p_{n_{1}} \otimes p_{n_{2}} \otimes \cdots
p_{n_{i}} \otimes 1)\right):~(\lambda_{n}^{i})\text{ is of rapid decay for all }i \right\}
\]
where
$p_{k}=S^{*k}p_0 S^{k}$. Let
$a=\sum_{i=0}^{\ell }(\sum_{n}\lambda_{n}^{i}(p_{n_{1}}\otimes p_{n_{2}}\otimes
\cdots p_{n_{i}} \otimes 1)$ be a homogeneous element of degree $0$ in
$C^{\infty}(S_{q}^{2\ell+1})$. Then 
\[
Trace(a| D_\ell |^{-z})=2 \sum_{i=0}^{\ell }\sum_{\substack{n\in\bbn^i,\,t\in\bbn\\m\in\bbn^{\ell-i}}}
\frac{\lambda_{n}^{i}}{(| n | + | m | + t)^{z}}+
\sum_{i=0}^{\ell }\sum_{\substack{n\in\bbn^i\\m\in\bbn^{\ell-i}}} 
   \frac{\lambda_{n}^{i}}{(|n | + | m |)^{z}}. 
\]
Now $\sum_{n \in \mathbb{N}^{\ell }}
\frac{\lambda_{n}^{\ell }}{(| n |)^{z}}$ is holomorphic and hence modulo a
holomorphic function
\[
Trace(a| D_\ell |^{-z}) =  2 \sum_{i=0}^{\ell }\left(
\sum_{\substack{n\in\bbn^i,\,t\in\bbn\\m\in\bbn^{\ell-i}}} 
\frac{\lambda_{n}^{i}}{(| n | + | m | +t)^{z}}\right) 
+\sum_{i=0}^{\ell -1}\left(\sum_{\substack{n\in\bbn^i\\m\in\bbn^{\ell-i}}}
\frac{\lambda_{n}^{i}}{(| n | + | m |)^{z}}\right)
\]
It follows from lemma~\ref{dimension} that $Trace(a| D_\ell |^{-z})$ is
meromorphic with simple poles in the set $\{1,2,\cdots, \ell +1\}$. Similarly one can show
that $Trace(aF_\ell| D_\ell |^{-z})$ is meromorphic with simple poles in
$\{1,2,\cdots,\ell\}$. Fix $0 \leq i \leq\ell$. Let $(\lambda_{n})_{n \in
\mathbb{N}^{i}}$ be such that $\sum_{n}\lambda_{n}=1$. Let $a=\sum_{n \in
\mathbb{N}^{i}} \lambda_{n}(p_{n_{1}}\otimes p_{n_{2}}\otimes p_{n_{i}} \otimes
1)$. Then one has $Res_{z=\ell+1-i}Trace(a| D_\ell |^{-z}) = \frac{2}{(\ell-i)!}$ by lemma
\ref{dimension} and by the above equation. Hence every $k \in
\{1,2,\cdots, \ell +1\}$ is in the dimension spectrum. This completes the proof.
\qed

%%%%%%%%%%%%%%%%%%%%%%%%%%%%%%%%%%%%%%%%%%%%%%%%%%%%%%%%%%%%%%%%%%%
%%%%%%%%%%%%%%%%%%%%%%%%%%%%%%%%%%%%%%%%%%%%%%%%%%%%%%%%%%%%%%%%%%%
%%%%%%%%%%%%%%%%%%%%%%%%%%%%%%%%%%%%%%%%%%%%%%%%%%%%%%%%%%%%%%%%%%%

%%%%%%%%%%%%%%%%%%%%%%%%%%%%%%%%%%%%%%%%%%%%%%%%%%%%%%%%%%%%%%%%%%%
%%%%%%%%%%%%%%%%%%%%%%%%%%%%%%%%%%%%%%%%%%%%%%%%%%%%%%%%%%%%%%%%%%%
\section{$SU_q(\ell+1)$ Equivariant spectral triple}
%%%%%%%%%%%%%%%%%%%%%%%%%%%%%%%%%%%%%%%%%%%%%%%%%%%%%%%%%%%%%%%%%%%

%%%%%%%%%%%%%%%%%%%%%%%%%%%%%%%%%%%%%%%%%%%%%%%%%%%%%%%%%%%%%%%%%%%
\subsection{Left multiplication operators}% on $L_2(SU_q(\ell+1))$}
%%%%%%%%%%%%%%%%%%%%%%%%%%%%%%%%%%%%%%%%%%%%%%%%%%%%%%%%%%%%%%%%%%%
Let us recall from~\cite{cha-pal-2008a} some basic facts
on representations of $C(SU_q(\ell+1))$ on $L_2(SU_q(\ell+1))$
by left multiplication.  Irreducible unitary representations of the quantum group
$SU_q(\ell+1)$ are indexed by Young tableaux $\lambda=(\lambda_1,\ldots,\lambda_{\ell+1})$
where $\lambda_i\in\bbn$ and $\lambda_1\geq\lambda_2\geq\cdots\geq\lambda_{\ell+1}=0$
(Theorem~1.5, \cite{wor-1998a}).
Denote by $u^\lambda$ the irreducible unitary indexed by $\lambda$.
Basis elements of the Hilbert space $\clh_\lambda$ on which
$u^\lambda$ acts can be parametrized by arrays of the form
\[
 \bldr=\left(\begin{matrix}r_{11}&r_{12} &\cdots&r_{1,\ell}&r_{1,\ell+1}\cr
                      r_{21}&r_{22}&\cdots &r_{2,\ell}&\cr
                          &\cdots&&&\cr
                      r_{\ell,1}&r_{\ell,2}&&&\cr
                      r_{\ell+1,1}&&&&
                   \end{matrix}\right),
\]
 where $r_{ij}$'s are integers satisfying
$r_{1j}=\lambda_j$ for $j=1,\ldots,\ell+1$,
$r_{ij}\geq r_{i+1,j}\geq r_{i,j+1}\geq 0$ for all $i$, $j$,
and the top row coincides with $\lambda$.
These are known as Gelfand-Tsetlin tableaux, to be abreviated
as GT tableaux now onwards.
Let $\{e(\lambda,\bldr): \bldr\text{ is a GT tableaux with top row }\lambda\}$
be an orthonormal basis for $\clh_\lambda$.
Denote the matrix entries of $u^\lambda$ with respect to this basis
by $u^\lambda_{\bldr,\blds}$. Note that the generators $u_{ij}$ of the
$C^*$-algebra $C(SU_q(\ell+1))$ are the matrix entries of the irreducible
$\one=(1,0,\ldots,0)$.
The collection $\{u^\lambda_{\bldr,\blds}:\lambda,\bldr,\blds\}$
form a complete orthogonal set of vectors in $L_2(SU_q(\ell+1))$.
Denote by $e^\lambda_{\bldr,\blds}$, or by $e_{\bldr,\blds}$ for short
(as $\bldr$ and $\blds$ specify $\lambda$), the normalized $u^\lambda_{\bldr,\blds}$'s,
i.e.\ $e_{\bldr,\blds}=\|u^\lambda_{\bldr,\blds}\|^{-1}u^\lambda_{\bldr,\blds}$.
Then $\{e_{\bldr,\blds}:\bldr,\blds\}$ form a complete orthonormal basis
for $L_2(SU_q(\ell+1))$.

Let $\rho$ be the half-sum of positive roots of $sl(\ell+1)$
and $\lambda(\bldr)$ is the weight of the weight vector
$e(\lambda,\bldr)$. Let $F_\lambda$ be the unique intertwiner
in $\mbox{Mor}\,(u^\lambda, (u^\lambda)^{cc})$ with 
$\mbox{trace}\, F_\lambda= \mbox{trace}\, F_\lambda^{-1}$ 
(here for a representation $u$, its contragradient representation
is denoted by $u^c$; see~\cite{kli-sch-1997a} for details).
Then one has
$\|u^\lambda_{\bldr\blds}\|=d_\lambda^{-\half}q^{-\psi(\bldr)}$,
where
\be\label{eq:kappa0}
\psi(\bldr)=(\rho,\lambda(\bldr))=-\frac{\ell}{2}\sum_{j=1}^{\ell+1}r_{1j}
   + \sum_{i=2}^{\ell+1}\sum_{j=1}^{\ell+2-i}r_{ij},
   \qquad
d_\lambda=\mbox{trace}\, F_\lambda=\sum_{\bldr:\bldr_1=\lambda} q^{2\psi(\bldr)}.
\ee
Write
\be
\kappa(\bldr,\bldm)=
 d_\lambda^\half d_\mu^{-\half}q^{\psi(\bldr)-\psi(\bldm)}.
\ee
From equation~(4.19) in~\cite{cha-pal-2008a}, we have
\be\label{left_mult}
\pi(u_{ij})e^\lambda_{\bldr\blds}
= \sum_{\mu,\bldm,\bldn}
  C_q(\one,\lambda,\mu;i,\bldr,\bldm)C_q(\one,\lambda,\mu;j,\blds,\bldn)
   \kappa(\bldr,\bldm)    e^\mu_{\bldm\bldn},
\ee
where $C_q$ denote the Clebsch Gordon coefficients.

For our subsequent analysis, we will compute the 
quantities $C_q(i,\bldr,\blds)$ and $\kappa(\bldr,\bldm)$
appearing in the above formula.  We will
use the formulae given in (\cite{kli-sch-1997a}, pp.\ 220), keeping in mind
that for our case (i.e.\ for $SU_q(\ell+1)$), the top right entry of
the GT tableaux is zero.

For a positive integer $j$ with $1\leq j\leq \ell+1$,
let 
\be
\bbm_j:=\{(m_1,m_2,\ldots,m_j)\in\bbn^j: 1\leq m_i\leq \ell+2-i\mbox{ for }1\leq
i\leq j\}
\ee
For $M=(m_1,m_2,\ldots,m_i)\in\bbm_i$, denote by $M(\bldr)$ the tableaux 
$\blds$ defined by
\be\label{movenotation}
s_{jk}=\begin{cases}
          r_{jk}+1 & \mbox{ if }k=m_j, 1\leq j\leq i,\cr
         r_{jk} & \mbox{ otherwise}.
       \end{cases}
\ee
With this notation, observe now that
$C_q(i,\bldr,\blds)$ will be zero unless $\blds$ is
$M(\bldr)$ for some $M\in\bbm_i$.
(One has to keep in mind though that not all tableaux of the form $M(\bldr)$
is a valid GT tableaux)

From (\cite{kli-sch-1997a}, pp.\ 220), we have
\be\label{cgc1}
C_q(i,\bldr,M(\bldr))=
\prod_{a=1}^{i-1}R(\bldr,a,m_a,m_{a+1}) \times R'(\bldr,i,m_i)
\ee
where 
% $e_k$ stands for a vector (in the appropriate space) whose
% $k$\raisebox{.4ex}{th} coordinate is 1 and the rest are all zero, and
\bea
R(\bldr,a,j,k) & =& \sgn(k-j)\,q^{\half\left(-r_{aj}+r_{a+1,k} -
k+j\right)}\nonumber \\
&\times&
\left(\prod_{{i=1}\atop{i\neq j}}^{\ell+2-a}
   \frac{[r_{a,i}-r_{a+1,k}-i+k]_q  }{[r_{a,i}-r_{a,j}-i+j]_q}
  \prod_{{i=1}\atop{i\neq k}}^{\ell+1-a}
   \frac{[r_{a+1,i}-r_{a,j}-i+j-1]_q 
}{[r_{a+1,i}-r_{a+1,k}-i+k-1]_q}\right)^\half \label{corrected_1}\\
R'(\bldr,a,j)
&=& q^{\half\left(1-j+\sum_{i=1}^{\ell+1-a}r_{a+1,i} -
           \sum_{{i=1}\atop{i\neq j}}^{\ell+2-a}r_{a,i}\right)} \nonumber \\
&&
\times \left(
\frac{\prod_{i=1}^{\ell+1-a}[r_{a+1,i}-r_{aj}-i+j-1]_q  }
{\prod_{{i=1}\atop{i\neq j}}^{\ell+2-a}[r_{a,i}-r_{aj}-i+j]_q  }\right)^\half,
\label{corrected_2}
\eea
where for an integer $n$, $[n]_q$ denotes the $q$-number
$(q^n-q^{-n})/(q-q^{-1})$
and $\sgn(k-j)$ is 1 if $k\geq j$ and is $-1$ if $k<j$.

\brmrk
Let us look at the denominators in the above expressions.
The integers $r_{a,i}-r_{a,j}$ and $j-i$ are of the same sign.
Therefore for $i\neq j$, the quantity $r_{a,i}-r_{a,j}-i+j$
is nonzero.
Similarly $r_{a+1,i}-r_{a+1,k}$ and $k-i$ are of the same sign.
So if $i\neq k$, then $r_{a+1,i}-r_{a+1,k}-i+k-1$ can be zero
only when $r_{a+1,i}=r_{a+1,k}$ and $k=i+1$.
Now if $\bldr$ and $M(\bldr)$ are GT tableaux, then
$M(\bldr)_{a+1,m_{a+1}}=r_{a+1,m_{a+1}}+1$ and 
$M(\bldr)_{a+1,i}=r_{a+1,i}$ for $i\neq m_{a+1}$.
Therefore if $m_{a+1}=i+1$, then
$r_{a+1,i}-(r_{a+1,m_{a+1}}+1)\geq 0$, i.e.\ $r_{a+1,i}-r_{a+1,m_{a+1}}\geq 1$.
Hence $r_{a+1,i}-r_{a+1,m_{a+1}}-i+m_{a+1}-1\geq 1$.
In other words, all the $q$-numbers
appearing in the denominator in equation~(\ref{cgc1}) are nonzero.
Thus no problem arises from division by zero.
\ermrk

\brmrk
This is essentially a repetition of remark~4.1 of \cite{cha-pal-2008a}.
The formulae (\ref{corrected_1}) and (\ref{corrected_2})
are obtained from equations~(45) and (46), page 220, \cite{kli-sch-1997a}
by replacing $q$ with $q^{-1}$. Equation~(45) is a special
case of the more general formula (48), page 221, \cite{kli-sch-1997a}.
However, there is a small error in equation~(48) there.
The correct form can be found in equations~(3.1, 3.2a, 3.2b)
in \cite{ali-smi-1994a}. Here we have incorporated 
that correction in equations~(\ref{corrected_1}) and (\ref{corrected_2}).
\ermrk

We next compute the quantities $R(\bldr,a,j,k)$ and $R'(\bldr,a,j)$.

For a positive integer $n$, denote by $Q(n)$ the number $(1-q^{2n})^{1/2}$.
Then for any two integers $m$ and $n$, one has
\[
 \left|\frac{[m]_q}{[n]_q}\right|=q^{-|m|+|n|}\left(\frac{Q(|m|)}{Q(|n|)}
\right)^2.
\]
The next two lemmas are obtained from equations (\ref{corrected_1})
and (\ref{corrected_2}) using the above equality repeatedly and
the fact that $r_{a,i}\geq r_{a+1,i}\geq r_{a,i+1}$ for all $a$ and $i$.
%%%%%%%%%%%%%%%%%%%%%%%%%%%%%%%%%%%%%%%%%%%%%%%%%%%%%%
\blmma\label{lem:comp_row1}
For a GT tableaux $\bldr=(r_{ab})$,
denote by $H_{ab}(\bldr)$ and $V_{ab}(\bldr)$ the following differences:
$H_{ab}(\bldr):=r_{a+1,b}-r_{a,b+1}$ and
$V_{ab}(\bldr):=r_{ab}-r_{a+1,b}$. Then one has
\be\label{eq:comp_row1}
R(\bldr,a,j,k)=\sgn(k-j) q^{P(\bldr,a,j,k) + S(\bldr,a,j,k) } L(\bldr,a,j,k),
\ee
where
\bea
 P(\bldr,a,j,k) &=& \sum_{j\wedge k\leq i < j\vee k}
H_{ai}(\bldr)+2\sum_{k<i<j}V_{ai}(\bldr),\\
S(\bldr,a,j,k) &=& \begin{cases}2(j-k-1)+1 & \mbox{ if }j>k,\cr
                          0 & \mbox{ if }j\leq k,\end{cases}\\
L(\bldr,a,j,k) &=& \prod_{{i=1}\atop{i\neq j}}^{\ell+2-a}
   \frac{Q(|r_{a,i}-r_{a+1,k}-i+k|) }{Q(|r_{a,i}-r_{a,j}-i+j|)}
  \prod_{{i=1}\atop{i\neq k}}^{\ell+1-a}
   \frac{Q(|r_{a+1,i}-r_{a,j}-i+j-1|)  }{Q(|r_{a+1,i}-r_{a+1,k}-i+k-1|)}.\nonumber\\
  && \label{eq:comp_row1c}
\eea
\elmma
%%%%%%%%%%%%%%%%%%%%%%%%%%%%%%%%%%%%%%%%%%%%%%%%%%%%%%

%%%%%%%%%%%%%%%%%%%%%%%%%%%%%%%%%%%%%%%%%%%%%%%
\blmma\label{lem:comp_row2}
One has
\be\label{eq:comp_row2}
R'(\bldr,a,j)= q^{P'(\bldr,a,j)} L'(\bldr,a,j),
\ee
where
\bea
 P'(\bldr,a,j) &=& \sum_{j\leq i < \ell+2-a} H_{ai}(\bldr),\label{eq:comp_row2a}\\
L'(\bldr,a,j) &=& 
   \frac{\prod_{i=1}^{\ell+1-a}Q(|r_{a+1,i}-r_{a,j}-i+j-1|) }
          {\prod_{{i=1}\atop{i\neq j}}^{\ell+2-a}Q(|r_{a,i}-r_{a,j}-i+j|)}.\label{eq:comp_row2b}
\eea
\elmma
%%%%%%%%%%%%%%%%%%%%%%%%%%%%%%%%%%%%%%%%%%%%%%%

Combining lemmas~\ref{lem:comp_row1} and \ref{lem:comp_row2}, we get
the following expression for the CG coefficient $C_q(i,\bldr,M)$.

%%%%%%%%%%%%%%%%%%%%%%%%%%%%%%%%%%%%%%%%%%%%%%%
\blmma\label{lem:cgcoeff}
For a move $M\in\bbm_i$, let $\sgn(M)$ denote
the product $\prod_{a=1}^{i-1}\sgn(m_{a+1}-m_a)$. Then one has
\be\label{eq:cgcoeff}
C_q(i,\bldr,M)=\sgn(M)
q^{B(M)+C(\bldr,M)}\left(\prod_{a=1}^{i-1}L(\bldr,a,m_a,m_{a+1})\right)L'(\bldr,
i,m_i),
\ee
where
\bea
B(M) &=& \sum_{j:m_j>m_{j+1}} \left(2(m_j - m_{j+1}-1)+1\right),\\
C(\bldr,M) &=& \sum_{a=1}^{i-1}\left(
   \sum_{m_a\wedge m_{a+1} \leq b < m_a\vee m_{a+1}}H_{ab}(\bldr)
   +2 \sum_{m_{a+1} < b < m_a}V_{ab}(\bldr)\right)
+\sum_{m_i \leq b < \ell+2-i}H_{ib}(\bldr)\nonumber\\
&&
\eea
\elmma
%%%%%%%%%%%%%%%%%%%%%%%%%%%%%%%%%%%%%%%%%%%%%%%

%%%%%%%%%%%%%%%%%%%%%%%%%%%%%%%%%%%%%%%%%%%%%%%
\blmma
\[
\left(\prod_{a=1}^{i-1}L(\bldr,a,m_a,m_{a+1})\right)L'(\bldr,i,m_i)=1+o(q).
\]
\elmma
%%%%%%%%%%%%%%%%%%%%%%%%%%%%%%%%%%%%%%%%%%%%%%%
\prf
This is a consequence of the following two  inequalities:
\[
 |1-(1-x)^\half|<x \quad \text{ for } 0\leq x \leq 1,
\]
and for $0<r<1$,
\[
 |1-(1-x)^{-\half}|<cx \quad \text{ for } 0\leq x \leq r,
\]
where $c$ is some fixed constant that depends on $r$.
\qed

Next we come to the computation of $\kappa(\bldr,\bldm)$.
Since $C_q(i,\bldr,\bldm)$ is 0 unless $\bldm$ is of the form
$M(\bldr)$ for some move $M=(m_1,\ldots,m_i)$, we need only to compute
$\kappa(\bldr,M(\bldr))$
which we will denote by $\kappa(\bldr,M)$.

Since
\[
\psi(\blds)
=-\frac{\ell}{2}\sum_{j=1}^{\ell+1}s_{1j}+
   \sum_{i=2}^{\ell+1}\sum_{j=1}^{\ell+2-i}s_{ij},
\]
we have
\[
q^{\psi(\bldr)-\psi(M(\bldr))}=
 q^{-\frac{\ell}{2}\sum_{j=1}^{\ell+1}r_{1j}+
   \sum_{i=2}^{\ell+1}\sum_{j=1}^{\ell+2-i}r_{ij}
    +\frac{\ell}{2}(\sum_{j=1}^{\ell+1}r_{1j}+1)
   - (\sum_{i=2}^{\ell+1}\sum_{j=1}^{\ell+2-i}r_{ij}+i-1) }
   =q^{\frac{\ell}{2}-i+1}.
\]
Let $\lambda=(\lambda_1,\ldots,\lambda_{\ell},0)$ be the top row of $\bldr$.
Then
\[
\min\{\psi(\blds):\blds_1=\lambda\}
  = -\frac{\ell}{2}\sum_1^\ell\lambda_i + \sum_{k=2}^\ell(k-1)\lambda_k.
\]
Hence
\[
d_\lambda=\sum_{\blds:\blds_1=\lambda}
q^{2\psi(\blds)}=q^{-\ell\sum_1^\ell\lambda_i + 2\sum_{k=2}^\ell(k-1)\lambda_k}
(1+q^2\phi(q^2)),
\]
where $\phi$ is a polynomial. Therefore 
\[
d_\lambda=q^{-\ell\sum_1^\ell\lambda_i +
2\sum_{k=2}^\ell(k-1)\lambda_k}(1+o(q)).
\]
It follows that
\[
\left(\frac{d_\lambda}{d_{\lambda+e_{m_1}}}\right)^\half=
q^{\frac{\ell}{2}-m_1+1}(1+o(q)).
\]
Thus
\be\label{eq:kappa1}
\kappa(\bldr,M(\bldr))=q^{\ell+2-i-m_1}(1+o(q)).
\ee

Next, observe that
\bean
\lefteqn{B(M)+\ell+2-i-m_1}\\
 &=& \sum_{j:m_j>m_{j+1}} \left(2(m_j - m_{j+1}-1)+1\right)-(m_1-m_i)
+\ell+2-i-m_i\\
&=&2\sum_{j:m_j>m_{j+1}} (m_j - m_{j+1})-\sum_{j:m_j>m_{j+1}}1
                -\sum_{j=1}^{i-1}(m_j-m_{j+1})+\ell+2-i-m_i\\
&=&2\sum_{j:m_j>m_{j+1}} (m_j -
m_{j+1})-\sum_{j=1}^{i-1}(m_j-m_{j+1})-\sum_{j:m_j>m_{j+1}}1
              +\ell+2-i-m_i\\
&=&\sum_{j=1}^{i-1}\left|m_j - m_{j+1}\right|-\#\{1\leq j\leq i-1:m_j>
m_{j+1}\}+\ell+2-i-m_i.
\eean
Thus if we write 
\bea
 A(M) &=& \sum_{j=1}^{i-1}\left|m_j - m_{j+1}\right|-\#\{1\leq j\leq
i-1:m_j> m_{j+1}\},\label{eq:A}\\
K(M) &=& \ell+2-i-m_i,\label{eq:K}
\eea
then both $A(M)$ and $K(M)$ are nonnegative and $B(M)+\ell+2-i-m_1=A(M)+K(M)$.
Thus we have
\bea
\pi(u_{ij})e^\lambda_{\bldr\blds}
&=& \sum_{{M\in\bbm_i}\atop {M'\in\bbm_j}}
C_q(i,\bldr,M(\bldr))\kappa(\bldr,M)C_q(j,\blds,M'(\blds))e_{M(\bldr),M'(\blds)}
\label{eq:left_mult3a}\\
&=&\sum_{{M\in\bbm_i}\atop
{M'\in\bbm_j}}\sgn(M)\sgn(M')q^{A(M)+K(M)+C(\bldr,M)+B(M')+C(\blds,M')}
(1+o(q))e_{M(\bldr),M'(\blds)}. \nonumber\\
&&\label{eq:left_mult3}
\eea

%%%%%%%%%%%%%%%%%%%%%%%%%%%%%%%%%%%%%%%%%%%%%%%%%%%%%%%%
%%%%%%%%%%%%%%%%%%%%%%%%%%%%%%%%%%%%%%%%%%%%%%%%%%%%%%%%

%%%%%%%%%%%%%%%%%%%%%%%%%%%%%%%%%%%%%%%%%%%%%%%%%%%%%%%%
\subsection{The spectral triple}
%%%%%%%%%%%%%%%%%%%%%%%%%%%%%%%%%%%%%%%%%%%%%%%%%%%%%%%%%
%%%%%%%%%%%%%%%%%%%%%%%%%%%%%%%%%%%%%%%%%%%%

Let us briefly recall from \cite{cha-pal-2008a} the description 
of the $L_2$ space of the sphere
sitting inside $L_2(SU_q(\ell+1))$.
%%%%%%%%%%%%%%%%%%%%%%%%%%%%%%%%%%%%%%%%%%%%
Let $u^\one$ denote the fundamental unitary for $SU_q(\ell+1)$,
i.\ e.\ the irreducible unitary representation corresponding to the
Young tableaux $\one=(1,0,\ldots,0)$.
Similarly write $v^\one$ for the fundamental unitary for $SU_q(\ell)$.
Fix some bases for the corresponding representation spaces.
Then $C(SU_q(\ell+1))$ is the $C^*$-algebra generated by the matrix
entries $\{u^\one_{ij}\}$ and $C(SU_q(\ell))$ is the
$C^*$-algebra generated by the matrix
entries $\{v^\one_{ij}\}$.
Now define $\phi$ by
\be
\phi(u^\one_{ij})=\begin{cases} I & \mbox{if $i=j=1$},\cr
        v^\one_{i-1,j-1} & \mbox{if $2\leq i,j\leq \ell+1$},\cr
        0 & \mbox{otherwise.}
       \end{cases}
\ee
Then $C(SU_q(\ell+1)\verb1\1SU_q(\ell))$ is the $C^*$-subalgebra of
$C(SU_q(\ell+1))$ generated by the entries $u_{1,j}$ for $1\leq j\leq \ell+1$.
Define $\psi:C(S_q^{2\ell+1})\rightarrow C(SU_q(\ell+1)\verb1\1SU_q(\ell))$
by 
\[
 \psi(z_i)=q^{-i+1}u^*_{1,i}.
\]
This gives an isomorphism between $C(SU_q(\ell+1)\verb1\1SU_q(\ell))$ and
$C(S_q^{2\ell+1})$
and the following diagram commutes:%\\[1ex]
\[
\def\labelstyle{\scriptstyle}
  \xymatrix@C=23pt@R=50pt{
 C(S_q^{2\ell+1})\ar[d]_{\psi}\ar[r]_-{\tau} & C(S_q^{2\ell+1})\otimes
C(SU_q(\ell+1))\ar[d]_{\psi\otimes\id}&\\
C(SU_q(\ell+1)\backslash SU_q(\ell))\ar[r]_-{\Delta}  & C(SU_q(\ell+1)\backslash SU_q(\ell))\otimes C(SU_q(\ell+1))\\
}
\]
where $\tau :C(S_q^{2\ell+1})\rightarrow C(S_q^{2\ell+1})\otimes C(SU_q(\ell+1))$
is the homomorphism given by $\tau(z_i)=\sum_k z_k\otimes u_{k,i}^*$
that gives an action of the quantum group $SU_q(\ell+1)$
on $S_q^{2\ell+1}$.
In other words,
$(C(S_q^{2\ell+1}), SU_q(\ell+1),\tau)$ is the quotient space $SU_q(\ell+1)\verb1\1SU_q(\ell)$.
This choice of $\psi$  makes $L_2(SU_q(\ell+1)\verb1\1SU_q(\ell))$ 
a span of certain rows of the $e_{\bldr,\blds}$'s
as the following two propositions say.

%%%%%%%%%%%%%%%%%%%%%%%%%%%%%%%%%%%%%%%%%%%%%%%%%%%%%%%%%
\bppsn[\cite{cha-pal-2008a}]
Assume $\ell>1$.
The right regular representation $u$ of $G$ keeps the subspace
$L_2(SU_q(\ell+1)\verb1\1SU_q(\ell))$ invariant, and the restriction of $u$ to
$L_2(SU_q(\ell+1)\verb1\1SU_q(\ell))$ decomposes as a direct sum of exactly one copy
of each of the irreducibles given by the young tableaux
$\lambda_{n,k}:=(n+k, k,k,\ldots, k,0)$, with $n,k\in\bbn$.
\eppsn
%%%%%%%%%%%%%%%%%%%%%%%%%%%%%%%%%%%%%%%%%%%%%%%%%%%%%%%%%

%%%%%%%%%%%%%%%%%%%%%%%%%%%%%%%%%%%%%%%%%%%%%%%%%%%%%%%%%
\bppsn[\cite{cha-pal-2008a}]
Let  $\bldr^{nk}$ denote the GT tableaux
given by
\[
r^{nk}_{ij}=\begin{cases} n+k & \mbox{if $i=j=1$},\cr
                0  & \mbox{if $i=1$, $j=\ell+1$},\cr
                k  & \mbox{otherwise},\end{cases}
\]
where $n,k \in \bbn$.
Let $\scrg_0^{n,k}$ be the set of all GT tableaux with
top row $(n+k,k,\ldots,k,0)$.
Then the family of vectors
\[
\{e_{\bldr^{nk},\blds}: n,k\in\bbn,\, \blds\in\scrg_0^{n,k}\}
\]
form a complete
orthonormal basis for $L_2(SU_q(\ell+1)\verb1\1SU_q(\ell))$.
\eppsn
%%%%%%%%%%%%%%%%%%%%%%%%%%%%%%%%%%%%%%%%%%%%%%%%%%%%%%%%%
We will denote $\cup_{n,k}\scrg_0^{n,k}$ by $\scrg_0$.
Since the top row of $\bldr^{nk}$ determines $\bldr^{nk}$ completely
and for $e_{\bldr^{nk},\blds}$, the top row of $\blds$ equals the top row of
$\bldr^{nk}$,
one can index the orthonormal basis $e_{\bldr^{nk},\blds}$ just by
$\blds\in\scrg_0$.
It was shown in \cite{cha-pal-2008a}
that the restriction of the left multiplication to 
$C(SU_q(\ell+1)\verb1\1SU_q(\ell))\cong C(S_q^{2\ell+1})$
keeps $L_2(SU_q(\ell+1)\verb1\1SU_q(\ell))\cong L_2(S_q^{2\ell+1})$ invariant.
We will continue to denote this restriction by $\pi$.
The operators $\pi(z_j)=q^{-j+1}\pi(u_{1,j}^*)$
will be denoted by $Z_{j,q}$. The $C^*$-algebra $\pi(C(S_q^{2\ell+1}))$
will be denoted by $C_\ell$.

The following theorem gives a generic equivariant spectral triple 
for the spheres $S_q^{2\ell+1}$ constructed in~\cite{cha-pal-2008a}.

%%%%%%%%%%%%%%%%%%%%%%%%%%%%%%%%%%%%%%%%%%%%%%%%%%%%%
\bthm[\cite{cha-pal-2008a}]\label{generic_d_sph}
Let $D_{eq}$ be the operator on $L_2(S_q^{2\ell+1})$ given by:
\be\label{eq_sphere1}
D_{eq} e_{\bldr^{nk},\blds}=\begin{cases}
                k e_{\bldr^{nk},\blds}& \mbox{if $n=0$},\cr
                 -(n+k)  e_{\bldr^{nk},\blds}& \mbox{if $n>0$}.
                   \end{cases}
\ee
Then $(\cla(S_{q}^{2\ell+1}),L_2(S_q^{2\ell+1}),D_{eq})$  is an 
equivariant nondegenerate $(2\ell+1)$-summable odd
spectral triple.
\ethm
%%%%%%%%%%%%%%%%%%%%%%%%%%%%%%%%%%%%%%%%%%%%%%%%%%%%%
Our main aim in the rest of the paper is to precisely formulate
the smooth function algebra for this spectral triple, establish
its regularity, and compute the dimension spectrum.

%%%%%%%%%%%%%%%%%%%%%%%%%%%%%%%%%%%%%%%%%%%%%%%%%%%%%%%%
%%%%%%%%%%%%%%%%%%%%%%%%%%%%%%%%%%%%%%%%%%%%%%%%%%%%%%%%

%%%%%%%%%%%%%%%%%%%%%%%%%%%%%%%%%%%%%%%%%%%%%%%%%%%%%%%%%
\subsection{The case $q=0$}
%%%%%%%%%%%%%%%%%%%%%%%%%%%%%%%%%%%%%%%%%%%%%%%%%%%%%%%%%
%%%%%%%%%%%%%%%%%%%%%%%%%%%%%%%%%%%%%%%%%%%%%%%%%%%%%%%%
The $L_2$ spaces $L_2(S_q^{2\ell+1})$ for different values of $q$
can be identified by identifying the elements of their 
canonical orthonormal bases which are parametrized by the same set.
Thus we will assume we are woring with one single Hilbert space
$\clh$ with orthonormal basis given by
$e_{\bldr^{n,k}\blds}$ where $\bldr^{n,k}$ is as defined earlier
and $\blds$ is given by
\be\label{eq:GTnotation}
\blds=\left(\begin{matrix}
  c_1=n+k& k  & k& \cdots& k & k & d_1=0 \cr
  c_2& k & k &\cdots & k & d_2&\cr
   \cdots    & &\cdots&&&\cr
   c_{\ell-1}& k& d_{\ell-1} &&&&\cr
   c_\ell & d_\ell &&&&&\cr
                     c_{\ell+1}=d_{\ell+1} &&&&&&
                   \end{matrix}\right)
\ee
where $c_1\geq c_2\geq\ldots\geq c_\ell\geq k$, $d_1\leq d_2\leq\ldots\leq d_\ell\leq k$
and $d_\ell\leq d_{\ell+1}\leq c_\ell$.
Since specifying the GT tableaux $\blds$ specifies $\bldr^{n,k}$ also
and thus completely specifies the basis element $e_{\bldr^{n,k}\blds}$,
we will sometimes use just $\blds$ in place of the basis element $e_{\bldr^{n,k}\blds}$.

Let us denote by $\bbm_j^\pm$ the following subsets of $\bbm_j$:
\bean
\bbm_j^+ &=& \{(m_1,\ldots,m_j)\in\bbm_j: m_i\in\{1,\ell+2-i\} 
                   \text{ for }1\leq i\leq j,\;m_1=1\},\\
\bbm_j^- &=& \{(m_1,\ldots,m_j)\in\bbm_j: m_i\in\{1,\ell+2-i\} 
                   \text{ for }1\leq i\leq j,\;m_1=\ell+1\}.
\eean
Let us denote by $N_{i,j}$ the following element of $\bbm_j$:
\[
 N_{i,j}=(\underbrace{1,\ldots,1}_{i},\ell+1-i,\ell-i,\ldots,\ell+2-j),\qquad
  0\leq i\leq j\leq \ell+1.
\]
We will denote $N_{i,\ell+1}$ by just $N_i$.
Then from (\ref{eq:left_mult3}), we get
\bea
\pi(u_{1j})e_{\bldr^{n,k}\blds}
&=&\sum_{M\in\bbm_j^+}\sgn(M)q^{\ell+k+B(M)+C(\blds,M)}
(1+o(q))e_{\bldr^{n+1,k},M(\blds)} \nonumber\\
&&  +
\sum_{M\in\bbm_j^-}\sgn(M)q^{B(M)+C(\blds,M)}(1+o(q))e_{
\bldr^{n,k-1},M(\blds)}\label{eq:left_mult4}
\eea
Therefore
\bea
Z_{j,q}^*e_{\bldr^{n,k}\blds}
&=&\sum_{M\in\bbm_j^+}\sgn(M)q^{-j+1+\ell+k+B(M)+C(\blds,M)}
(1+o(q))e_{\bldr^{n+1,k},M(\blds)} \nonumber\\
&&  +
\sum_{M\in\bbm_j^-}\sgn(M)q^{-j+1+B(M)+C(\blds,M)}
(1+o(q))e_{\bldr^{n,k-1},M(\blds)}\label{eq:left_mult5}
\eea

Let us first look at the cases $1\leq j\leq \ell$.
In this case, the power of $q$ in the first summation is
positive. Therefore none of the terms would survive
for $q=0$.
For terms in the second summation, assume $M\in\bbm_j$ with
$m_1=\ell+1$ and $m_i=1$ for some $i\leq j$.
Let $a=\min\{2\leq i\leq j: m_i=1\}$.
Then $m_i=\ell+2-i$ for $1\leq i\leq a-1$ so that
\bean
B(M) &\geq& \sum_{i=1}^{a-2}(2((\ell+2-i)-(\ell+1-i)-1)+1)+2(\ell+3-a-1-1)+1\\
   &=&a-2+2(\ell-a+1)+1\\
 &=& 2\ell-a+1.
\eean
Hence $B(M)+1-j>0$ and so such terms will not survive for $q=0$.
Thus the only term that will survive is the one corresponding
to $M=N_{0,j}=(\ell+1,\ell,\ell-1,\ldots,\ell+2-j)$.
In this case we have
$B(M)=j-1$, $C(\blds,M)=d_j$ and $\sgn(M)=(-1)^{j-1}$.
Therefore
\be\label{eq:q0case1}
Z_{j,0}^*e_{\bldr^{n,k}\blds}=\begin{cases}
(-1)^{j-1}e_{\bldr^{n,k-1},N_{0,j}(\blds)} & \mbox{ if }d_j=0,\cr
    0  &         \mbox{ if }d_j>0.
                              \end{cases}
\ee

Next let us look at the case $j=\ell+1$.
Here the first sum will be over all $M$ with $m_1=1=m_{\ell+1}$.
If $m_i\neq 1$ for some $i$, then $B(M)>0$ and
therefore the power of $q$ will be positive, so that
the term will not survive for $q=0$.
If $m_i=1$ for all $i$, i.e.\ if $M=N_\ell$,
then we have $B(M)=0=C(\blds,M)$ and $\sgn(M)=1$.
Therefore for $q=0$, the first summation will become
$e_{\bldr^{n+1,k},N_{\ell}(\blds)}$ provided $k=0$.

The second sum is over all $M$ with $m_1=\ell+1$.
Let $a=\min\{2\leq i\leq \ell+1: m_i=1\}$.
Then as before, $B(M) \geq 2\ell-a+1$.
Therefore if $a\leq\ell$, then
$-\ell+B(M)\geq\ell-a+1>0$, so that the term will
not survive for $q=0$. If $a=\ell+1$, i.e.\ if
$M=N_0$, then $B(M)=\ell$, $C(\blds,M)=d_{\ell+1}$ and
$\sgn(M)=(-1)^{\ell}$.
So for $q=0$, the second summation will become
$(-1)^\ell e_{\bldr^{n,k-1},N_{0}(\blds)}$ if  $k>0$ and $d_{\ell+1}=0$.
Thus we have
\be\label{eq:q0case2}
Z_{\ell+1,0}^*e_{\bldr^{n,k}\blds}=\begin{cases}
e_{\bldr^{n+1,k},N_{\ell}(\blds)} &  \mbox{ if }k=0,\cr
(-1)^\ell e_{\bldr^{n,k-1},N_{0}(\blds)} &  \mbox{ if } k>0,\; d_{\ell+1}=0,\cr
    0  &         \mbox{ if } k>0, \;d_{\ell+1}>0.
                              \end{cases}
\ee

Next we will establish a natural unitary map between
$L_2(S_q^{2\ell+1})$ and
\[
\clh_{\Sigma}\equiv \underbrace{\ell_2(\bbn)\otimes\cdots\otimes\ell_2(\bbn)}_{
\mbox{$\ell$ copies}}\otimes\ell_2(\bbz)
\otimes\underbrace{\ell_2(\bbn)\otimes\cdots\otimes\ell_2(\bbn)}_{\mbox{$\ell$
copies}}.
\]
For $t\in\bbr$, let  $t_+$ denote the positive part $\max\{t, 0\}$ and 
let $t_-$ denote the negative part $\max\{(-t), 0\}$ of $t$.
Let us now observe that for any $\gamma\in\Gamma_\Sigma$,
the tableaux
\[
\blds(\gamma):=\left(\begin{matrix}
  \sum_{1}^{2\ell+1}|\gamma_i|& \sum_{1}^\ell \gamma_i +(\gamma_{\ell+1})_+ 
                      & \cdots
                      & \sum_{1}^\ell \gamma_i +(\gamma_{\ell+1})_+ & 0 
                               \cr
  \sum_{1}^{2\ell} |\gamma_i| & \sum_{1}^\ell \gamma_i +(\gamma_{\ell+1})_+ &\cdots &
 \gamma_1&\cr
                 \cdots         &\cdots&&&\cr
    \sum_{1}^{\ell+3} |\gamma_i|& \sum_{1}^\ell \gamma_i +(\gamma_{\ell+1})_+&
\sum_{1}^{\ell-2} \gamma_i&&&\cr
   \sum_{1}^{\ell+2} |\gamma_i| & \sum_{1}^{\ell-1} \gamma_i&&&&\cr
                      \sum_{1}^\ell \gamma_i  +(\gamma_{\ell+1})_{-}  &&&&&
                   \end{matrix}\right)
\]
is in $\scrg_0$.
Conversely, let $\blds\in\scrg_0^{n,k}$ for some $n,k\in\bbn$
so that $e_{\bldr^{n,k}\blds}$ is a basis element of $L_2(S_q^{2\ell+1})$.
Note that $\blds$ is of the form (\ref{eq:GTnotation}).
Define $\gamma\in\Gamma_\Sigma$ as follows:
\begin{enumerate}
\item 
if $k> d_{\ell+1}$, then
\[
\begin{array}{ll}
\gamma_i=d_{i+1}-d_i &\text{ for } 1\leq i\leq \ell-1,\\
\gamma_i=c_{2\ell+2-i}-c_{2\ell+3-i} &\text{ for } \ell+3\leq i\leq 2\ell+1,\\
\gamma_\ell=d_{\ell+1}-d_\ell, \quad \gamma_{\ell+1}=k-d_{\ell+1},\quad  \gamma_{\ell+2}=c_\ell-k,&
\end{array}
\]
\item 
if $k\leq d_{\ell+1}$, then
\[
\begin{array}{ll}
\gamma_i=d_{i+1}-d_i &\text{ for }  1\leq i\leq \ell-1,\\
\gamma_i=c_{2\ell+2-i}-c_{2\ell+3-i} &\text{ for }  \ell+3\leq i\leq 2\ell+1,\\
\gamma_\ell=k-d_\ell, \quad \gamma_{\ell+1}=k-d_{\ell+1}, \quad \gamma_{\ell+2}=c_\ell-d_{\ell+1}.&
\end{array}
\]
\end{enumerate}
Then $\blds(\gamma)=\blds$. Thus we have a bijective correspondence
between $\scrg_0$ and $\Gamma_\Sigma$. We will often denote
a basis element $e_{\bldr^{n,k}\blds}$ by $\xi_\gamma$
using this bijective correspomdence.

%%%%%%%%%%%%%%%%%%%%%%%%%%%%%%%%%%%%%%%%%%%%%%%%%%%
\blmma\label{lem:newonb}
Let $\gamma\in\Gamma_\Sigma$. For $n\in\bbz$, let
\[
 Z_{\ell+1,0}^{(n)}:=\begin{cases}
                      Z_{\ell+1,0}^n & \mbox{ if } n\geq 0,\cr
                   (Z_{\ell+1,0}^*)^{-n} & \mbox{ if } n< 0.
                     \end{cases}
\]
Define 
\[
 \xi'_{\gamma}:=
    Z_{1,0}^{\gamma_1}\ldots Z_{\ell,0}^{\gamma_\ell}Z_{\ell+1,0}^{(\gamma_{\ell+1})}
         \left(\begin{matrix}\sum_{\ell+2}^{2\ell+1} \gamma_i&0 &\cdots& 0 & 0 \cr
                      \sum_{\ell+2}^{2\ell} \gamma_i& 0 &\cdots & 0&\cr
                          &\cdots&&&\cr
                      \gamma_{\ell+2} & 0&&&\cr
                      0&&&&
                   \end{matrix}\right).
\]
Then $\{\xi'_{\gamma}: \gamma\in\Gamma_\Sigma\}$
is an orthonormal basis for $L_2(S_q^{2\ell+1})$.
\elmma
%%%%%%%%%%%%%%%%%%%%%%%%%%%%%%%%%%%%%%%%%%%%%%%%%%%
\prf
It follows from equations~(\ref{eq:q0case1}) and (\ref{eq:q0case2})
that the actions of $Z_{j,0}$ for $1\leq j\leq\ell$
on the basis elements $e_{\bldr^{n,k}\blds}$
are as follows:
\bean
\lefteqn{Z_{j,0}:\left(\begin{matrix}
  n+k& k  & & \cdots&& k & k & 0 \cr
  c_2& k &  &\cdots && k & 0&\cr
   \cdots    & &\cdots&&&&&&\cr
  c_j &k&\cdots& k & 0 &&&&\cr
  c_{j+1}&k&\cdots & d_{j+1} &&&&&\cr
   \cdots    & &\cdots&&&&&&\cr
   c_\ell & d_\ell &&&&&&&\cr
      d_{\ell+1}& &&&&&&&
                   \end{matrix}\right)
\longrightarrow}\\
&& (-1)^{j-1}\left(\begin{matrix}
  1+n+k& 1+k  & & \cdots&& 1+k & 1+k & 0 \cr
  1+c_2& 1+k &  &\cdots && 1+k & 0&\cr
   \cdots    & &\cdots&&&&&&\cr
  1+c_j &1+k&\cdots& 1+k & 0 &&&&\cr
  1+c_{j+1}&1+k&\cdots & 1+d_{j+1} &&&&&\cr
   \cdots    & &\cdots&&&&&&\cr
   1+c_\ell & 1+d_\ell &&&&&&&\cr
     1+ d_{\ell+1}& &&&&&&&
                   \end{matrix}\right)
\eean
and is 0 for $\blds$ with $d_j>0$.

Similarly the action of $Z_{\ell+1,0}$ on the basis elements
are as follows:
\begin{align*}
\left(\begin{matrix}
  n& 0  & & \cdots&&  0 & 0 \cr
  c_2& 0 &  &\cdots &  & 0&\cr
   \cdots    & &\cdots&&&&&\cr
  c_{\ell-1}&0& 0 & &&&&\cr
   c_\ell & 0 &&&&&&&\cr
      d_{\ell+1}& &&&&&&&
                   \end{matrix}\right)
&\longrightarrow
 \left(\begin{matrix}
  n-1& 0  & & \cdots&&  0 & 0 \cr
  c_2-1& 0 &  &\cdots &  & 0&\cr
   \cdots    & &\cdots&&&&&\cr
  c_{\ell-1}-1&0& 0 & &&&&\cr
   c_\ell -1 & 0 &&&&&&&\cr
      d_{\ell+1}-1& &&&&&&&
                   \end{matrix}\right)
\\
\intertext{if $d_{\ell+1}>0$, and}
\left(\begin{matrix}
  n+k& k  & & \cdots&&  k & 0 \cr
  c_2& k &  &\cdots &  & 0&\cr
   \cdots    & &\cdots&&&&\cr
  c_{\ell-1}&k& 0 & &&&\cr
   c_\ell & 0 &&&&&&\cr
      0& &&&&&&
                   \end{matrix}\right)
&\longrightarrow
 (-1)^{\ell}\left(\begin{matrix}
  1+n+k& 1+k  & & \cdots&&  1+k & 0 \cr
  1+c_2& 1+k &  &\cdots &  & 0&\cr
   \cdots    & &\cdots&&&&\cr
  1+c_{\ell-1}&1+k& 0 & &&&\cr
   1+c_\ell & 0 &&&&&&\cr
      0& &&&&&&
                   \end{matrix}\right)
\\
\intertext{if $d_{\ell+1}=0$.
Similarly the action of $Z_{\ell+1,0}^*$ on the basis elements
are as follows:}
\left(\begin{matrix}
  n& 0  & & \cdots&&  0 & 0 \cr
  c_2& 0 &  &\cdots &  & 0&\cr
   \cdots    & &\cdots&&&&&\cr
  c_{\ell-1}&0& 0 & &&&&\cr
   c_\ell & 0 &&&&&&&\cr
      d_{\ell+1}& &&&&&&&
                   \end{matrix}\right)
&\longrightarrow
 \left(\begin{matrix}
  1+n& 0  & & \cdots&&  0 & 0 \cr
  1+c_2& 0 &  &\cdots &  & 0&\cr
   \cdots    & &\cdots&&&&&\cr
  1+c_{\ell-1}&0& 0 & &&&&\cr
   1+c_\ell  & 0 &&&&&&&\cr
      1+d_{\ell+1}& &&&&&&&
                   \end{matrix}\right)
\\
\intertext{and for $k>0$,}
\left(\begin{matrix}
  n+k& k  & & \cdots&&  k & 0 \cr
  c_2& k &  &\cdots &  & 0&\cr
   \cdots    & &\cdots&&&&\cr
  c_{\ell-1}&k& 0 & &&&\cr
   c_\ell & 0 &&&&&&\cr
      0& &&&&&&
                   \end{matrix}\right)
&\longrightarrow
 (-1)^{\ell}\left(\begin{matrix}
  n+k-1& k-1  & & \cdots&&  k-1 & 0 \cr
  c_2-1& k -1&  &\cdots &  & 0&\cr
   \cdots    & &\cdots&&&&\cr
  c_{\ell-1}-1&k-1& 0 & &&&\cr
  c_\ell-1 & 0 &&&&&&\cr
      0& &&&&&&
                   \end{matrix}\right)
\end{align*}

Then it follows from the above
that
\bea
\lefteqn{ Z_{1,0}^{\gamma_1}\ldots Z_{\ell,0}^{\gamma_\ell}Z_{\ell+1,0}^{(\gamma_{\ell+1})}
         \left(\begin{matrix}\sum_{\ell+2}^{2\ell+1} \gamma_i&0 &\cdots& 0 & 0 \cr
                      \sum_{\ell+2}^{2\ell} \gamma_i& 0 &\cdots & 0&\cr
                          &\cdots&&&\cr
                      \gamma_{\ell+2} & 0&&&\cr
                      0&&&&
                   \end{matrix}\right)}\nonumber \\
&=& (-1)^{\eta(\gamma)}\left(\begin{matrix}
  \sum_{1}^{2\ell+1}|\gamma_i|& \sum_{1}^\ell \gamma_i +(\gamma_{\ell+1})_+ 
                      & \cdots
                      & \sum_{1}^\ell \gamma_i +(\gamma_{\ell+1})_+ & 0 
                               \cr
  \sum_{1}^{2\ell} |\gamma_i| & \sum_{1}^\ell \gamma_i +(\gamma_{\ell+1})_+ &\cdots &
 \gamma_1&\cr
                 \cdots         &\cdots&&&\cr
    \sum_{1}^{\ell+3} |\gamma_i|& \sum_{1}^\ell \gamma_i +(\gamma_{\ell+1})_+&
\sum_{1}^{\ell-2} \gamma_i&&&\cr
   \sum_{1}^{\ell+2} |\gamma_i| & \sum_{1}^{\ell-1} \gamma_i&&&&\cr
                      \sum_{1}^\ell \gamma_i  +(\gamma_{\ell+1})_{-}  &&&&&
                   \end{matrix}\right),\nonumber \\
&&\label{eq:newonb2}
\eea
where
$\eta(\gamma):=\sum_{i=1}^{\ell}(i-1)\gamma_i +\ell (\gamma_{\ell+1})_{+}$.
Thus $\xi'_\gamma=(-1)^{\eta(\gamma)}\xi_\gamma$.
Therefore it follows that $\{\xi'_\gamma:\gamma\in\Gamma_\Sigma\}$
is an orthonormal basis for $L_2(S_q^{2\ell+1})$.
\qed

The map $U:L_2(S_q^{2\ell+1})\rightarrow \clh_\Sigma$ 
given by
$U\xi'_\gamma=e_\gamma$
sets up a unitary isomorphism between $L_2(S_q^{2\ell+1})$ and
$\clh_\Sigma$.
Let $P$ denote the projection onto the span of $e_0\otimes\cdots\otimes e_0$ in
$\ell_2(\bbn^\ell)$. Then we have
\be\label{eq:decomposition1}
 UZ_{j,0}U^*=Y_{j,0}\otimes I= Y_{j,0}\otimes P + Y_{j,0}\otimes (I-P),
\ee
and
\be\label{eq:decomposition2}
 UD_{eq}U^*=D_{\ell}\otimes P -|D_{\ell}|\otimes (I-P) - I\otimes\widetilde{N},
\ee
where 
$\widetilde{N}$ is the operator 
$e_{m_1}\otimes\cdots\otimes e_{m_\ell}\mapsto \left(\sum m_i\right)
e_{m_1}\otimes\cdots\otimes e_{m_\ell}$.
In other words,
with respect to the decomposition 
\[
\clh_{\Sigma}=\clh_\ell \oplus \left(\clh_\ell\otimes
\ell_2(\bbn^\ell\backslash\{0,\ldots,0\})\right),
\]
one has
\[
 UZ_{j,0}U^*=Y_{j,0}\oplus (Y_{j,0}\otimes I),
\]
and
\[
 UD_{eq}U^*=D_{\ell}\oplus \left(-|D_{\ell}|\otimes I - I\otimes
\widetilde{N}\right).
\]

Next we will define the smooth function algebra $C^\infty_{eq}(S_0^{2\ell+1})$
and prove that the spectral triple $( C^\infty_{eq}(S_0^{2\ell+1}),\clh,D_{eq})$
is regular with simple dimension spectrum $\{1,2,\ldots,2\ell+1\}$.

It follows from decomposition~(\ref{eq:decomposition1}) that
if we identify $L_2(S_q^{2\ell+1})$ with $\clh_{\Sigma}$,
then the $C^*$-algebra generated by the $Z_{j,0}$'s
is $A_\ell\otimes I$, where $A_\ell$ is the $C^*$-algebra generated
by the $Y_j$'s in $\cll(\clh_{\ell})$.
% In view of the decomposition~(\ref{eq:decomposition1}--\ref{eq:decomposition2}),
Therefore it is natural to define
\be\label{eq:smoothq0}
C^\infty_{eq}(S_0^{2\ell+1})=\{a\otimes I: a\in\cla_\ell^\infty\}.
\ee

%%%%%%%%%%%%%%%%%%%%%%%%%%%%%%%%%%%%%%%%%%%%%%
\bthm\label{thm:regularityq=0}
The triple $(C^\infty_{eq}(S_0^{2\ell+1}),\clh_{\Sigma},D_{eq})$
is a regular spectral triple with simple dimension spectrum
$\{1,2,\ldots,2\ell+1\}$.
\ethm
%%%%%%%%%%%%%%%%%%%%%%%%%%%%%%%%%%%%%%%%%%%%%%
\prf
Since $\cla^\infty_\ell$ is closed under holomorphic function calculus
in $A_\ell$, it follows that  $C^\infty_{eq}(S_0^{2\ell+1})$
is closed under holomorphic function calculus in 
$C^*(\{Z_{j,0}:1\leq j\leq \ell+1\})=A_\ell\otimes I$.
In order to show regularity, let us introduce the algebra
%%%%%%%%%%%%%%%%%%%%%%%%%%%%%%%%%%%%%%%%%%%%%%
\be\label{eq:enlargedalg1}
\clb_{eq}:=\{a\otimes P + b\otimes (I-P): a,b\in\clb_\ell\}
\ee
%%%%%%%%%%%%%%%%%%%%%%%%%%%%%%%%%%%%%%%%%%%%%%
Clearly $\clb_{eq}$ contains $C^\infty_{eq}(S_0^{2\ell+1})$.
We will show that $\clb_{eq}$ is closed under derivations with
both $|D_{eq}|$ as well as $D_{eq}$. This will prove regularity
of  the spectral triple $( C^\infty_{eq}(S_0^{2\ell+1}),\clh,D_{eq})$.

Note that $|D_{eq}|=|D_\ell|\otimes I + I\otimes \widetilde{N}$.
Since $I\otimes \widetilde{N}$ commutes with every element of $\clb_{eq}$,
we get
$\delta(a\otimes P + b\otimes (I-P))=[|D_\ell|,a]\otimes P+ [|D_\ell|,b]\otimes(I-P)$
and
$[D_{eq}, a\otimes P + b\otimes (I-P)]=[D_\ell,a]\otimes P - [|D_\ell|,b]\otimes(I-P)$.
Since $\clb_\ell$ is closed under derivations with $|D_\ell|$
and $D_\ell$, 
it follows that $\clb_{eq}$ is closed under derivations
with $|D_{eq}|$ and $D_{eq}$.

Next we compute the dimension spectrum of the spectral triple.
For $w \in \mathbb{T}^{\ell+1}$, let $\tilde{U}_{w}:=U_{w}\otimes I$. Then
$|D_{eq}|$ commutes with $\tilde{U}_{w}$. Hence again it is enough to consider
homogeneous elements of degree $0$. Now by lemma \ref{dimension} it follows that
for $b \in \mathcal{B}_{eq}$ with $b$ homogeneous of degree $0$, 
the function $Trace(b|D_{eq}|^{-z})$ is meromorphic with simple poles and the poles lie in
$\{1,2,\cdots,2\ell +1\}$. To show that every point of $\{1,2,\cdots,2\ell +1\}$ 
is in the dimension spectrum, observe
that 
\begin{equation}
  \label{poles for Deq}
 Trace(|D_{eq}|^{-z})= \Sigma_{k=0}^{2\ell }(2c^{2\ell }_{k}-c^{2\ell -1}_{k})\zeta(z-k) 
\end{equation}
where $c_{k}^{r}$ is defined as the coefficient of $N^{k}$ in $\binom{N+r}{r}$.
Note that for $0 \leq k \leq r$ one has $c_{k}^{r} > 0$.  Also note the
recurrence $rc_{k}^{r}= c_{k-1}^{r-1} + rc_{k}^{r-1}$. Hence $c_{k}^{r} \geq
c_{k}^{r-1}$. Now from equation~(\ref{poles for Deq}) it  follows that
$Res_{z=k+1}Trace(|D_{eq}|^{-z})=2c_{k}^{2\ell }-c_{k}^{2\ell -1} > 0$ for $0 \leq k \leq
2\ell $. This proves that every point of $\{1,2,\cdots,2\ell +1\}$ is in the dimension
spectrum. This completes the proof. 
\qed

We will need the fact that $Trace(|D_{eq}|^{-z})$ is meromorphic with simple poles
at $\{1,2,\cdots,2\ell +1\}$ with non-zero residue and hence we state it as a
separate lemma.
%%%%%%%%%%%%%%%%%%%%%%%%%%%%%%%%%%%%%%%%%%%%%%
\blmma
\label{poles for D_{eq}}
 The function $Trace(|D_{eq}|^{-z})$ is meromorphic with simple poles at
$\{1,2,\cdots,2\ell +1\}$. Also for $k \in \{1,2,\cdots,2\ell +1\}$, the residue
$Res_{z=k}Trace(|D_{eq}|^{-z})$ is non-zero.
\elmma
%%%%%%%%%%%%%%%%%%%%%%%%%%%%%%%%%%%%%%%%%%%%%%

%%%%%%%%%%%%%%%%%%%%%%%%%%%%%%%%%%%%%%%%%%%%%%
%%%%%%%%%%%%%%%%%%%%%%%%%%%%%%%%%%%%%%%%%%%%%%

%%%%%%%%%%%%%%%%%%%%%%%%%%%%%%%%%%%%%%%%%%%%%%%%%%%%%%%%%%%%%%%%%%%
\subsection{Regularity and dimension spectrum for $q\neq 0$}
%%%%%%%%%%%%%%%%%%%%%%%%%%%%%%%%%%%%%%%%%%%%%%%%%%%%%%%%%%%%%%%%%%%
%%%%%%%%%%%%%%%%%%%%%%%%%%%%%%%%%%%%%%%%%%%%%%
Consider the smooth subalgebra of the Toeplitz algebra defined as:
\bea
\scrt^\infty&=&\left\{\sum_{j,k \in \mathbb{N}}\lambda_{jk}~S^{*j}p_0 S^{k}+
\sum_{k\geq 0}\lambda_{k}S^{k}+ \sum_{k>0}\lambda_{-k}S^{*k}:
                 \lambda_{jk}, (\lambda_{k}) \text{ are ~rapidly decreasing } \right\} 
\nonumber
 \eea
For $a:= \sum_{j,k \in \mathbb{N}}\lambda_{jk}~S^{*j}p_0 S^{k}+  \sum_{k\geq
0}\lambda_{k}S^{k}+  \sum_{k>0}\lambda_{-k}S^{*k} \in \scrt^{\infty}$,
define the seminorm $\|\cdot\|_{m}$ by  $\| a \|_{m}:=
 \sum (1+|j|+|k|)^{m} |\lambda_{kl} | +  \sum (1+| k |)^{m}
| \lambda_{k} |$. Equipped with this family of seminorms,
$\scrt^{\infty}$ is a Fr\'echet algebra.
We will denote by $\scrt^\infty_k$ the $k$-fold tensor
product of $\scrt^\infty$. 

%%%%%%%%%%%%%%%%%%%%%%%%%%%%%%%%%%%%%%%%%%%%%%
\blmma
The triple $(\scrt^\infty,\ell_{2}(\mathbb{N}),N)$ is a regular spectral
triple. More precisely, $\scrt^\infty$ is contained in $Dom(\delta)$ where
$\delta$ is the unbounded derivation $[N,.]$ and $\delta$ leaves the algebra
$\scrt^{\infty}$ invariant. Also the map $\delta:\scrt^{\infty}\to
\scrt^{\infty}$ is continuous. 
\elmma
%%%%%%%%%%%%%%%%%%%%%%%%%%%%%%%%%%%%%%%%%%%%%%
\prf
Note that $[N,S]=-S$ and $[N,p]=0$. Now the lemma follows from
the fact that the unbounded derivation $\delta$ is closed.
\qed

For $\alpha \in \mathbb{N}^{2} \cup \mathbb{Z}$, let
\begin{displaymath}
\begin{array}{lll}
W_{\alpha} &= &\left \{ \begin{array}{lll}
                        S^{*m}p_0 S^{n} & \text{if} & \alpha=(m,n) \\
                        S^{r}        & \text{if} & \alpha=r \geq 0 \\
                        S^{*r}       & \text{if} & \alpha=r<0 \\
                       \end{array} \right. 

\end{array}  
\end{displaymath}
For $\alpha \in \mathbb{N}^{2} \cup \mathbb{Z}$, define $| \alpha |$ to be 
$|m|+| n|$ if $\alpha=(m,n) \in \mathbb{N}^{2}$ and  the
usual absolute value $|\alpha|$ if $\alpha \in \mathbb{Z}$. For an $\ell$ tuple
$\alpha=(\alpha_{1},\alpha_{2},\cdots,\alpha_{\ell})$ in $(\bbn^2\cup\bbz)^\ell$, 
let $| \alpha | = \sum | \alpha_{i}|$ and $W_{\alpha}:=W_{\alpha_{1}}\otimes
W_{\alpha_{2}}\otimes \cdots W_{\alpha_{\ell}}$. We need the following simple lemma
whose proof we omit as it is easy to prove.

%%%%%%%%%%%%%%%%%%%%%%%%%%%%%%%%%%%%%%%%%%%%%
\blmma
The natural tensor product representation of $\scrt_\ell^\infty $
on $\ell_{2}(\mathbb{N})^{\otimes \ell }$ is injective. Thus we identify
$\scrt_\ell^\infty $ with it's range which is $\{  \sum
x_{\alpha}W_{\alpha}:    \sum (1+| \alpha |)^{p} | x_{\alpha}|<\infty$ for
every $p\}$.
\elmma
%%%%%%%%%%%%%%%%%%%%%%%%%%%%%%%%%%%%%%%%%%%%%

%%%%%%%%%%%%%%%%%%%%%%%%%%%%%%%%%%%%%%%%%%%%%
\brmrk
 The tensor product representation of $OP^{-\infty}_{D_{\ell} }\otimes
\scrt_\ell^\infty $ on $\cll(\mathcal{H}_{\Sigma})$ is injective since
$OP^{-\infty}_{D_{\ell} }:= \mathcal{S}(\clh_\ell)$
and hence we identify $OP^{-\infty}_{D_{\ell} }\otimes \scrt_\ell^\infty $
with its image.
\ermrk
%%%%%%%%%%%%%%%%%%%%%%%%%%%%%%%%%%%%%%%%%%%%%

For an operator $T$, let $L_{T}$ denote the left
multiplication map $X\mapsto TX$.
Then for $T\in OP^0_{D_\ell}$, the map 
$L_T:OP^{-\infty}_{D_{\ell} }\rightarrow OP^{-\infty}_{D_{\ell} }$
is continuous. 
Note that if $A$ is a Fr\'echet algebra  and
$a \in A$, then $L_{a}$  is a continuous linear operator.
%%%%%%%%%%%%%%%%%%%%%%%%%%%%%%%%%%%%%%%%%%%%%%%%%%%%%
\blmma
\label{left multiplication}
Let $T \in OP^{0}_{D_{\ell} }$ and $a \in \scrt^{\infty}_\ell$.
Then the map
$L_{T\otimes a}$ leaves the algebra
$OP^{-\infty}_{D_{\ell} }\otimes \scrt^{\infty}_{\ell}$ invariant.
Moreover $L_{T\otimes a}=L_{T}\otimes L_{a}$ on the algebra
$OP^{-\infty}_{D_{\ell} } \otimes\scrt^{\infty}_{\ell}$.
\elmma
%%%%%%%%%%%%%%%%%%%%%%%%%%%%%%%%%%%%%%%%%%%%%%%%%%%%%
\prf
 Clearly $L_{T \otimes a}=L_{T}\otimes L_{a}$ on the algebraic
tensor product $OP^{-\infty}_{D_{\ell} }\otimes_{alg}\scrt^{\infty}_{\ell}$. Now let
$a \in OP^{-\infty}_{D_{\ell} }\otimes\scrt^{\infty}_{\ell}$. Then there
exists a sequence 
$a_{n} \in OP^{-\infty}_{D_{\ell} }\otimes_{alg}\scrt^{\infty}_{\ell}$ 
which converges to $a$ in
$OP^{-\infty}_{D_{\ell} }\otimes\scrt^{\infty}_{\ell}$. Also $a_{n}$
converges to $a$ in the operator norm. Now the result follows from the
continuity of $L_{T\otimes a}$ and $L_{T}\otimes L_{a}$.
\qed

%%%%%%%%%%%%%%%%%%%%%%%%%%%%%%%%%%%%%%%%%%%%%
\bppsn
\label{TheAlgebraB}
Let 
\be
\mathscr{B}:= \mathcal{B}_{eq}+OP^{-\infty}_{D_{\ell} } \otimes
\scrt_\ell^\infty 
\ee
Then one has the following.
\begin{enumerate}
\item The vector space $\mathscr{B}$ is an algebra.
\item The algebra $\mathscr{B}$ is invariant under the derivations
$\delta:=[|D_{eq}|,.]$ and $[D_{eq},.]$.
\item For $b \in \mathscr{B}$, the commutator $[F_{eq},b] \in
OP^{-\infty}_{D_{eq}}$.
\item For $b \in \mathscr{B}$, the function $Trace(b|D_{eq}|^{-z})$ is meromorphic
with only simple poles and the poles lie in  $\{1,2,\cdots,2\ell +1\}$.
\end{enumerate}
\eppsn
%%%%%%%%%%%%%%%%%%%%%%%%%%%%%%%%%%%%%%%%%%%%%
\prf
Lemma~\ref{left multiplication} and the fact that
$\mathcal{B}_{\ell} \subset OP^{0}$ implies that $\mathscr{B}$ is an algebra. As
seen in Theorem~\ref{thm:regularityq=0}, it follows that $\mathcal{B}_{eq}$ is
invariant under $\delta$ and $[D_{eq},.]$. Also $(3)$ and $(4)$ holds for $b \in
\mathcal{B}_{eq}$. Hence to complete the proof it is enough to consider
$(2),(3)$ and $(4)$ for the algebra $OP^{-\infty}_{D_{\ell} }\otimes
\scrt_\ell^\infty $.

Lemma~\ref{ProductofSpectralTriples} and  the decomposition
$|D_{eq}|=|D_{\ell} |\otimes 1 + 1\otimes \tilde{N}$ implies that $\delta$ leaves
the algebra $OP^{-\infty}_{D_{\ell} } \otimes \scrt_\ell^\infty $
invariant. Now note that $P \in OP^{-\infty}_{\tilde{N}}$ , it follows that left
and right multiplication by $F_{\ell}\otimes P$ and $1\otimes P$ sends
$OP^{-\infty}_{D_{\ell} }\otimes \scrt_\ell^\infty $ to
$OP^{-\infty}_{D_{eq}}\equiv OP^{-\infty}_{D_{\ell} }\otimes OP^{-\infty}_{\tilde{N}}$.
Now since $F_{eq}=F_{\ell}\otimes P - I\otimes(I-P)$, it follows that $[F_{eq},b]$
is smoothing for every $b \in OP^{-\infty}_{D_{\ell} }\otimes \scrt^{\infty}$. Now
the invariance of $OP^{-\infty}_{D_{\ell} }\otimes \scrt_\ell^\infty $
under $[D_{eq},.]$ follows from the equation $[D_{eq},b]=\delta(b)F_{eq}+
|D_{eq}|[F_{eq},b]$ and the fact that $OP^{-\infty}_{D_{eq}}:=
OP^{-\infty}_{D_{\ell} }\otimes OP^{-\infty}_{\tilde{N}}$ is contained in
$OP^{-\infty}_{D_{\ell} }\otimes \scrt_\ell^\infty $.   

Now we will prove that for $b \in OP^{-\infty}_{D_{\ell} }\otimes
\scrt^{\infty}_\ell$, the function 
$Trace(b|D_{eq}|^{-z})$ is meromorphic with simple
poles and the poles lie in $\{1,2,\cdots,\ell\}$. For $w \in \mathbb{T}^{2\ell +1}$,
let $U_{w}= U_{w_{1}}\otimes U_{w_{2}} \otimes \cdots U_{w_{2\ell +1}}$ be the
unitary operator on $\mathcal{H}_{\Sigma}$. Clearly $U_{w}|D_{eq}|U_{w}^{*}=
|D_{eq}|$ for $w \in \mathbb{T}^{2\ell +1}$. Hence it is enough to consider
$Trace(b|D_{eq}|^{-z})$ with $b$ homogeneous of degree $0$. 

An element $b$ is homogeneous if and only if it commutes with
the operators $U_w$ for all $w\in\mathbb{T}^{2\ell +1}$.
This implies $b$ must be of the form
$ e_\gamma\mapsto \phi(\gamma)e_\gamma$ for some function $\phi$,
i.e.\ $b=\sum_\gamma \phi(\gamma)p_\gamma$.
An operator of the form $\sum_{\gamma\in\Gamma_{\Sigma_\ell}}\phi(\gamma)p_\gamma$
is in $OP^{-\infty}_{D_\ell}$ if and only if $\phi(\gamma)$ is rapidly
decaying on $\Gamma_{\Sigma_\ell}$.
Also, using the description of $\scrt^\infty$, it follows
that an operator of the form 
$\sum_{n\in\bbn}\phi(n)p_n$ belongs to $\scrt^\infty$
if and only if $\phi(\cdot)-\lim_{n\rightarrow\infty}\phi(n)$
is rapidly decreasing.
Thus combining these, one can see that
the operator $\sum_\gamma \phi(\gamma)p_\gamma$ belongs to 
$OP_{D_\ell}^{-\infty}\otimes\scrt_\ell^\infty$
if and only if $\phi$ is a linear combination of
$\phi_A$ with $A$ varying over subsets of $\Sigma$ containing
$\Sigma_\ell$, where each 
$\phi_A(\gamma)$ depends only on $\gamma_A$ 
and $\phi_A(\gamma_A)$ is rapidly decreasing
on $\Gamma_A$.
For an element $b=\sum_\gamma \phi_A(\gamma)p_\gamma$, one has
\[
 Trace(b|D_{eq}|^{-z})=\sum_\gamma \frac{\phi_A(\gamma)}{|\gamma|^z}
  =\sum_\gamma \frac{\phi_A(\gamma_A)}{(|\gamma_A|+|\gamma_{\Sigma\backslash A}|)^z}.
\]
By lemma~\ref{dimension} it follows that $Trace(b|D_{eq}|^{-z})$ is meromorphic
with simple poles and the poles lie in 
$\{1,2,\cdots,|\Sigma\backslash A|\}\seq \{1,2,\ldots,\ell\}$. This completes the
proof.
% % %%%%%%%%%%%%%%%%%%%%%%%%%%%%%%%%%%%%%%%
\qed

%%%%%%%%%%%%%%%%%%%%%%%%%%%%%%%%%%%%%%%%%%%%%%%%%%%%%%%%%%%%
%%%%%%%%%%%%%%%%%%%%%%%%%%%%%%%%%%%%%%%%%%%%%%%%%%%%%%%%%%%%

%%%%%%%%%%%%%%%%%%%%%%%%%%%%%%%%%%%%%%%%%%%%%%%%%%%%%%%%%%%%
\subsection{The smooth function algebra $C^{\infty}(S_{q}^{2\ell +1})$}
%%%%%%%%%%%%%%%%%%%%%%%%%%%%%%%%%%%%%%%%%%%%%%%%%%%%%%5%%%%%%
%%%%%%%%%%%%%%%%%%%%%%%%%%%%%%%%%%%%%%%%%%%%%%%%%%%%%%%%%%%%
In this subsection we will define a dense $*$ Fr\'echet algebra
$C^{\infty}(S_{q}^{2\ell +1})$ of $C_{\ell}=\pi(C(S_{q}^{2\ell +1}))$ and show that it is closed under
holomorphic functional calculus.  Let $B_{\ell}$ be the $C^{*}$ algebra generated
by $A_{\ell}$ and $F_{\ell}$. 
% Note that $B_{\ell}$ contains 
% $\mathcal{K}(\ell_{2}(\mathbb{N})^{\otimes \ell }\otimes \ell_{2}(\mathbb{Z}))$.
Recall that $\cle$ denotes the $C^{*}$ algebra generated by $C(\mathbb{T})$ and $F_0$.
%%%%%%%%%%%%%%%%%%%%%%%%%%%%%%%%%%%%%%%%%%%%%
\blmma
\label{l=0}
 The $C^{*}$ algebra $\cle$ contains $\mathcal{K}$ and $\cle/\mathcal{K}$ is
isomorphic to the $C^{*}$ algebra $C(\mathbb{T})\oplus C(\mathbb{T})$.
\elmma
%%%%%%%%%%%%%%%%%%%%%%%%%%%%%%%%%%%%%%%%%%%%%
\prf
Let $|e_{m}\rangle\langle e_{n}|$ be the matrix units in $\mathcal{K}(\ell_2(\bbz))$. Note that
$[F_0, S^*]S=2|e_{0}\rangle\langle e_{0}|$. Hence $p_0\equiv|e_{0}\rangle\langle e_{0}| \in \cle$. 
Now $S^{*m}p_0 S^{n}=|e_{m}\rangle\langle e_{n}|$. Hence $\mathcal{K}
\subset \cle$. Let $P_0:=\frac{1+F_0}{2}$. Then $[P_0,f]$ is compact for every $f$. Thus
$\cle/\mathcal{K}$ is generated by $C(\mathbb{T})$ and a projection $P_0$ which is in
the center of $\cle/\mathcal{K}$. Now consider the map 
\begin{displaymath}
 C(\mathbb{T})\oplus C(\mathbb{T})\ni (f,g)  \mapsto fP_0+g(1-P_0)\quad(mod ~\mathcal{K}) \in
\cle/\mathcal{K}.
\end{displaymath}
 We claim that this map is an isomorphism. To prove this, we need to show that if $fP_0$
is compact then $f=0$ and if $g(1-P_0)$ is compact then $g=0$. 

Assume that $fP_0$ is compact for $f \in C(\mathbb{T})$. Fix an $r \in
\mathbb{Z}$. Since $fP_0$ is compact, it follows that $|\langle fP_0(e_{n}),e_{n+r}\rangle
| = | \hat{f}(r) | $ converges to $0$ as $n \to +\infty$. Hence
$\hat{f}(r)=0$ for every $r$. This proves that $f=0$. Similarly one can show
that if $g(1-P_0)$ is compact then $g=0$. This completes the proof.
\qed

%%%%%%%%%%%%%%%%%%%%%%%%%%%%%%%%%%%%%%%%%%%%%
\blmma
\label{extension}
 The $C^*$-algebra $B_{\ell}$ contains $\mathcal{K}(\clh_\ell)$ and the map $(a,b) \mapsto  aP_\ell+b(1-P_\ell)~mod ~\mathcal{K}$ 
 from  $C(S_{q}^{2\ell +1}) \oplus C(S_{q}^{2\ell +1})$
to $B_{\ell}/\mathcal{K}(\clh_\ell)$ is an
isomorphism.
\elmma
%%%%%%%%%%%%%%%%%%%%%%%%%%%%%%%%%%%%%%%%%%%%%
\prf
For $\ell=0$ this is just lemma~\ref{l=0}. 
So let us prove the statement for $\ell\geq 1$.
Since $A_\ell$ contains 
$\mathcal{K}(\ell_2(\bbn^\ell)) \otimes C(\mathbb{T})$, 
it follows that $B_{\ell}$ contains $\mathcal{K}(\clh_\ell)$.
Observe that $[P_{\ell},\alpha_{i}]=0$ for $1 \leq i \leq \ell$ 
and $[P_{\ell},\alpha_{\ell+1}]$ is compact.
Therefore it follows that $[P_{\ell},a]$ is compact for every $a \in A_\ell$.
Hence the map $(a,b)\mapsto aP_{\ell}+b(1-P_{\ell})~mod~ \mathcal{K}$
from $A_\ell\oplus A_\ell$ to $B_{\ell}/\mathcal{K}$ 
is a $*$ algebra homomorphism onto $B_{\ell}/\mathcal{K}$. We will
show that the map is one-one. For that we have to show if $aP_{\ell}$ is compact with
$a \in A_\ell$ then $a=0$ and if $b(1-P_{\ell})$ is compact with 
$b \in A_\ell$ then $b=0$.

Suppose now that $aP_{\ell}$ is compact. Observe that 
$B_{\ell} \subset \scrt_{\ell}\otimes \cle$
and $aP_\ell=a(I\otimes P_0)$.
Since $aP_\ell$ is compact, if we apply the symbol map $\sigma$
on the $\ell$\raisebox{.4ex}{th} copy of $\scrt$,
we get $\sigma_\ell(a)\otimes P_0=0$.
Hence $a$ is in the ideal 
$\mathcal{K}(\ell_{2}(\mathbb{N})^{\otimes \ell }) \otimes C(\mathbb{T})$. For $m,n
\in \mathbb{N}^{\ell}$, let $e_{mn}$ be the ``matrix'' units. Let
$a_{mn}=(e_{mm}\otimes 1)a(e_{nn} \otimes 1)$. Then $a_{mn}=e_{mn} \otimes
f_{mn}$ for some $f_{mn} \in C(\mathbb{T})$. Since $aP_\ell$ is compact, it follows
that $f_{mn}P_0$ is compact as $P_\ell=I\otimes P_0$ 
commutes with $e_{nn} \otimes I$. By the $\ell=0$
case, it follows that $f_{mn}=0$ and hence $a_{mn}=0$ for every $m,n$. Thus
$a=0$. Similarly one can show that if $b(1-P_\ell)$ is compact then $b=0$. This
completes the proof.
\qed

Let $\mathbb{B}$ be the $C^*$ algebra on $\mathcal{H}_{\Sigma}$ generated by
$A_\ell\otimes I$, $P_\ell\otimes 1$ and $1\otimes P$ and 
$J:=\mathcal{K}(\clh_\ell)\otimes \scrt_\ell$. 
Note that $J$ is an ideal since $\mathbb{B}_{\ell}$ is contained
in $\scrt_{\ell }\otimes \cle\otimes \scrt_{\ell }$. The next proposition
identifies the quotient $\mathbb{B}/J$.
%%%%%%%%%%%%%%%%%%%%%%%%%%%%%%%%%%%%%%%%%%%%%
\bppsn
\label{exactness}
Let  $\rho:A_\ell\oplus A_\ell\oplus A_\ell\oplus A_\ell\to \mathbb{B}/J$ be
the map 
\begin{equation*}
(a_{1},a_{2},a_{3},a_{4})\mapsto 
     a_{1}P_\ell\otimes P + a_{2}P_\ell\otimes (1-P) +
a_{3}(1-P_\ell)\otimes P + a_{4}(1-P_\ell)\otimes (1-P) %mod J.
\end{equation*}
from $A_\ell\oplus A_\ell\oplus A_\ell\oplus A_\ell$ into $\mathbb{B}$
composed with the canonical projection from $\mathbb{B}$ onto $\mathbb{B}/J$.
Then $\rho$ is an isomorphism.
\eppsn
%%%%%%%%%%%%%%%%%%%%%%%%%%%%%%%%%%%%%%%%%%%%%
\prf
First note that since $[P_\ell,a] \in \mathcal{K}$ 
for $a \in A_\ell$, it follows that $P_\ell\otimes I$ and $I\otimes P$ are in the
center of $\mathbb{B}/J$. Hence the map $\rho$ is an algebra
homomorphism. By the definition of $\mathbb{B}$ it follows that $\rho$ is
onto. Thus we have to show $\rho$ is one-one.

Suppose that 
$a=a_{1}P_\ell\otimes P + a_{2}P_\ell\otimes(1-P)+ 
   a_{3}(1-P_\ell)\otimes P + a_{4}(1-P_\ell)\otimes(1-P) \in J$. 
Let ${\epsilon}:\scrt \to \mathbb{C}$
be the map $ev_1\circ\sigma$, where $ev_1$ is evaluation at the point $1$. Now consider
the map $id \otimes {\epsilon}^{\otimes \ell }:\scrt_{\ell } \otimes \cle \otimes
\scrt_{\ell } \to \scrt_{\ell }\otimes \cle$. Note that $I \otimes
{\epsilon}^{\otimes \ell }$ sends $J$ to $\mathcal{K}(\mathcal{H}_{\ell})$.
Hence $(I \otimes {\epsilon}^{\otimes \ell })(a)=a_{2}P_\ell+ a_{4}(1-P_\ell) \in
\mathcal{K}(\mathcal{H}_{\ell})$. Hence by lemma \ref{extension}, it follows that $a_{2}=0=a_{4}$.
Since left multiplication by $I\otimes P$ sends the ideal $J$ to
$\mathcal{K}(\mathcal{H}_{\Sigma})$. It follows that $(I\otimes P)a = a_{1}P_\ell
\otimes P + a_{3}(1-P_\ell) \otimes P$ is compact. Hence $a_{1}P_\ell + a_{3}(1-P_\ell)$ is
compact. Thus again by lemma \ref{extension}, it follows that $a_{1}=0=a_{3}$.
This completes the proof.
\qed

Now we prove that $\mathscr{B}$ is closed under holomorphic functional calculus
in $\mathbb{B}$. Let $\mathscr{J}:=OP^{-\infty}_{D_{\ell} }\otimes \scrt_{\ell}^{\infty}$. 
Note that
\begin{align*}
\mathscr{B}:=& \{a_{1}P_\ell\otimes P + a_{2}P_\ell\otimes (1-P) + a_{3}(1-P_\ell)\otimes P +
a_{4}(1-P_\ell)\otimes (1-P)+R :\\
            &   a_{1},a_{2},a_{3},a_{4} \in A_{\ell}^{\infty},  R \in \mathscr{J}
\} \end{align*} 

%%%%%%%%%%%%%%%%%%%%%%%%%%%%%%%%%%%%%%%%%%%%%%%%%%%%%%%%%%%%%
\bppsn
\label{h.f.c}
The algebra $\mathscr{B}$ has the following properties:
\begin{enumerate}
\item If $a_{1}P_\ell \otimes P + a_{2}P_\ell\otimes(1-P) + a_{3}(1-P_\ell)\otimes P +
a_{4}(1-P_\ell)\otimes (1-P) \in \mathscr{J}$ then $a_{i}=0$ for $i=1,2,3,4$. Hence
$\mathscr{B}$ is isomorphic to the direct sum $A_{\ell}^{\infty}\oplus
A_{\ell}^{\infty}\oplus A_{\ell}^{\infty}\oplus A_{\ell}^{\infty} \oplus \mathscr{J}$.
Equip $\mathscr{B}$ with the Fr\'echet space structure coming from this direct sum
decomposition.
\item The algebra $\mathscr{B}$ is a $*$-Fr\'echet algebra contained in
$\mathbb{B}$. Moreover the inclusion $\mathscr{B} \subset \mathbb{B}$ is
continuous.
\item The algebra $\mathscr{B}$ is closed under holomorphic functional calculus
in $\mathbb{B}$.
\end{enumerate}
\eppsn
%%%%%%%%%%%%%%%%%%%%%%%%%%%%%%%%%%%%%%%%%%%%%%%%%%%%%%%%%%%%%
\prf
Proposition \ref{exactness} implies $(1)$. Parts $(2)$ and
$(3)$ follows from proposition \ref{regularity}. Now by proposition
\ref{exactness} one has the exact sequence
\begin{displaymath}
 0 \rightarrow J \rightarrow \mathbb{B} \rightarrow A_{\ell}\oplus A_{\ell} \oplus
A_{\ell} \oplus A_{\ell} \rightarrow 0.
\end{displaymath}
At the smooth algebra level we have the following  exact sequence
\begin{displaymath}  
0 \rightarrow \mathscr{J} \rightarrow \mathscr{B}\stackrel{\theta}{\rightarrow} A_{\ell}^{\infty}
\oplus A_{\ell}^{\infty} \oplus A_{\ell}^{\infty} \oplus A_{\ell}^{\infty} 
     \rightarrow 0.
\end{displaymath}
Since $\mathscr{J}\subset J$ and $A_{\ell}^{\infty} \subset A_{\ell}$ 
are closed under holomorphic functional calculus,
it follows from theorem~3.2, part~2, \cite{sch-1993a} that $\scrb$ is 
spectrally invariant in $\mathbb{B}$. 
Since by part~(2), the Fr\'echet topology of $\scrb$ is finer
than the norm topology, it follows that $\scrb$ is
closed in the holomorphic function calculus of $\mathbb{B}$.
\qed

%%%%%%%%%%%%%%%%%%%%%%%%%%%%%%%%%%%%%%%%%%%%
\brmrk
 One can prove that $OP^{-\infty}_{D_{\ell} }\otimes \scrt_\ell^\infty $
is closed under holomorphic functional calculus in 
$\mathcal{K}(\clh_\ell) \otimes \scrt_{\ell}$ in the same manner by applying 
theorem~3.2, part~2, \cite{sch-1993a} and by using
the extension (after tensoring suitably)  
 \[
  0 \rightarrow \mathcal{K} \rightarrow \scrt \rightarrow C(\mathbb{T})
\rightarrow 0
\] 
at the $C^{*}$ algebra level and the extension 
\[
 0 \rightarrow \mathcal{S}(\ell_{2}(\mathbb{N})) \rightarrow \scrt^{\infty}
\rightarrow C^{\infty}(\mathbb{T}) \rightarrow 0
\]
at the Fr\'echet algebra level.
\ermrk
%%%%%%%%%%%%%%%%%%%%%%%%%%%%%%%%%%%%%%%%%%%%

%%%%%%%%%%%%%%%%%%%%%%%%%%%%%%%%%%%%%%%%%%%%
\bcrlre
Define the smooth function algebra $C^{\infty}(S_{q}^{2\ell +1})$
by
\[
 C^{\infty}(S_{q}^{2\ell +1})=\{a\in\scrb\cap C_\ell:\theta(a)\in\iota(\cla_\ell^\infty)\},
\]
where $\theta$ is as in the proof of proposition~\ref{h.f.c}
and $\iota:A_\ell\rightarrow A_\ell\oplus A_\ell\oplus A_\ell\oplus A_\ell$
is the inclusion map $a\mapsto a\oplus a \oplus a\oplus a$.
 Then the algebra $C^{\infty}(S_{q}^{2\ell +1})$ is closed in $\mathscr{B}$ and it is
closed under holomorphic functional calculus in $C_{\ell}$.
\ecrlre
%%%%%%%%%%%%%%%%%%%%%%%%%%%%%%%%%%%%%%%%%%%%
\prf
Let $j:\mathscr{B} \to \cll(\mathcal{H}_{\Sigma})$ denote the
inclusion map. Then by definition $C^{\infty}(S_{q}^{2\ell +1})=
\theta^{-1}(\iota(A_{\ell}^{\infty})) \cap j^{-1}(C_{\ell})$. Since $\theta$ and
$j$ are continuous and as $\iota(A_{\ell}^{\infty})$ and $C_{\ell}$ are closed, it follows
that $C^{\infty}(S_{q}^{2\ell +1})$ is closed in $\mathscr{B}$. Hence
$C^{\infty}(S_{q}^{2\ell +1})$ is a Fr\'echet algebra. Also $C^{\infty}(S_{q}^{2\ell +1})$
is $*$-closed as $\rho$ is $*$-preserving. Now let $a \in C^\infty(S_{q}^{2\ell +1})$ be
invertible in $C_{\ell}$. Then $a$ is invertible in $\cll(\mathcal{H}_{\Sigma})$.
By proposition \ref{h.f.c}, it follows that $a^{-1} \in \mathscr{B}$. By the
closedness of $A_{\ell}^{\infty}$ under holomorphic functional calculus, it follows
that $\theta(a^{-1}) \in \iota(A_{\ell}^{\infty})$. Thus one has 
$a^{-1} \in C^{\infty}(S_{q}^{2\ell +1})$. 
We have already seen that the Fr\'echet topology
of $\scrb$ is finer than the norm topology. Same is therefore true
for the topology of $C^{\infty}(S_{q}^{2\ell +1})$.
Hence it is closed under holomorphic functional calculus in $C_{\ell}$.
\qed

%%%%%%%%%%%%%%%%%%%%%%%%%%%%%%%%%%%%%%%%%%%%
\bppsn\label{pro:zjq}
The operators $Z_{j,q} \in C^{\infty}(S_{q}^{2\ell +1})$. Hence
$C^{\infty}(S_{q}^{2\ell +1})$ is a dense subalgebra of $C_{\ell}$
that contains $\pi(\cla(S_q^{2\ell+1}))$.
\eppsn
%%%%%%%%%%%%%%%%%%%%%%%%%%%%%%%%%%%%%%%%%%%%
The proof of this proposition will be given in the next subsection. 

We are now in a position to prove the main theorem.
%%%%%%%%%%%%%%%%%%%%%%%%%%%%%%%%%%%%%%%%%%%%
\bthm
The triple $(C^{\infty}(S_{q}^{2\ell +1}),\mathcal{H}_{\Sigma},D_{eq})$ is a
regular spectral triple with simple dimension spectrum $\{1,2,\cdots,2\ell +1\}$.
\ethm
%%%%%%%%%%%%%%%%%%%%%%%%%%%%%%%%%%%%%%%%%%%%
\prf
Since $C^{\infty}(S_{q}^{2\ell +1}) \subset \mathscr{B}$ the regularity of the
spectral triple $(C^{\infty}(S_{q}^{2\ell +1}),\mathcal{H}_{\Sigma},D_{eq})$
follows from the regularity of the spectral triple
$(\mathscr{B},\mathcal{H}_{\Sigma},D_{eq})$ which is proved in 
proposition~\ref{TheAlgebraB}. Proposition~\ref{TheAlgebraB} also 
implies that the spectral triple has simple dimension spectrum 
which is a subset of $\{1,2,\cdots,2\ell +1\}$. The fact
that every point in $\{1,2,\cdots,2\ell +1\}$ is in the dimension spectrum follows
from lemma \ref{poles for D_{eq}}. This completes the proof.
\qed

%%%%%%%%%%%%%%%%%%%%%%%%%%%%%%%%%%%%%%%%%%%%%%
%%%%%%%%%%%%%%%%%%%%%%%%%%%%%%%%%%%%%%%%%%%%%%

%%%%%%%%%%%%%%%%%%%%%%%%%%%%%%%%%%%%%%%%%%%%%%
%%%%%%%%%%%%%%%%%%%%%%%%%%%%%%%%%%%%%%%%%%%%%%

%%%%%%%%%%%%%%%%%%%%%%%%%%%%%%%%%%%%%%%%%%%%%%
\subsection{The operators $Z_{j,q}$}
%%%%%%%%%%%%%%%%%%%%%%%%%%%%%%%%%%%%%%%%%%%%%%
We will give a proof of proposition~\ref{pro:zjq}
in this subsection. The main idea will be to exploit the
isomorphism between the Hilbert spaces $L_2(S_q^{2\ell+1})$
and $\clh_\Sigma$ and a detailed analysis of the 
operators $Z_{j,q}$ to show that certain parts of these
operators can be ignored for the purpose of establishing
regularity and computing dimension spectrum.
Deciding and establishing which parts of these operators
can be ignored is the key step here.
It should be noted here that a similar analysis has been
done by D'Andrea in~\cite{dan-2008a}, where he embeds
$L_2(S_q^{2\ell+1})$ in a bigger Hilbert space and
proves certain approximations for the operators $Z_{j,q}$.
But the approximation there is not strong enough to enable
the computation of dimension spectrum.
Here we prove stronger versions of those approximations,
which have made it possible to use them to compute the dimension
spectrum in the previous subsection. 

We start with a few simple lemmas that we will use 
repeatedly during the computations in this subsection.
%%%%%%%%%%%%%%%%%%%%%%%%%%%%%%%%%%%%%%%%%%%%%%
\blmma\label{lem:opinclusion}
Let $A\seq B\seq \Sigma$.
Then one has
$OP_{D_B}^{-\infty}\otimes \cle_{\Sigma\backslash B}^\infty
  \seq OP_{D_A}^{-\infty}\otimes \cle_{\Sigma\backslash A}^\infty$.
\elmma
%%%%%%%%%%%%%%%%%%%%%%%%%%%%%%%%%%%%%%%%%%%%%%
\prf
Since 
\[
OP_{D_B}^{-\infty}=\cls(\clh_B)=\cls(\clh_A)\otimes\cls(\clh_{B\backslash A})
=OP_{D_A}^{-\infty}\otimes\cls(\clh_{B\backslash A}),
\]
and $\cls(\clh_{B\backslash A})\seq \cle_{B\backslash A}^\infty$,
we have the required inclusion.
\qed

%%%%%%%%%%%%%%%%%%%%%%%%%%%%%%%%%%%%%%%%%%%%%%

Let $A\seq\Sigma$. Let $\clp$ be a polynomial in $|A|$ variables
and let $T$ be the operator on $\clh_A$ given by
\[
Te_{\gamma}=\clp(\{\gamma_{i},i\in A\})q^{|\gamma_{A}|}e_\gamma.
\]
Since the function
$\gamma\mapsto \clp(\{\gamma_{i},i\in A\})q^{|\gamma_{A}|}$
is a rapid decay function on $\Gamma_A$, it follows that $T\in OP_{D_A}^{-\infty}$. 
%%%%%%%%%%%%%%%%%%%%%%%%%%%%%%%%%%%%%%%%%%%%%%

%%%%%%%%%%%%%%%%%%%%%%%%%%%%%%%%%%%%%%%%%%%%%%
\blmma\label{lem:neglect2}
Let $A\seq\Sigma$. 
Let $T$ and $T_0$ be the following operators on $\clh_{A}$:
\[
 Te_\gamma=q^{\phi(\gamma_A)}Q(\psi(\gamma_A))e_\gamma,\qquad
T_0e_\gamma=q^{\phi(\gamma_A)}e_\gamma,
\]
where $\phi$ and $\psi$ are some nonnegative functions.
If $\phi(\gamma_A)+\psi(\gamma_A)>|\gamma_A|$, then
Then $T-T_0\in OP_{D_{A}}^{-\infty}$. %\otimes \cle_{\Sigma\backslash A\cup B}^\infty$.
\elmma
%%%%%%%%%%%%%%%%%%%%%%%%%%%%%%%%%%%%%%%%%%%%%%
\prf
This is a consequence of the inequality
 $|1-(1-x)^\half|<x$  for  $0\leq x \leq 1$.
\qed

%%%%%%%%%%%%%%%%%%%%%%%%%%%%%%%%%%%%%%%%%%%%%%
\blmma\label{lem:neglect3}
Let $A\seq\Sigma$.
Let $T$ and $T_0$ be  operators on $\clh_{A}$ given by:
\[
 Te_\gamma=q^{\phi(\gamma_A)}Q(\psi(\gamma_A))^{-1}e_\gamma,\qquad
T_0e_\gamma=q^{\phi(\gamma_A)}e_\gamma
\]
for some nonnegative functions $\phi$ and $\psi$. 
If $\phi(\gamma_A)+\psi(\gamma_A)>|\gamma_A|$, then
Then $T-T_0\in OP_{D_{A}}^{-\infty}$. %\otimes \cle_{\Sigma\backslash A\cup B}^\infty$.
\elmma
%%%%%%%%%%%%%%%%%%%%%%%%%%%%%%%%%%%%%%%%%%%%%%
\prf
For $0<r<1$, one has
\[
 |1-(1-x)^{-\half}|<cx \quad \text{ for } 0\leq x \leq r,
\]
where $c$ is some fixed constant that depends on $r$.
Using this, it follows that the map
$\gamma\mapsto q^{\phi(\gamma)}|1-(1-q^{2\psi(\gamma)})^{-\half}|$
is a rapid decay function on $\Gamma_{A}$.
\qed

For $j\in\Sigma$, we will denote by $\cle_j$ the $C^*$-algebra
$\scrt$ if $j\neq\ell+1$ and the $C^*$-algebra $\cle$ if $j=\ell+1$.
Thus $\cle^\infty_j$ will be $\scrt^\infty$ for $j\neq \ell+1$
and $\cle^\infty=\clb$ for $j=\ell+1$. Thus $\cle_\Sigma^\infty$
will stand for the space 
$\scrt_\ell^\infty\otimes\cle^\infty\otimes\scrt_\ell^\infty$.
Note that for any subset $A$ of $\Sigma$, one has 
$OP_{D_A}^{-\infty}\seq\cle_\Sigma^\infty$.
%%%%%%%%%%%%%%%%%%%%%%%%%%%%%%%%%%%%%%%%%%%%%%
\blmma\label{lem:neglect4}
Let $A\seq\Sigma$, $a,b,m,n\in\bbn$ and  $n>0$. 
Let $T_1$ and $T_2$  be the operators on $\clh_{\Sigma}$ given by
\[
 T_1 e_\gamma=Q(|\gamma_{A}|+a(\gamma_{\ell+1})_{+}+b(\gamma_{\ell+1})_{-}+m ) e_\gamma,\qquad
  T_2 e_\gamma=Q(|\gamma_{A}|+a(\gamma_{\ell+1})_{+}+b(\gamma_{\ell+1})_{-} +n)^{- 1} e_\gamma.
\]
Then $T_1$ and $T_2$ are in $\cle_\Sigma^\infty$.
\elmma
%%%%%%%%%%%%%%%%%%%%%%%%%%%%%%%%%%%%%%%%%%%%%%
\prf
First note that if $T'_1$ and $T''_1$ are operators given by
\[
 T'_1 e_\gamma=Q(|\gamma_{A}|+a|\gamma_{\ell+1}|+m ) e_\gamma,\qquad
  T''_1 e_\gamma=Q(|\gamma_{A}|+b|\gamma_{\ell+1}| +m) e_\gamma,
\]
then $T_1= P_\Sigma T'_1 + (I-P_\Sigma)T''_1$,
where $P_\Sigma=\frac{I+F_\Sigma}{2}$.
By the two previous lemmas, $I-T'_1$ and $I-T'_2$ are in 
$OP_{D_B}^{-\infty}$ where $B=A\cup\{\ell+1\}$.
Since  $OP_{D_B}^{-\infty}$ is contained in $\cle_\Sigma^\infty$,
it follows that $T'_1,T'_2\in\cle_\Sigma^\infty$.
Since $P_\Sigma\in\cle_\Sigma^\infty$, we get $T_1\in\cle_\Sigma^\infty$.

Proof for $T_2$ is exactly similar.
\qed

We next proceed with a detailed analysis of the
operators $Z_{j,q}$.
First recall that 
\be\label{eq:ustar1}
 U^* e_\gamma=\xi'_\gamma=
 (-1)^{\sum_{i=1}^\ell (i-1)\gamma_i +\ell(\gamma_{\ell+1})_{+}}e_{\bldr^{n,k},\blds},
\ee
where $\blds$ is given by
\be\label{eq:ustar2}
n=(\gamma_{\ell+1})_{-} +\sum_{i=\ell+2}^{2\ell+1}\gamma_i,\quad
k=\sum_{i=1}^\ell \gamma_i +(\gamma_{\ell+1})_+,
\ee
\be\label{eq:ustar3}
d_m=\sum_{i=1}^{m-1} \gamma_i,\qquad
c_m=\sum_{i=1}^{\ell} \gamma_i +|\gamma_{\ell+1}|+\sum_{i=\ell+2}^{2\ell+2-m} \gamma_i
\quad     \text{ for } 1\leq m\leq \ell.
\ee
\be\label{eq:ustar4}
d_{\ell+1}=c_{\ell+1}=\sum_{i=1}^\ell \gamma_i +(\gamma_{\ell+1})_-.
\ee
We will use this correspondence between $e_{\bldr^{n,k},\blds}$
and $\xi'_\gamma$ freely in what follows.

From equation~(\ref{eq:left_mult3a}), we get
\bean
 \pi(u_{1j})e_{\bldr^{n,k}\blds} &=&
   \sum_{M\in\bbm_j^+} C_q(1,\bldr^{n,k},N_{1,1})C_q(j,\blds,M)
      \kappa(\bldr^{n,k},N_{1,1}) e_{\bldr^{n+1,k},M(\blds)}\\
&& + \sum_{M\in\bbm_j^-} C_q(1,\bldr^{n,k},N_{0,1})C_q(j,\blds,M)
      \kappa(\bldr^{n,k},N_{0,1}) e_{\bldr^{n,k-1},M(\blds)}.
\eean
Therefore
\bean
 Z_{j,q}^*e_{\bldr^{n,k}\blds}  &=&
  q^{-j+1}\sum_{M\in\bbm_j^+} C_q(1,\bldr^{n,k},N_{1,1})C_q(j,\blds,M)
      \kappa(\bldr^{n,k},N_{1,1}) e_{\bldr^{n+1,k},M(\blds)}\\
&& + q^{-j+1}\sum_{M\in\bbm_j^-} C_q(1,\bldr^{n,k},N_{0,1})C_q(j,\blds,M)
      \kappa(\bldr^{n,k},N_{0,1}) e_{\bldr^{n,k-1},M(\blds)}.
\eean
Thus we have $Z_{j,q}^*= \sum_{M\in\bbm_j^+}S_M^+ T_M^+  + \sum_{M\in\bbm_j^-}S_M^-T_M^-$,
where the operators $S_M^\pm$ and $T_M^\pm$ are given by
\bea
S_M^+e_{\bldr^{n,k}\blds} &=& e_{\bldr^{n+1,k},M(\blds)},\quad M\in \bbm_j^+,
     \label{eq:smplus}\\
S_M^-e_{\bldr^{n,k}\blds} &=&  e_{\bldr^{n,k-1},M(\blds)},\quad M\in \bbm_j^-,
       \label{eq:smminus} \\
T_M^+ e_{\bldr^{n,k}\blds} &=& q^{-j+1}C_q(1,\bldr^{n,k},N_{1,1})C_q(j,\blds,M)
      \kappa(\bldr^{n,k},N_{1,1}) e_{\bldr^{n,k}\blds},\quad M\in \bbm_j^+,
       \label{eq:tmplus}  \\ 
T_M^- e_{\bldr^{n,k}\blds} &=& q^{-j+1}C_q(1,\bldr^{n,k},N_{0,1})C_q(j,\blds,M)
      \kappa(\bldr^{n,k},N_{0,1}) e_{\bldr^{n,k}\blds},\quad M\in \bbm_j^-.
     \label{eq:tmminus}
\eea

%%%%%%%%%%%%%%%%%%%%%%%%%%%%%%%%%%%%%%%%%%%%%%
\blmma\label{lem:shift}
Let $S_M^\pm$ be as above. Then $U S_M^\pm U^*\in\cle_\Sigma^\infty$.
\elmma
%%%%%%%%%%%%%%%%%%%%%%%%%%%%%%%%%%%%%%%%%%%%%%
\prf
Let us first look at the case $M\in \bbm_j^{\pm}$ where $1\leq j\leq \ell$.
In this case, one has $S_M^\pm\xi_\gamma=\xi_{\gamma'}$
where $\gamma'$ is given by
\[
 \gamma'_i=\begin{cases}
            \gamma_i+1 & \text{ if } \begin{cases}
              m_i=1\text{ and }m_{i+1}=\ell+1-i,\cr
              m_{2\ell+2-i}=1\text{ and }m_{2\ell+3-i}=\ell+2-(2\ell+3-i),%\cr
%               i=j \text{ and } m_j=1,
                                   \end{cases}\cr
   &\cr
          \gamma_i-1 & \text{ if }\begin{cases}
              m_i=\ell+2-i\text{ and }m_{i+1}=1,\cr
               m_i=\ell+2-i\text{ and }i=j,\cr
              m_{2\ell+2-i}=\ell+2-(2\ell+2-i)\text{ and }m_{2\ell+3-i}=1,
                                   \end{cases}\cr
   &\cr
            \gamma_i & \text{ otherwise.}
           \end{cases}
\]
Note that since $1\leq j\leq\ell$, we have 
$\gamma'_{\ell+1}=\gamma_{\ell+1}$,
and $\eta(\gamma')-\eta(\gamma)$ depends just
on $M$ and not on $\gamma$.
Therefore $US_M^\pm U^*$
is a constant times simple tensor product of
shift operators. Thus in this case 
$U S_M^\pm U^*\in\scrt_\ell^\infty\otimes I\otimes \scrt_\ell^\infty
                    \seq \cle_\Sigma^\infty$.

Next we look at the case $M\in\bbm_{\ell+1}^\pm$.
In this case, define $\gamma'$ and $\gamma''$ as follows:
\[
 \gamma'_i=\begin{cases}
            \gamma_i+1 & \text{ if } \begin{cases}
              m_i=1\text{ and }m_{i+1}=\ell+1-i,\cr
              m_{2\ell+2-i}=1\text{ and }m_{2\ell+3-i}=\ell+2-(2\ell+3-i),
                                   \end{cases}\cr
   &\cr
          \gamma_i-1 & \text{ if }\begin{cases}
              m_i=\ell+2-i\text{ and }m_{i+1}=1,\cr
              m_{2\ell+2-i}=\ell+2-(2\ell+2-i)\text{ and }m_{2\ell+3-i}=1,\cr
               i=\ell+1
                                   \end{cases}\cr
   &\cr
            \gamma_i & \text{ otherwise.}
           \end{cases}
\]
\[
 \gamma''_i=\begin{cases}
            \gamma_i+1 & \text{ if } \begin{cases}
              m_i=1\text{ and }m_{i+1}=\ell+1-i,\cr
              m_{2\ell+2-i}=1\text{ and }m_{2\ell+3-i}=\ell+2-(2\ell+3-i),\cr
              i=\ell,
                                   \end{cases}\cr
   &\cr
          \gamma_i-1 & \text{ if }\begin{cases}
              m_i=\ell+2-i\text{ and }m_{i+1}=1,\cr
              m_{2\ell+2-i}=\ell+2-(2\ell+2-i)\text{ and }m_{2\ell+3-i}=1,\cr
               i=\ell+1
                                   \end{cases}\cr
   &\cr
            \gamma_i & \text{ otherwise.}
           \end{cases}
\]
Then one has
\[
 S_M^\pm\xi_\gamma=\begin{cases}
                    \xi_{\gamma'} & \text{ if } \gamma_{\ell+1}\leq 0,\cr
                   \xi_{\gamma''} & \text{ if } \gamma_{\ell+1}>0.
                   \end{cases}
\]
Therefore in this case, one will
have
$U S_M^\pm U^*\in\scrt_\ell^\infty\otimes \cle^\infty \otimes \scrt_\ell^\infty
                    \seq \cle_\Sigma^\infty$.

\qed

We will next take a closer look at the operators $T_M^\pm$.
For this, we will need to compute the quantitites involved
in equations~(\ref{eq:tmplus}) and (\ref{eq:tmminus})
more precisely than we have done earlier.
We start with the computation of $\kappa$.
From equation~(\ref{eq:kappa0}), we get
\bean
\psi(\bldr^{n,k}) &=& -\frac{\ell}{2}\left(n+k+(\ell-1)k\right)
                         +\frac{\ell(\ell+1)}{2} k\\
&=& -\frac{\ell}{2}(n-k).
\eean
Therefore
\be\label{eq:kappa2}
 \psi(\bldr^{n,k})-\psi(N_{1,1}(\bldr^{n,k}))=\psi(\bldr^{n,k})-\psi(\bldr^{n+1,k})
   = \frac{\ell}{2},
\ee
\be\label{eq:kappa3}
 \psi(\bldr^{n,k})-\psi(N_{0,1}(\bldr^{n,k}))=\psi(\bldr^{n,k})-\psi(\bldr^{n,k-1})
   = \frac{\ell}{2}.
\ee
Let us write $\lambda=(n+k,k,\ldots,k,0)$.
We will next compute $d_\lambda$, where $d_\lambda$ is given by (\ref{eq:kappa0}).
One has $d_\lambda=\sum_{\blds}q^{2\psi(\blds)}$
where the sum is over all those $\blds$ for which the top row is $\lambda$.
Such an $\blds$ is of the form~(\ref{eq:GTnotation})
and one has 
\[
\psi(\blds)=-\half\ell(n+\ell k)+\half(\ell-1)(\ell-2)k+ \sum_{i=2}^\ell (c_i+d_i) + d_{\ell+1}.
\]
Thus we have
\[
 d_\lambda=q^{-\ell(n+k)-2(\ell-1)k}
   \sum_{\substack{k\leq c_\ell\leq c_{\ell-1}\leq\ldots\leq c_2\leq n+k \\
         0\leq d_2\leq d_3\leq\ldots\leq d_\ell\leq k \\ d_\ell\leq d_{\ell+1}\leq c_\ell}}
   q^{2\left(\sum_{i=2}^\ell (c_i+d_i) + d_{\ell+1}\right)}
\]
Now for any $x$, we have
\bea
\lefteqn{ \sum_{\substack{k\leq c_\ell\leq c_{\ell-1}\leq\ldots\leq c_2\leq n+k \\
         0\leq d_2\leq d_3\leq\ldots\leq d_\ell\leq k \\ d_\ell\leq d_{\ell+1}\leq c_\ell}}
   x^{\left(\sum_{i=2}^\ell (c_i+d_i) + d_{\ell+1}\right)} }\nonumber\\
&=&\left(\sum_{k\leq d_{\ell+1}\leq c_\ell\leq c_{\ell-1}\leq\ldots\leq c_2\leq n+k}
   x^{\left(\sum_{i=2}^\ell c_i + d_{\ell+1}\right)}\right)
\left(\sum_{0\leq d_2\leq d_3\leq\ldots\leq d_\ell\leq k}
   x^{\left(\sum_{i=2}^\ell d_i \right)}\right)\nonumber\\
&& + \left(\sum_{k\leq c_\ell\leq c_{\ell-1}\leq\ldots\leq c_2\leq n+k}
   x^{\left(\sum_{i=2}^\ell c_i \right)}\right)
\left(\sum_{0\leq d_2\leq d_3\leq\ldots\leq d_\ell\leq d_{\ell+1}< k}
   x^{\left(\sum_{i=2}^\ell d_i + d_{\ell+1}\right)}\right).
\eea
If we now use the identity
\[
 \sum_{k\leq t_1\leq t_2\leq\ldots\leq t_j\leq n}
   x^{\left(\sum_{i=1}^j t_i \right)}
= x^{jk}\prod_{i=1}^j\left(\frac{1-x^{n-k+i}}{1-x^i}\right),
\]
we get
\bean
\lefteqn{ \sum_{\substack{k\leq c_\ell\leq c_{\ell-1}\leq\ldots\leq c_2\leq n+k \\
         0\leq d_2\leq d_3\leq\ldots\leq d_\ell\leq k \\ d_\ell\leq d_{\ell+1}\leq c_\ell}}
   x^{\left(\sum_{i=2}^\ell (c_i+d_i) + d_{\ell+1}\right)} }\\
&=&
x^{\ell k}\prod_{i=1}^\ell\left(\frac{1-x^{n+i}}{1-x^i}\right)
       \prod_{i=1}^{\ell-1}\left(\frac{1-x^{k+i}}{1-x^i}\right)
+
x^{(\ell-1)k}\prod_{i=1}^{\ell-1}\left(\frac{1-x^{n+i}}{1-x^i}\right)
    \prod_{i=1}^\ell\left(\frac{1-x^{k-1+i}}{1-x^i}\right)\\
&=&
 x^{(\ell-1)k}\prod_{i=1}^{\ell-1}\left(\frac{1-x^{n+i}}{1-x^i}\right)
    \prod_{i=1}^{\ell-1}\left(\frac{1-x^{k+i}}{1-x^i}\right)  %\prod_{i=1}^\ell \frac{1}{1-x^i}
\frac{1}{1-x^\ell}\left(x^k(1-x^{n+\ell})+1-x^k\right)\\
&=&
  x^{(\ell-1)k}\prod_{i=1}^{\ell-1}\left(\frac{1-x^{n+i}}{1-x^i}\right)
    \prod_{i=1}^{\ell-1}\left(\frac{1-x^{k+i}}{1-x^i}\right)
   \left(\frac{1-x^{n+k+\ell}}{1-x^\ell}\right).
\eean
Thus
\be
d_\lambda^\half=q^{-\frac{\ell(n+k)}{2}}
  \prod_{i=1}^{\ell-1}\left(\frac{Q(n+i)}{Q(i)}\frac{Q(k+i)}{Q(i)}\right)
   \frac{Q(n+k+\ell)}{Q(\ell)}.
\ee
Write  
\[
 \lambda'=(n+1+k,k,\ldots,k,0),  \qquad
   \lambda''=(n+k-1,k-1,\ldots,k-1,0).
\]
Then one has
\bean
 d_\lambda^\half d_{\lambda'}^{-\half} &=&
   q^{\ell/2}\frac{Q(n+1)}{Q(n+\ell)}\frac{Q(n+k+\ell)}{Q(n+k+\ell+1)},\\
 d_\lambda^\half d_{\lambda''}^{-\half} &=&
   q^{-\ell/2}\frac{Q(k+\ell-1)}{Q(k)}\frac{Q(n+k+\ell)}{Q(n+k+\ell-1)}.
\eean
Combining these with (\ref{eq:kappa2}) and (\ref{eq:kappa3}),
we get
\bea
\kappa(\bldr^{n,k},N_{1,1}(\bldr^{n,k})) &=& 
    q^{\ell}\frac{Q(n+1)}{Q(n+\ell)}\frac{Q(n+k+\ell)}{Q(n+k+\ell+1)},\label{eq:kappa4}\\
\kappa(\bldr^{n,k},N_{0,1}(\bldr^{n,k})) &=& 
    \frac{Q(k+\ell-1)}{Q(k)}\frac{Q(n+k+\ell)}{Q(n+k+\ell-1)}.\label{eq:kappa5}
\eea

%%%%%%%%%%%%%%%%%%%%%%%%%%%%%%%%%%%%%%%%%%%%%
\blmma\label{lem:tmplusignore}
Let $M\in\bbm_j^+$ and $T_M^+$ be as in equation~(\ref{eq:tmplus}). Then
$U T_M^+ U^*\in OP_{D_\ell}^{-\infty}\otimes\scrt_\ell^\infty$
if $j\leq \ell$ or if $j=\ell+1$ and $M\neq N_\ell$.
\elmma
%%%%%%%%%%%%%%%%%%%%%%%%%%%%%%%%%%%%%%%%%%%%%
\prf
From lemma~\ref{lem:cgcoeff} and equations (\ref{eq:tmplus}) and (\ref{eq:kappa4}),
we get, for $M=(m_1,\ldots,m_j)\in \bbm_j^+$,
\bea
T_M^+ e_{\bldr^{n,k}\blds} &=& 
  \sgn(M)q^{\ell-j+1+C(\bldr^{n,k},N_{1,1})+B(N_{1,1})+C(\blds,M)+B(M)}
    \frac{Q(n+1)}{Q(n+\ell)}\frac{Q(n+k+\ell)}{Q(n+k+\ell+1)} \nonumber\\
&\times &  
  L'(\bldr^{n,k},1,1)\left(\prod_{a=1}^{j-1}L(\blds,a,m_a,m_{a+1})\right)L'(\blds,j,m_j)
 \;e_{\bldr^{n,k}\blds}.
       \label{eq:tmplus1}
\eea
Since $C(\bldr^{n,k},N_{1,1})=k$ and $B(N_{1,1})=0$, we get
\[
 T_M^+ e_{\bldr^{n,k}\blds} =\sgn(M)q^{\ell-j+1+B(M)+k+C(\blds,M)}\phi(\blds,M)e_{\bldr^{n,k}\blds},
\]
with $\phi(\blds,M)$  a product of terms of the form
$Q(\psi(\gamma))^{\pm 1}$ where $\psi(\gamma)=|\gamma_A|+c(\gamma_{\ell+1})_\pm +m$
for some subset $A\seq\Sigma$, $c\in\{0,1\}$ and some integer $m$ that does not depend
on $\blds$. 
Therefore
\[
 UT_M^+ U^* e_\gamma =\sgn(M)q^{\ell-j+1+B(M)+k+C(\blds,M)}\phi(\blds,M)e_\gamma,
\]
where
$k$ and $\blds$ are given by equations~(\ref{eq:ustar2}--\ref{eq:ustar4}).
Since $\phi(\blds,M)$  a product of terms of the form
$Q(\psi(\gamma))^{\pm 1}$,
it follows from lemma~\ref{lem:neglect4} that the operator $e_\gamma\mapsto \phi(\blds,M)e_\gamma$
is in $\cle_{\Sigma}^\infty$. Next look at the operator
$e_\gamma\mapsto q^{k+C(\blds,M)}e_\gamma$.
Assume that there is some $i\leq j$ such that $m_i\neq 1$.
Let $p=\min\{2\leq i\leq j: m_i\neq 1\}$.
Then $C(\blds,M)\geq H_{p-1,1}(\blds)\geq (\gamma_{\ell+1})_-$.
Therefore
\[
 k+C(\blds,M)  \geq k+(\gamma_{\ell+1})_- = \sum_{i=1}^\ell\gamma_i+|\gamma_{\ell+1}|.
\]
Hence $UT_M^+ U^*\in  OP_{D_\ell}^{-\infty}\otimes\scrt_\ell^\infty$.
Next assume that $j\leq \ell$ and $m_i=1$ for all $i\leq j$.
In this case,
$C(\blds,M)\geq H_{j,1}(\blds)\geq (\gamma_{\ell+1})_-$.
Therefore again we have
\[
k+C(\blds,M)  \geq k+(\gamma_{\ell+1})_- = \sum_{i=1}^\ell\gamma_i+|\gamma_{\ell+1}|.
\]
and hence
$UT_M^+ U^*\in  OP_{D_\ell}^{-\infty}\otimes\scrt_\ell^\infty$.
Combining the two cases, we have the required result.
\qed

%%%%%%%%%%%%%%%%%%%%%%%%%%%%%%%%%%%%%%%%%%%%%
\blmma\label{lem:tmminusignore}
Let $M\in\bbm_j^-$ and $T_M^-$ be as in equation~(\ref{eq:tmminus}). Then
$U T_M^- U^*\in OP_{D_\ell}^{-\infty}\otimes\scrt_\ell^\infty$
if $M\neq N_{0,j}$.
\elmma
%%%%%%%%%%%%%%%%%%%%%%%%%%%%%%%%%%%%%%%%%%%%%
\prf
From lemma~\ref{lem:cgcoeff} and equations (\ref{eq:tmminus}) and (\ref{eq:kappa5}),
we get, for $M=(m_1,\ldots,m_j)\in \bbm_j^-$,
\bea
T_M^- e_{\bldr^{n,k}\blds} &=& \sgn(M)q^{-j+1+C(\bldr^{n,k},N_{0,1})+B(N_{0,1})+C(\blds,M)+B(M)}
   \frac{Q(k+\ell-1)}{Q(k)}\frac{Q(n+k+\ell)}{Q(n+k+\ell-1)} \nonumber\\
&\times &  
  L'(\bldr^{n,k},1,\ell+1)\left(\prod_{a=1}^{j-1}L(\blds,a,m_a,m_{a+1})\right)L'(\blds,j,m_j)
 e_{\bldr^{n,k}\blds}.
       \label{eq:tmminus1}
\eea
Since $C(\bldr^{n,k},N_{0,1})=0$ and $B(N_{0,1})=0$, we get
\[
 T_M^- e_{\bldr^{n,k}\blds} =\sgn(M)q^{-j+1+C(\blds,M)+B(M)}\phi(\blds,M)e_{\bldr^{n,k}\blds},
\]
with $\phi(\blds,M)$  a product of terms of the form
$Q(\psi(\gamma))^{\pm 1}$ where $\psi(\gamma)=|\gamma_A|+c(\gamma_{\ell+1})_\pm +m$
for some subset $A\seq\Sigma$, $c\in\{0,1\}$ and some integer $m$ that does not depend
on $\blds$. 
Therefore
\[
 UT_M^- U^* e_\gamma =\sgn(M)q^{-j+1+C(\blds,M)+B(M)}\phi(\blds,M)e_\gamma,
\]
where
$k$ and $\blds$ are given by equations~(\ref{eq:ustar2}--\ref{eq:ustar4}).
As in the proof of lemma~\ref{lem:tmplusignore},
it is now enough to prove that 
$C(\blds,M)\geq \sum_{i=1}^\ell\gamma_i +|\gamma_{\ell+1}|$.
Now assume that $m_i=1$ for some $i\leq\ell$.
Let $p=\min\{2\leq i\leq j: m_i=1\}$. Then $p\leq \ell$.
We then have 
% $B(M)> (p-2)+2(\ell+2-p+1-1)+1=2\ell-p+3$, and
\bean
C(\blds,M) &\geq&\sum_{i=1}^{p-2}H_{i,\ell+1-i}(\blds)
      +H_{p-1,1}(\blds)+H_{p-1,\ell+2-p}(\blds)+V_{p-1,\ell+2-p}(\blds)\\
&\geq&\sum_{i=1}^{p-2}\gamma_i +(\gamma_{\ell+1})_- + \gamma_{p-1}
    +\left(\sum_{i=1}^\ell\gamma_i +(\gamma_{\ell+1})_{+} -\sum_{i=1}^{p-1}\gamma_i\right)\\
&=& \sum_{i=1}^\ell\gamma_i +|\gamma_{\ell+1}|.
\eean
So the result follows.
\qed

%%%%%%%%%%%%%%%%%%%%%%%%%%%%%%%%%%%%%%%%%%
\brmrk\label{rem:dandrea}
As mentioned in the beginning of this subsection,
weaker versions of the two lemmas above have been proved
by D'Andrea in \cite{dan-2008a}. In our notation,  he
proves that the part of $Z_{j,q}$ that be ignored is
of the order $q^k=q^{\sum_{i=1}^\ell\gamma_i +(\gamma_{\ell+1})_{+}}$,
whereas we prove here that one can actually ignore terms
of a slightly higher order, namely 
$q^{\sum_{i=1}^\ell\gamma_i +|\gamma_{\ell+1}|}$,
which makes it possible to compute $Z_{j,q}$ modulo
the ideal $OP_{D_\ell}^{-\infty}\otimes\scrt_\ell^\infty$.
\ermrk
%%%%%%%%%%%%%%%%%%%%%%%%%%%%%%%%%%%%%%%%%%

%%%%%%%%%%%%%%%%%%%%%%%%%%%%%%%%%%%%%%%%%%
\blmma \label{lem:xj}
Define operators $X_j$ on $L_2(S_q^{2\ell+1})$ by
\be\label{eq:zjq3}
e_{\bldr^{n,k},\blds}  \mapsto    \begin{cases}
  (-1)^{j-1} q^{d_j}Q(d_{j+1}-d_j)  e_{\bldr^{n,k-1},N_{0,j}(\blds)} & \text{ if }1\leq j\leq\ell-1,\cr
    (-1)^{\ell-1} q^{d_\ell}Q(d_{\ell+1}-d_\ell)Q(k-d_\ell)  e_{\bldr^{n,k-1},N_{0,\ell}(\blds)}
                   & \text{ if } j=\ell. 
                            \end{cases}
\ee
Then one has
\[
UZ_{j,q}^* U^*- UX_j U^*\in OP_{D_\ell}^{-\infty}\otimes \scrt_\ell^\infty.
\]
\elmma
%%%%%%%%%%%%%%%%%%%%%%%%%%%%%%%%%%%%%%%%%%
\prf
In view of the two forgoing lemmas, it is  enough to show that
\bea
U S_{N_{0,j}}^- T_{N_{0,j}}^- U^* - U X_{j} U^* &\in & 
   OP_{D_\ell}^{-\infty}\otimes\scrt_\ell^\infty,\quad\text{for }1\leq j\leq\ell,
              \label{eq:equality1} 
\eea
Let us first look at the case $1\leq j\leq \ell-1$.
Observe that 
\[
\sgn(N_{0,j})=(-1)^{j-1}, \quad
C(\bldr^{n,k},N_{0,1})=0=B(N_{0,1}),\quad
C(\blds,N_{0,j})=d_j,\quad
B(N_{0,j})=j-1.
\]
Therefore from (\ref{eq:tmminus1}), we get
\bea
 UT_{N_{0,j}}^- U^* e_\gamma
   &=&(-1)^{j-1}q^{d_j}
   \frac{Q(k+\ell-1)}{Q(k)}\frac{Q(n+k+\ell)}{Q(n+k+\ell-1)} L'(\bldr^{n,k},1,\ell+1)\nonumber\\
&\times &  
  \left(\prod_{a=1}^{j-1}L(\blds,a,\ell+2-a,\ell+1-a)\right)L'(\blds,j,\ell+2-j)\;
 e_\gamma.\label{eq:tmminus2}
\eea
From (\ref{eq:comp_row2b}), one gets
\bea
L'(\bldr^{n,k},1,\ell+1) 
&=& \left(\prod_{i=2}^{\ell}\frac{Q(|k-0-i+\ell+1-1|)}{Q(|k-0-i+\ell+1|)}\right)
            \frac{Q(|k-0-1+\ell+1-1|)}{Q(|n+k-0-1+\ell+1|)} \nonumber\\
&=& 
\left(\prod_{i=2}^{\ell}\frac{Q(k+\ell-i)}{Q(k+\ell-i+1)}\right)\frac{
Q(k+\ell-1)}{Q(n+k+\ell)}\nonumber\\
&=& \frac{Q(k)}{Q(n+k+\ell)}.\label{eq:Lprime1}
\eea
Similarly, from (\ref{eq:comp_row1c}) one gets, for $1\leq a\leq \ell-1$,
\bea
\lefteqn{L(\blds,a,\ell+2-a,\ell+1-a)} \nonumber\\ 
&=&
\prod_{i=1}^{\ell+1-a}\frac{Q(|s_{a,i}-s_{a+1,\ell+1-a}-i+\ell+1-a|)}{Q(|s_{a,i}
-s_{a,\ell+2-a}-i+\ell+2-a|)}
\prod_{i=1}^{\ell-a}\frac{Q(|s_{a+1,i}-s_{a,\ell+2-a}-i+\ell+2-a-1|)}{Q(|s_{a+1,
i}-s_{a+1,\ell+1-a}-i+\ell+1-a-1|)} \nonumber\\
&=&\frac{Q(c_a-d_{a+1}+\ell-a)}{Q(c_{a}-d_{a}+\ell+1-a)}
       \frac{Q(c_{a+1}-d_{a}+\ell-a)}{Q(c_{a+1}-d_{a+1}+\ell-a-1)} \nonumber\\
&& \times
\prod_{i=2}^{\ell+1-a}\frac{Q(k-d_{a+1}-i+\ell+1-a)}{Q(k-d_{a}-i+\ell+2-a)}
     \prod_{i=2}^{\ell-a}\frac{Q(k-d_{a}-i+\ell+1-a)}{Q(k-d_{a+1}-i+\ell-a)} \nonumber\\
&=&\frac{Q(c_a-d_{a+1}+\ell-a)}{Q(c_{a}-d_{a}+\ell+1-a)}
       \frac{Q(c_{a+1}-d_{a}+\ell-a)}{Q(c_{a+1}-d_{a+1}+\ell-a-1)} \nonumber\\
&& \times
\prod_{i=1}^{\ell-a}\frac{Q(k-d_{a+1}-i+\ell-a)}{Q(k-d_{a}-i+\ell+1-a)}
     \prod_{i=2}^{\ell-a}\frac{Q(k-d_{a}-i+\ell+1-a)}{Q(k-d_{a+1}-i+\ell-a)} \nonumber\\
&=& \frac{Q(c_a-d_{a+1}+\ell-a)}{Q(c_{a}-d_{a}+\ell+1-a)}
       \frac{Q(c_{a+1}-d_{a}+\ell-a)}{Q(c_{a+1}-d_{a+1}+\ell-a-1)}
            \frac{Q(k-d_{a+1}+\ell-a-1)}{Q(k-d_{a}+\ell-a)},\label{eq:Lnoprime1}
\eea
and from (\ref{eq:comp_row2b}), for $j\leq\ell-1$,
\bean
\lefteqn{L'(\blds,j,\ell+2-j)}\\ 
&=&
\prod_{i=1}^{\ell+1-j}\frac{Q(|s_{j+1,i}-s_{j,\ell+2-j}-i+\ell+2-j-1|)}{Q(|s_{j,
i}-s_{j,\ell+2-j}-i+\ell+2-j|)}\\
&=&  \frac{Q(c_{j+1}-d_j+\ell-j)}{Q(c_j-d_{j}+\ell+1-j)}
\left(\prod_{i=2}^{\ell-j}\frac{Q(k-d_j+\ell+1-j-i)}{Q(k-d_{j}+\ell+2-j-i)}
\right)
          \frac{Q(d_{j+1}-d_j)}{Q(k-d_{j}+1)}\\
&=&
\frac{Q(c_{j+1}-d_j+\ell-j)}{Q(c_j-d_{j}+\ell+1-j)}\frac{Q(d_{j+1}-d_j)}{Q(k-d_{
j}+\ell-j)}
\eean
From the above two equations, we get
\bean
\lefteqn{\left(\prod_{a=1}^{j-1}L(\blds,a,\ell+2-a,\ell+1-a)\right)L'(\blds,j,\ell+2-j)}\\
 &=& 
\frac{Q(d_{j+1}-d_j)}{Q(k+\ell-1)}
  \left(\prod_{a=1}^{j-1}\frac{Q(c_a-d_{a+1}+\ell-a)}{Q(c_{a+1}-d_{a+1}+\ell-a-1)}\right)
\left(\prod_{a=1}^{j}\frac{Q(c_{a+1}-d_a+\ell-a)}{Q(c_a-d_a+\ell+1-a)}
\right).
\eean
Now substituting all these in equation~(\ref{eq:tmminus2}), we get
\bea
 UT_{N_{0,j}}^- U^* e_\gamma
   &=&(-1)^{j-1}q^{d_j}
   \frac{Q(d_{j+1}-d_j)}{Q(n+k+\ell-1)}
\left(\prod_{a=1}^{j-1}\frac{Q(c_a-d_{a+1}+\ell-a)}{Q(c_{a+1}-d_{a+1}+\ell-a-1)}\right)
\nonumber \\
&&\times
 \left(\prod_{a=1}^{j}\frac{Q(c_{a+1}-d_a+\ell-a)}{Q(c_a-d_a+\ell+1-a)}
\right)\;e_\gamma.\label{eq:tmminus3}
\eea
Now note that for $1\leq a\leq j-1$,
\[
 d_j+c_a-d_{a+1}+\ell-a  \geq \sum_{i=1}^\ell \gamma_i+|\gamma_{\ell+1}|,\qquad
d_j+c_{a+1}-d_{a+1}+\ell-a-1 \geq \sum_{i=1}^\ell \gamma_i+|\gamma_{\ell+1}|,
\]
and for $1\leq a\leq j$,
\[
 d_j+c_{a+1}-d_a+\ell-a \geq \sum_{i=1}^\ell \gamma_i+|\gamma_{\ell+1}|,\qquad
d_j+c_a-d_a+\ell+1-a \geq \sum_{i=1}^\ell \gamma_i+|\gamma_{\ell+1}|,
\]
and $d_j+n+k+\ell-1\geq  \sum_{i=1}^\ell \gamma_i+|\gamma_{\ell+1}|$.
Therefore by using lemmas~\ref{lem:neglect2} and \ref{lem:neglect3},
we can write, modulo an operator in $OP^{-\infty}_{D_\ell}\otimes\scrt_\ell^\infty$,
\[
UT_{N_{0,j}}^- U^* e_\gamma
   = (-1)^{j-1}q^{d_j}Q(d_{j+1}-d_j)\;e_\gamma.
 \]
Using equation~(\ref{eq:smminus}), we get
\[
 U S_{N_{0,j}}^- U^* e_\gamma
   = (-1)^{j-1} e_{\gamma'},
\]
where 
\[
\gamma'_i=\begin{cases}
                \gamma_i & \text{ if }i\neq j,\cr
              \gamma_i-1 & \text{ if }i=j.
               \end{cases}
\]
Observe also that
\[
 UX_j U^* e_\gamma=q^{d_j}Q(d_{j+1}-d_j)\;e_{\gamma'},
\]
where $\gamma'$ is as above. Therefore
we get (\ref{eq:equality1}) for $j\leq\ell-1$.

In the case $j=\ell$, one has
\[
 L'(\blds,\ell,2)=
   \frac{Q(|s_{\ell+1,1}-s_{\ell,2}|)}{Q(|s_{\ell,1}-s_{\ell,2}+1|)}
   =\frac{Q(d_{\ell+1}-d_{\ell})}{Q(c_\ell-d_{\ell}+1)}.
\]
and as a result, one has
\bean
\lefteqn{\left(\prod_{a=1}^{\ell-1}L(\blds,a,\ell+2-a,\ell+1-a)\right)L'(\blds,\ell,2)}\\
 &=& 
\frac{Q(k-d_\ell)}{Q(k+\ell-1)}
  \left(\prod_{a=1}^{\ell-1}\frac{Q(c_a-d_{a+1}+\ell-a)}{Q(c_{a+1}-d_{a+1}+\ell-a-1)}\right)
\left(\prod_{a=1}^{\ell}\frac{Q(c_{a+1}-d_a+\ell-a)}{Q(c_a-d_a+\ell+1-a)}
\right).
\eean
As before, substituting all these in equation~(\ref{eq:tmminus2}), one gets
\bea
 UT_{N_{0,\ell}}^- U^* e_\gamma
   &=&(-1)^{\ell-1}q^{d_\ell}
   \frac{Q(k-d_\ell)}{Q(n+k+\ell-1)}
\left(\prod_{a=1}^{\ell-1}\frac{Q(c_a-d_{a+1}+\ell-a)}{Q(c_{a+1}-d_{a+1}+\ell-a-1)}\right)
\nonumber \\
&&\times
 \left(\prod_{a=1}^{\ell}\frac{Q(c_{a+1}-d_a+\ell-a)}{Q(c_a-d_a+\ell+1-a)}
\right)\;e_\gamma.\label{eq:tmminus4}
\eea
Application of lemmas~\ref{lem:neglect2} and \ref{lem:neglect3},
now enable us to write the following equality modulo an operator in $OP^{-\infty}_{D_\ell}\otimes\scrt_\ell^\infty$:
\[
UT_{N_{0,\ell}}^- U^* e_\gamma
   = (-1)^{\ell-1}q^{d_\ell}Q(k-d_\ell)Q(d_{\ell+1}-d_\ell)\;e_\gamma.
 \]
Using equation~(\ref{eq:smminus}), we get
\[
 U S_{N_{0,\ell}}^- U^* e_\gamma
   = (-1)^{\ell-1} e_{\gamma'},
\]
where 
\[
\gamma'_i=\begin{cases}
                \gamma_i & \text{ if }i\neq \ell,\cr
              \gamma_i-1 & \text{ if }i=\ell.
               \end{cases}
\]
Observe also that
\[
 UX_\ell U^* e_\gamma=q^{d_\ell}Q(k-d_\ell)Q(d_{\ell+1}-d_\ell)\;e_{\gamma'},
\]
where $\gamma'$ is as above. Therefore
we get (\ref{eq:equality1}) for $j=\ell$.
\qed

%%%%%%%%%%%%%%%%%%%%%%%%%%%%%%%%%%%%%%%%%%%%%%
\blmma\label{lem:zjq2}
Let $X_j$ be as in lemma~\ref{lem:xj}. 
Then  for $1\leq j\leq\ell$, one has
$UX_j U^*-Y_{j,q}^*\otimes I\in OP_{D_\ell}^{-\infty}\otimes \scrt_\ell^\infty$.
\elmma
%%%%%%%%%%%%%%%%%%%%%%%%%%%%%%%%%%%%%%
\prf
It follows from equations~(\ref{eq:ustar1}--\ref{eq:ustar4}) 
that  for $j\leq \ell-1$, one in fact has
$UX_j U^*-Y_{j,q}^*\otimes I=0$.
For $j=\ell$, one has
\[
\left(UX_j U^*-Y_{j,q}^*\otimes I\right)e_\gamma
   = \left(q^{\sum_{i=1}^{\ell-1}\gamma_i}
 Q\left(\gamma_\ell+(\gamma_{\ell+1})_-\right)
   Q\left(\gamma_\ell+(\gamma_{\ell+1})_+\right)\right)e_{\hat{\gamma}},
\]
where $\hat{\gamma}_i=\gamma_i-1$ if $i=\ell$ and $\hat{\gamma}_i=\gamma_i$
for all other $i$.
Thus 
\[
\left|UX_j U^*-Y_{j,q}^*\otimes I\right|\in OP_{D_\ell}^{-\infty}\otimes \scrt_\ell^\infty
\qquad  
\sgn \left(UX_j U^*-Y_{j,q}^*\otimes I\right)\in \cle^\infty_\Sigma.
\]
Therefore $UX_j U^*-Y_{j,q}^*\otimes I\in OP_{D_\ell}^{-\infty}\otimes \scrt_\ell^\infty$.
\qed

From the two lemmas above (lemmas~\ref{lem:xj} and \ref{lem:zjq2}),
it follows that for $1\leq j\leq \ell$,
one has $UZ_{j,q}^* U^*\in C^\infty(S_q^{2\ell+1})$.
Thus we now need only to take care of the case $j=\ell+1$.

%%%%%%%%%%%%%%%%%%%%%%%%%%%%%%%%%%%%%%%%%%%
\blmma\label{lem:proofpart3}
$UZ_{\ell+1,q}^*U^*\in C^\infty(S_q^{2\ell+1})$.
\elmma
%%%%%%%%%%%%%%%%%%%%%%%%%%%%%%%%%%%%%%%%%%%
\prf
Using lemmas~\ref{lem:tmplusignore} and \ref{lem:tmminusignore},
it is enough to show that
\bea
U(S_{N_\ell}^+T_{N_\ell}^+ + S_{N_{0,\ell+1}}^-
       T_{N_{0,\ell+1}}^-)U^*
 &\in &
C^\infty(S_q^{2\ell+1}).
\label{eq:equality2}
\eea

From (\ref{eq:tmplus1}), we get
\be
T_{N_\ell}^+ e_{\bldr^{n,k}\blds} = 
  q^{k}\frac{Q(n+1)}{Q(n+\ell)}\frac{Q(n+k+\ell)}{Q(n+k+\ell+1)}   
  L'(\bldr^{n,k},1,1)\left(\prod_{a=1}^{\ell}L(\blds,a,1,1)\right)\;
 e_{\bldr^{n,k}\blds}.
       \label{eq:tmplus2}
\ee
From (\ref{eq:comp_row1c}), we get
for $1\leq a\leq \ell-1$,
\bean
\lefteqn{L(\blds,a,1,1)}\\ 
&=&
\prod_{i=2}^{\ell+2-a}\frac{Q(|s_{a,i}-s_{a+1,1}-i+1|)}{Q(|s_{a,i}-s_{a,1}-i+1|)
}
\prod_{i=2}^{\ell+1-a}\frac{Q(|s_{a+1,i}-s_{a,1}-i+1-1|)}{Q(|s_{a+1,i}-s_{a+1,1}
-i+1-1|)}\\
&=&\prod_{i=2}^{\ell+1-a}\frac{Q(c_{a+1}-k+i-1)}{Q(c_{a}-k+i-1)}
     \prod_{i=2}^{\ell-a}\frac{Q(c_{a}-k+i)}{Q(c_{a+1}-k+i)}\\
&& \times \frac{Q(c_{a+1}-d_{a}+\ell+1-a)}{Q(c_{a}-d_{a}+\ell+1-a)}
       \frac{Q(c_{a}-d_{a+1}+\ell+1-a)}{Q(c_{a+1}-d_{a+1}+\ell+1-a)}\\
&=&\prod_{i=1}^{\ell-a}\frac{Q(c_{a+1}-k+i)}{Q(c_{a}-k+i)}
     \prod_{i=2}^{\ell-a}\frac{Q(c_{a}-k+i)}{Q(c_{a+1}-k+i)}\\
&& \times \frac{Q(c_{a+1}-d_{a}+\ell+1-a)}{Q(c_{a}-d_{a}+\ell+1-a)}
       \frac{Q(c_{a}-d_{a+1}+\ell+1-a)}{Q(c_{a+1}-d_{a+1}+\ell+1-a)}\\
&=& \frac{Q(c_{a+1}-k+1)}{Q(c_{a}-k+1)}
     \frac{Q(c_{a+1}-d_{a}+\ell+1-a)}{Q(c_{a}-d_{a}+\ell+1-a)}
       \frac{Q(c_{a}-d_{a+1}+\ell+1-a)}{Q(c_{a+1}-d_{a+1}+\ell+1-a)},
\eean
and for $a=\ell$,
\[
L(\blds,\ell,1,1)=\frac{Q(|s_{\ell,2}-s_{\ell+1,1}-2+1|)}{Q(|s_{\ell,2}-s_{\ell,
1}-2+1|)}
   =\frac{Q(d_{\ell+1}-d_{\ell}+1)}{Q(c_\ell-d_{\ell}+1)}.
\]
Also from (\ref{eq:comp_row2b}), we have
\bean
L'(\bldr^{n,k},1,1) 
&=&
\left(\frac{\prod_{i=1}^{\ell}Q(|k-n-k-i+1-1|)}{\prod_{i=2}^{\ell}Q(|k-n-k-i+1|)
}\right)
                  \frac{1}{Q(|0-n-k-\ell-1+1|)}\\
&=& 
\left(\prod_{i=2}^{\ell}\frac{Q(n+i)}{Q(n+i-1)}\right)\frac{Q(n+1)}{Q(n+k+\ell)}
\\
&=& \frac{Q(n+\ell)}{Q(n+k+\ell)}.
\eean
Plugging these in equation~(\ref{eq:tmplus2}) and using (\ref{eq:ustar1}), we get
\bean
UT_{N_\ell}^+ U^* e_\gamma &=& 
  q^{k}\frac{Q(n+1)}{Q(n+\ell)}\frac{Q(n+k+\ell)}{Q(n+k+\ell+1)}   
  \frac{Q(n+\ell)}{Q(n+k+\ell)}  \\
&&
\times
 \left(\prod_{a=1}^{\ell-1}\frac{Q(c_{a+1}-k+1)}{Q(c_{a}-k+1)}
     \frac{Q(c_{a+1}-d_{a}+\ell+1-a)}{Q(c_{a}-d_{a}+\ell+1-a)}
       \frac{Q(c_{a}-d_{a+1}+\ell+1-a)}{Q(c_{a+1}-d_{a+1}+\ell+1-a)}\right)\\
&&
\times  \frac{Q(d_{\ell+1}-d_{\ell}+1)}{Q(c_\ell-d_{\ell}+1)}\;
 e_\gamma.
\eean
Thus as earlier, modulo an operator in $OP_{D_\ell}^{-\infty}\otimes\scrt_\ell^\infty$,
we have the equality
\be\label{eq:tplus}
 UT_{N_\ell}^+ U^* e_\gamma =q^{k}e_\gamma.
\ee

Next note that
$B(N_0)=\ell$, $C(\blds,N_0)=d_{\ell+1}$ and $\sgn(N_0)=(-1)^\ell$
so that  we get from~(\ref{eq:tmminus1})
\bea
T_{N_{0}}^- e_{\bldr^{n,k}\blds} &=& (-1)^\ell
    q^{d_{\ell+1}}
   \frac{Q(k+\ell-1)}{Q(k)}\frac{Q(n+k+\ell)}{Q(n+k+\ell-1)} \nonumber\\
&\times &  
  L'(\bldr^{n,k},1,\ell+1)\left(\prod_{a=1}^{\ell}L(\blds,a,\ell+2-a,\ell+1-a)\right)
 e_{\bldr^{n,k}\blds}.
       \label{eq:tmminus5}
\eea
Now using (\ref{eq:ustar1}), (\ref{eq:Lprime1}), (\ref{eq:Lnoprime1}) and
the fact that 
 \[
L(\blds,\ell,2,1)=\frac{Q(|s_{\ell,1}-s_{\ell+1,1}-1+\ell+1-\ell|)}{Q(|s_{\ell,1
}-s_{\ell,2}-1+\ell+2-\ell|)}
   =\frac{Q(c_\ell-d_{\ell+1})}{Q(c_\ell-d_{\ell}+1)},
\]
we get 
\bean
UT_{N_{0}}^- U^* e_\gamma &=& (-1)^\ell
    q^{d_{\ell+1}}
   \frac{Q(k+\ell-1)}{Q(k)}\frac{Q(n+k+\ell)}{Q(n+k+\ell-1)} \frac{Q(k)}{Q(n+k+\ell)}\nonumber\\
&\times &  
  \left(\prod_{a=1}^{\ell-1}\frac{Q(c_a-d_{a+1}+\ell-a)}{Q(c_{a}-d_{a}+\ell+1-a)}
       \frac{Q(c_{a+1}-d_{a}+\ell-a)}{Q(c_{a+1}-d_{a+1}+\ell-a-1)}
            \frac{Q(k-d_{a+1}+\ell-a-1)}{Q(k-d_{a}+\ell-a)}\right)\nonumber\\
&\times &\frac{Q(c_\ell-d_{\ell+1})}{Q(c_\ell-d_{\ell}+1)}
 \;e_\gamma.
\eean
Thus modulo $OP_{D_\ell}^{-\infty}\otimes\scrt_\ell^\infty$,
we have the equality
\be\label{eq:tminus}
 UT_{N_0}^- U^* e_\gamma =(-1)^\ell q^{d_{\ell+1}}e_\gamma.
\ee
Define operators $T^\pm$ on $L_2(S_q^{2\ell+1})$ by
\[
 T^+\xi_\gamma=q^k\xi_\gamma,\qquad
T^-\xi_\gamma=(-1)^\ell q^{d_{\ell+1}}\xi_\gamma.
\]
By equations~(\ref{eq:tplus}) and (\ref{eq:tminus}),
it is enough to look at the operators $S_{N_\ell}^+T^{+} + S_{N_0}^{-}T^{-}$.

Now observe that
\[
 S_{N_0}^- \xi_\gamma = \begin{cases}
                         \xi_{\gamma'} & \text{ if }\gamma_{\ell+1}>0,\cr
                       \xi_{\gamma''} & \text{ if }\gamma_{\ell+1}\leq 0,
                        \end{cases}
\qquad
S_{N_\ell}^+ \xi_\gamma =  \begin{cases}
                         \xi_{\gamma'''} & \text{ if }\gamma_{\ell+1}>0,\cr
                       \xi_{\gamma'} & \text{ if }\gamma_{\ell+1}\leq 0,
                        \end{cases}
\]
where
\[
\gamma'_i=\begin{cases}
           \gamma_i -1 &\text{ if }i=\ell+1,\cr
           \gamma_i & \text{ otherwise},
          \end{cases}
\qquad
 \gamma''_i=\begin{cases}
                 \gamma_i-1 & \text{ if } \ell\leq i \leq \ell+2,\cr
                \gamma_i & \text{ otherwise}.
                \end{cases}
\]
and
\[
 \gamma'''_i=\begin{cases}
                 \gamma_i+1 & \text{ if } i=\ell \text{ or } i=\ell+2,\cr
                \gamma_i-1 & \text{ if }i=\ell+1,\cr
                \gamma_i & \text{ otherwise}.
                \end{cases}
\]
Therefore
\[
 \left(S_{N_\ell}^+T^{+} + S_{N_0}^{-}T^{-}\right)\xi_\gamma
  = \begin{cases}
 q^k \xi_{\gamma'''}+(-1)^\ell q^{d_{\ell+1}}\xi_{\gamma'} & \text{ if }\gamma_{\ell+1}>0,\cr
 q^k \xi_{\gamma'}+(-1)^\ell q^{d_{\ell+1}}\xi_{\gamma''} & \text{ if }\gamma_{\ell+1}\leq 0.
    \end{cases}
\]
So if we now define 
\[
 T\xi_\gamma=
\begin{cases}
  (-1)^\ell q^{\sum_{i=1}^\ell \gamma_i}\xi_{\gamma'} & \text{ if } \gamma_{\ell+1}>0,\cr
q^{\sum_{i=1}^\ell \gamma_i}\xi_{\gamma'} & \text{ if } \gamma_{\ell+1} \leq 0,
\end{cases}
\]
then one gets from the above equation that 
$U\left(S_{N_\ell}^+T^{+} + S_{N_0}^{-}T^{-}  - T\right)U^*$
is in  $OP_{D_\ell}^{-\infty}\otimes\scrt_\ell^\infty$.
Thus it is enough to show that $UTU^*\in C^\infty(S_q^{2\ell+1})$.
Now note that
\[
 \eta(\gamma)-\eta(\gamma')=\begin{cases}
                                 \ell &\text{ if } \gamma_{\ell+1}>0,\cr
                              0 & \text{ if } \gamma_{\ell+1} \leq 0.
                                \end{cases}
\]
Therefore it follows that 
$UTU^* e_\gamma=  q^{\sum_{i=1}^\ell \gamma_i}e_{\gamma'}$,
i.e.\ $UTU^*= Y_{\ell+1,q}^*\otimes I$.
Thus we get the required result.
\qed

Putting together lemmas \ref{lem:xj}, \ref{lem:zjq2} and \ref{lem:proofpart3},
we get proposition~\ref{pro:zjq}.

\vspace{1ex}

\noindent
\begin{footnotesize}\textbf{Acknowledgement}: 
We would like to thank Partha Sarathi Chakraborty
for making us work on the problem.

In addition, SS would like to  thank his advisor Partha Sarathi Chakraborty
for  the support that he has given him.
He would also like to thank ISI Delhi Centre for its hospitality during his stay
there. 
\end{footnotesize}

%%%%%%%%%%%%%%%%%%%%%%%%%%%%%%%%%%%%%%%%%%%%%%%%%%%%
%%%%%%%%%%%%%%%%%%%%%%%%%%%%%%%%%%%%%%%%%%%%%%%%%%%%

%%%%%%%%%%%%%%%%%%%%%%%%%%%%%%%%%%%%%%%%%%%%%%%%%%%%

% \bibliographystyle{plain}
% \bibliography{master.bib}

%%%%%%%%%%%%%%%%%%%%%%%%%%%%%%%%%%%%%%%%%%%%%%%%%%%%%%%%%
\noindent{\sc Arupkumar Pal} (\texttt{arup@isid.ac.in})\\
         {\footnotesize Indian Statistical
Institute, 7, SJSS Marg, New Delhi--110\,016, INDIA}\\[1ex]
{\sc S. Sundar}
(\texttt{ssundar@imsc.res.in})\\
         {\footnotesize  Institute of Mathematical Sciences, 
CIT Campus, Chennai--600\,113, INDIA}

%%%%%%%%%%%%%%%%%%%%%%%%%%%%%%%%%%%%%%%%%%%%%%%%%%%%%%%%%%%

\end{document}